\documentclass[12pt]{article}

\usepackage{mathrsfs}
\usepackage{amsmath}
\usepackage{amsfonts}
\usepackage{amssymb}
\usepackage[all,cmtip]{xy}
\usepackage{hyperref}
\usepackage{enumerate}  
\usepackage{ulem}
\usepackage{eufrak}
\usepackage[all,cmtip]{xy}
\usepackage[toc,page]{appendix}
\usepackage{amsmath, amsthm, amsfonts, amssymb}
\usepackage{graphicx}
\usepackage{subcaption}
\usepackage{tikz-cd}

\def\<{\langle}
\def\>{\rangle}

\def\g{\gamma}

\def\-{\overline}

\def\h{\hbox}

\def\endpf{\hbox{\vrule height1.5ex width.5em}}

\def\endpf{\hbox{\vrule height1.5ex width.5em}}

\def\-{\overline}

\def\h{\hbox}

\def\endpf{\hbox{\vrule height1.5ex width.5em}}

\def\endpf{\hbox{\vrule height1.5ex width.5em}}

\def\g{\gamma}
\def\-{\overline}

\DeclareMathOperator{\Ima}{Im}
\newtheorem{fact}{Fact}[section]
\newtheorem{question}{Question}[section]
\newtheorem{dt}{Definition-Theorem}[section]
\newtheorem{theorem}{Theorem}[section]
\newtheorem{lemma}[theorem]{Lemma}
\newtheorem{corollary}[theorem]{Corollary}
\newtheorem{proposition}[theorem]{Proposition}

\newtheorem*{theorem*}{Theorem}
\date{}
\theoremstyle{definition}
\newtheorem{definition}[theorem]{Definition}
\newtheorem{example}[theorem]{Example}
\newtheorem{remark}[theorem]{Remark}

\begin{document}

\title{\bf Construct holomorphic invariants in \v Cech cohomology by a combinatorial formula}

\author{Hanlong Fang}

\maketitle

\begin{abstract}
Let $E$ be a holomorphic vector bundle over a complex manifold $X$. We construct some holomorphi invariants, which is called the $T$ invariants of $E$, in certain \v Cech cohomology groups by a combinatorial formula.  Then, we show that  when $X$ is a compact complex manifold the $T$ invariants of $E$ coincide with the power sums of the Chern roots of $E$, up to certain normalized factors. Next, we refine the first $T$ invariants for all holomorphic vector bundles and the $k$-th $T$ invariants for all holomoprhic vector bundles with a full flag structure.  At last, we generalize the notion of the $T$ invariants and the refined $T$ invariants to the locally free sheaves of schemes over a general field.
\end{abstract}
\tableofcontents
\bigskip\bigskip



\section{Introduction}
Mathematicians have been interested in investigating  the  invariants of fiber bundles  for a long time. According to E. Cartan, S. Chern, A. Weil and many other great geometers' philosophy, the properties of an underlying space can be better understood in the geometry of its fiber bundles.  By the works of Stiefel, Whitney, Pontrjagin, Steenrod, Chern, Weil, etc., to name a few, a great number of topological invariants of fiber bundles were discovered which also have been widely used in all fields of mathematics since then.  Nowadays, even in physics, it becomes more and more important to investigate and apply those invariants. For instance, Chern numbers and Chern-Simons forms are fundamental tools in studying the topological order of the quantum Hall effect.  

Among these invariants,  the Chern classes of vector bundles  play an extremely important role. In Chern's original work \cite{C}, he constructed the characteristic classes by integrating the invariant polynomials of a curvature form.  This approach has a far reaching influence on differential geometry and complex geometry, for it relates the topological properties of a manifold to its geometric properties by simple differentiation and integration.  People also extend the concept of the Chern classes largely in  various contexts.  For instance, in the famous Grothendieck-Riemann-Roch theorem, Chern classes are defined for the elements in the Grothendieck group, and the classes are in the Chow rings.

Because of its importance, mathematicians have been constructing the Chern classess of vector bundles by various methods. Since a vector bundle is  determined by a system of transition functions, one of the natural approaches would be to construct the Chere classes  by  certain formulas  only involving transition functions. In 1973, Bott \cite{B} initiated the study of constructing Chern classes of topological vector bundles by transition functions.  Later on, Sharygin \cite{Sh} prove that the Chern classes of smooth principle  $GL_n$-bundle can be represented by a \v Cech cocycle in the \v Cech-de Rham complex.  He also computed some interesting examples in certain lower dimensional cases.

In this paper, we will continue the study of \cite{B} and \cite{Sh} but focus more on its holomoprhic aspect instead of the topological aspect. To be more precise, we will  represent the Chern characters of the holomorphic vector bundles over complex manifolds or the locally free sheaves over schemes as \v Cech cocycles. Similar to \cite{B} and \cite{Sh}, our construction involves transition functions only. However, there are also some distinct features among our approach and theirs.  Firstly, by restricting our attention to the holomophic category, we produce \v Cech cocycles in the usual \v Cech cohomology group instead of the \v Cech-de Rham complex. Secondly, since our  proof of the invariance of the \v Cech cocycles is purely algebraic, namely, no curvature, connection or integration involved, it can be easily generalized  for the locally free sheaves on schemes over general fields. Thirdly, our formula is explicit and hence computable.

Moreover, by an observation motivation by \cite{F}, we construct some new holomorphic invariants refining the Chern classes, which is called the refined $T$ invariants or $\mathcal Q$-flat classes in the lowest degree case.
\medskip


We start with our main results by the following  key formula, which constructs invariant \v Cech cocycles by a combinatorial formula. 
\begin{theorem}[The combinatorial formula]\label{tcf}
Let $X$ be a complex manifold.  Denote by $\Omega^k$ the sheaf of germs of holomorphic differential $k$-forms. Let $E$ be a holomorphic vector bundle over $X$. Let $\mathcal U:=\{U_i\}_{i\in I}$ be an open cover of $X$, where $I$ is an ordered set, and $g:=\{g_{ij}\}$ be a system of transition functions associated with a certain trivialization of $E$ with respect to $\mathcal U$.  For every positive integer $k$ and indices $i_1,i_2,\cdots, i_{k+1}\in I$, define $t_{i_1\cdots i_{k+1}}\in\Gamma(U_{i_1\cdots i_{k+1}},\Omega^k)$ by
	\begin{equation}\label{t11}
	\begin{split}
t_{i_1\cdots i_{k+1}}:=\sum_{\sigma\in S_{k+1}}\frac{sgn(\sigma)}{(k+1)!}&\cdot tr\big(g^{-1}_{i_{\sigma(1)}i_{\sigma(k+1)}}dg_{i_{\sigma(1)}i_{\sigma(k+1)}}g^{-1}_{i_{\sigma(2)}i_{\sigma(k+1)}}dg_{i_{\sigma(2)}i_{\sigma(k+1)}}\cdot\\
	&\cdot g^{-1}_{i_{\sigma(3)}}i_{i_{\sigma(k+1)}}dg_{i_{\sigma(3)}i_{\sigma(k+1)}}\cdots g^{-1}_{i_{\sigma(k)}i_{\sigma(k+1)}}dg_{i_{\sigma(k)}i_{\sigma(k+1)}}\big);
	\end{split}
	\end{equation}
	define a \v Cech $k$-cochain $\widehat f_{k,\,\mathcal U}(E)$ by
	\begin{equation}\label{t123}
	\widehat f_{k,\,\mathcal U}(E):=\bigoplus\limits_{ i_1<\cdots<i_{k+1}}t_{i_1\cdots i_{k+1}}\in\bigoplus\limits_{i_1<\cdots<i_{k+1}}\Gamma(U_{i_1\cdots i_{k+1}},\Omega^k).
	\end{equation}
	Here $S_{k+1}$ is the permutation group of $\{1,2,\cdots,k+1\}$; $sgn(\sigma)=1$ for $\sigma$ an even permutation and $sgn(\sigma)=-1$ for $\sigma$ an odd permutation; $tr$ is the trace operator of differential form-valued matrices, where the product between differential forms is understood as the wedge product.
	
	Then $\widehat f_{k,\,\mathcal U}(E)$ is a \v Cech $k$-cocycle and hence defines an element $f_{k,\,\mathcal U}(E)\in\check H^k(\mathcal U,\Omega^k)$; moreover, $f_{k,\,\mathcal U}(E)$ is independent of the choice of $g$.  
\end{theorem}

Since the invariant \v Cech cocycles constructed in Theorem \ref{tcf} is compatible with the natural restriction $I_{\mathcal U\mathcal V}:\check H^k(\mathcal U,\Omega^k)\rightarrow \check H^k(\mathcal V,\Omega^k)$ where $\mathcal V$ is a refinement of $\mathcal U$ (see \S 2 for the details), we have the following theorem.

\begin{theorem}\label{cfd}
	Denote by  $f_k(E)$ the image of $f_{k,\,\mathcal U}(E)$ in the \v Cech cohomology group $\check H^k(X,\Omega^k)$ under the canonical homomorphism.
Then $f_k(E)$ is independent of the choice of $\mathcal U$. We call  $f_k(E)$ the $k$-th $T$ invariant.
\end{theorem}

To compute the $T$ invariants by sheaf cohomology, we introduce the following definition for the purpose of comparing the sheaf cohomology and the \v Cech cohomology of sheaf $d\Omega^*$.
\begin{definition}
	An open cover $\mathcal U:=\{U_i\}_{i\in I}$ is called a Stein cover if each $U_i$ is Stein, and is called a  good cover if each $U_i$ is Stein and each  nonempty intersection $U_{i_1\cdots i_p}$ is contractible. 
\end{definition}

We then have the following corollary, which identifies the \v Cech cohomology with the sheaf cohomology when $X$ is compact.

\begin{corollary}\label{comt}
Let $X$ be a compact complex manifold. Let $\pi:E\rightarrow X$ be a holomorphic vector bundle on $X$. 
Then $\check H^k(X,\Omega^k)\cong  H^k(X,\Omega^k)$. Moreover, there is a finite Stein cover $\mathcal U$ of $X$ such that $f_k(E)=\widehat f_{k,\,\mathcal U}(E)$ under the natural isomorphisms $\check H^k({\mathcal U},\Omega^k)\cong \check H^k(X,\Omega^k)\cong  H^k(X,\Omega^k)$, where $\widehat f_{k,\,\mathcal U}(E)$ is the \v Cech $k$-cocycle defined by formulas (\ref{t11}) and (\ref{t123}).	
\end{corollary}
\begin{remark}\label{rcomt}
When $X$ is a compact complex manifold,  we naturally view $f_k(E)$ as an element of $H^k(X,\Omega^k)$  by a slight abuse of notation.
\end{remark}

The following theorem shows that the invariant cocycles constructed in Theorem \ref{tcf} coincide with the power sums of Chern roots when $X$ is a compact complex manifold. 
\begin{theorem}\label{compare}
Let $X$ be a compact complex manifold and $E$ be a holomorphic vector bundle over $X$. Let ch$(E)$ be the Chern character of the holomophic vector bundle $E$ in the ring $\oplus_{i=0}^{\infty} H^k(X,\Omega^k)$ (see \S 4). Then the following equality holds in the ring $\oplus_{i=0}^{\infty} H^k(X,\Omega^k)$:
	\begin{equation}
	{\rm ch}(E)=\sum_{k=0}^{\infty}\frac{1}{k!\cdot(2\pi\sqrt{-1})^{k}}\cdot f_{k}(E),
	\end{equation}
    where $f_0(E)=1$ by convention.
\end{theorem}
\begin{remark}
The  Chern classes can be represented as symmetric polynomials of $T$ invariants. In particular,  the first Chern class in the Dolbeault cohomology is the first $T$ invariant, up to a certain normalized factor.  
\end{remark}

Next,  we derive a new holomorphic invariant which refines the first $T$ invariant by considering the \v Cech cohomology group $\check H^1(X,d\Omega^0)$ instead of $\check H^1(X,\Omega^1)$ (see \S 2 for the definitions).
\begin{theorem}\label{fla}
Let $X$ be a complex manifold.  Let $\pi:E\rightarrow X$ be a holomorphic vector bundle on $X$.   Let $\mathcal U:=\{U_i\}_{i\in I}$ be an open cover of $X$ and $g:=\{g_{ij}\}$ be a system of transition functions associated with a certain trivialization of $E$ with respect to $\mathcal U$.  For indices $i_1,i_2\in I$, define an element $t_{i_1i_{2}}\in\Gamma(U_{i_1i_{2}},\Omega^k)$ by
	\begin{equation}\label{t12}
	\begin{split}
	t_{i_1 i_{2}}:=\sum_{\sigma\in S_{2}}&\frac{sgn(\sigma)}{2!}\cdot tr\big(g^{-1}_{i_{\sigma(1)}i_{\sigma(2)}}dg_{i_{\sigma(1)}i_{\sigma(2)}}\big)=tr\big(g^{-1}_{i_1i_2}dg_{i_1i_2}\big);
	\end{split}
	\end{equation}
	define  a \v Cech $1$-cochain $\widehat f^r_{1,\,\mathcal U}(E)$ by 
	\begin{equation}\label{f1u}
	\widehat f^r_{1,\,\mathcal U}(E):=\bigoplus\limits_{i_1<i_{2}}t_{i_1i_{2}}\in\bigoplus\limits_{ i_1<i_{2}}\Gamma(U_{i_1i_{2}},d\Omega^0).
	\end{equation}
	Then $\widehat f^r_{1,\,\mathcal U}(E)$ is a \v Cech $1$-cocycle and hence defines an element $f^r_{1,\,\mathcal U}(E)\in\check H^1(\mathcal U,d\Omega^0)$; moreover,  $f^r_{1,\,\mathcal U}(E)$ is independent of the choice of $g$.  Denote by $f^r_1(E)$ the image of $f^r_{1,\,\mathcal U}(E)$ in the \v Cech cohomology group $\check H^1(X,d\Omega^0)$ under the canonical homomorphism; then $f^r_1(E)$ is independent of the choice of $\mathcal U$. We call  $f^r_1(E)$ the the refined first $T$ invariant.	
\end{theorem}



When $X$ is compact, we can compute the refined first $T$ invariant by the following corollary.	

\begin{corollary}\label{comqf}
Let  $X$ be a compact complex manifold. Let $\pi:E\rightarrow X$ be a holomorphic vector bundle on $X$.  Then $\check H^1(X,d\Omega^0)\cong H^1(X,d\Omega^0)$. Moreover, there is a finite good cover $\mathcal U$ of $X$ such that $f_1^r(E)=\widehat f^r_{1,\,\mathcal U}(E)$ under the natural isomorphisms $\check H^1(\mathcal U,d\Omega^0)\cong\check H^1(X,d\Omega^0)\cong H^1(X,d\Omega^0)$, where  $\widehat f^r_{1,\,\mathcal U}(E)$ is the \v Cech $1$-cocycle defined by formulas (\ref{t12}) and (\ref{f1u}).	
\end{corollary} 
\begin{remark}
 When $X$ is a compact complex manifold,  we naturally view $f_1^r(E)$ as an element of $H^1(X,d\Omega^0)$ by  a slight abuse of notation. Sometimes, we even simply denote $f_1^r(E)$ by  $f_1(E)$ when there is no confusion.
\end{remark}

Let $E$ be  a holomorphic line bundle.  We say $E$ is flat if its transition functions can be taken as constant functions with repsect to a certain trivialization; we say $E$ is $\mathcal Q$-flat if there is a positive integer $m$ such that $mE$ is flat.

We have the following theorem which explains the geometric meaning of the refined first $T$ invariants.

\begin{theorem}\label{flat}
	Let $X$ be a compact complex manifold.  Let $\pi:E\rightarrow X$ be a holomorphic line bundle over $X$. Then $E$ is $\mathcal Q$-flat if and only if the refined first $T$ invariant of $E$ is trivial in $H^1(X,d\Omega^0)$.
\end{theorem}

Notice that in many situations the refined first $T$ invariants coincide with the first $T$ invariants, for instance, when $X$ is a compact K\"ahler manifold. The following result gives a criterion for determining whether a manifold has a line bundle whose refined first $T$ invariant is strictly finer than its first $T$ invariant. (See  \S 4.3 or \cite{Fr} for the notion of the Fr\"olicher spectral sequences.) 

\begin{theorem}\label{criterion}
 Let $X$ be a compact complex manifold. There is a line bundle of $X$ with trivial first $T$ invariant but non-trivial refined first $T$ invariant if only if in the Fr\"olicher spectral sequence of $X$,
	\begin{equation}
	E_2^{0,1}\neq E_3^{0,1}\,\,.
	\end{equation}	
\end{theorem}
\begin{remark}
Together with an example by \cite{RB} (see Example \ref{exam}), there is a compact complex manifold with a holomorphic line bundle such that the refined first $T$ invariant of the line bundle is non-trivial but its first $T$ invariant is trivial.
\end{remark}

We say a holomorphic vector bundle $E$ over $X$ of rank $M$ has a full flag structure if $E$ has a filtration by holomoprhic subbundles
\begin{equation}
E=E_M\supset E_{M-1}\supset\cdots\supset E_2\supset E_1\supset E_0=0
\end{equation}
with line bundle quotients $L_i=E_i/E_{i-1}$. Denote the full flag structure by $\mathcal F$.  Notice that for a Stein open cover $\mathcal U$ of $X$ one can choose a trivialization of $E$ with respect to $\mathcal F$ and $\mathcal U$ such that the tansition functions $g_{ij}^{\prime}s$ are upper triangular matrices. 
We call such system of transition functions a system of transition functions with respect to $\mathcal F$ and $\mathcal U$. 

Noticing that formula (\ref{cf}) does not give a \v Cech cocycle with coefficients in the sheaf $d\Omega^k$ for $k\geq 2$, it is difficult to define the refined higher $T$ invariants for generic holomorphic vector bundles. For a holomorphic vector bundle with a full flag structure, we define its refined higher $T$ invariants  as follows.

\begin{theorem}\label{f}
Let $X$ be a complex manifold.  Let $\pi:E\rightarrow X$ be a holomorphic vector bundle over $X$ with a full flag structure $\mathcal F$. Let $\mathcal U:=\{U_i\}_{i\in I}$ be an open cover of $X$ and $g:=\{g_{ij}\}$ be a system of transition functions associated with a certain trivialization of $E$ with respect to $\mathcal F$ and $\mathcal U$.  Then, there is a well-defined \v Cech $k$-cocycle $\widehat f^r_{k,\,\mathcal U}(E,\mathcal F)$,
	\begin{equation}\label{tr1}
	\widehat f^r_{k,\,\mathcal U}(E,\mathcal F):=\bigoplus\limits_{i_1<\cdots<i_{k+1}}t_{i_1\cdots i_{k+1}}\in\bigoplus\limits_{i_1<\cdots<i_{k+1}}\Gamma(U_{i_1\cdots i_{k+1}},d\Omega^{k-1}),
	\end{equation}
	where
	\begin{equation}\label{tr2}
	\begin{split}
	t_{i_1\cdots i_{k+1}}:=\sum_{\sigma\in S_{k+1}}&\frac{sgn(\sigma)}{(k+1)!}\cdot tr\big(g^{-1}_{i_{\sigma(1)}i_{\sigma(k+1)}}dg_{i_{\sigma(1)}i_{\sigma(k+1)}}g^{-1}_{i_{\sigma(2)}i_{\sigma(k+1)}}dg_{i_{\sigma(2)}i_{\sigma(k+1)}}\cdot\\
	&\cdot g^{-1}_{i_{\sigma(3)}i_{\sigma(k+1)}}dg_{i_{\sigma(3)}i_{\sigma(k+1)}}\cdots g^{-1}_{i_{\sigma(k)}i_{\sigma(k+1)}}dg_{i_{\sigma(k)}i_{\sigma(k+1)}}\big).
	\end{split}
	\end{equation}
	Denote  by $f^r_k(E,\mathcal F)$ the image of $\widehat f^r_{k,\,\mathcal U}(E,\mathcal F)$ in $\check H^k(X,d\Omega^{k-1})$ under the canonical homomorphism.  Then $f_k^r(E,\mathcal F)$ is independent of the choice of $\mathcal U$. We call  $f_k^r(E,\mathcal F)$ the refined $k$-th $T$ invariant.

\end{theorem}

For compact complex manifolds, we have the following corollary. 
\begin{corollary}\label{comf}
	Let  $X$ be a compact complex manifold. Let $\pi:E\rightarrow X$ be a holomorphic vector bundle on $X$.  Then $\check H^k(X,d\Omega^{k-1})\cong H^k(X,d\Omega^{k-1})$. Moreover, there is a finite good cover $\mathcal U$ of $X$ such that $f_k^r(E,\mathcal F)=\widehat f^r_{k,\,\mathcal U}(E,\mathcal F)$ under the natural isomorphisms $\check H^k(\mathcal U,d\Omega^{k-1})\cong\check H^k(X,d\Omega^{k-1})\cong H^k(X,d\Omega^{k-1})$, where  $\widehat f^r_{k,\,\mathcal U}(E,\mathcal F)$ is the \v Cech $k$-cocycle defined by formulas (\ref{tr1}) and (\ref{tr2}).
\end{corollary}

Finally, we generalize the notion of $T$ invariants for the locally free sheaves on schemes over a field $K$.  First, let us recall some notation in \cite{Har} (see Exercises II.5.16 and II.5.18 and \S 2.8 therein for reference). Let $X$ be a scheme over a field $K$. Denote by $\Omega_{X/{\rm Spec\,}K}$ the sheaf of K\"ahler differential of $X$ over Spec $K$. Denote by $\Omega^k$ the $k$-th exterior power sheaf of $\Omega_{X/{\rm Spec\,}K}$ for $k\geq 1$.  Then similar to formulas $(\ref{tran1})$, $(\ref{tran2})$ and $(\ref{tran3})$ in \S 2, we can define the notion of a system of transition functions associated with a certain trivialization with respect to a open cover $\mathcal U$ in the Zariski topology.

Under the above notation, we have the following theorem.
\begin{theorem}\label{sch}
Let $X$ be a scheme over a field $K$. Suppose that $k\geq 1$ and $(k+2)!$ is not divisible by the characteristic of $K$. Let $E$ be a locally free sheaf over $X$. Let $\mathcal U:=\{U_i\}_{i\in I}$ be an open cover of $X$ (in the Zariski topology) and $g:=\{g_{ij}\}$ be a system of transition functions associated with a certain trivialization of $E$ with respect to $\mathcal U$. 
For any indices $i_1,i_2,\cdots, i_{k+1}\in I$, define $t_{i_1\cdots i_{k+1}}\in\Gamma(U_{i_1\cdots i_{k+1}},\Omega^k)$ by
 \begin{equation}\label{scheme}
 \begin{split}
 t_{i_1\cdots i_{k+1}}:=\sum_{\sigma\in S_{k+1}}\frac{sgn(\sigma)}{(k+1)!}&\cdot tr\big(g^{-1}_{i_{\sigma(1)}i_{\sigma(k+1)}}dg_{i_{\sigma(1)}i_{\sigma(k+1)}}g^{-1}_{i_{\sigma(2)}i_{\sigma(k+1)}}dg_{i_{\sigma(2)}i_{\sigma(k+1)}}\cdot\\
 &\cdot g^{-1}_{i_{\sigma(3)}i_{\sigma(k+1)}}dg_{i_{\sigma(3)}i_{\sigma(k+1)}}\cdots g^{-1}_{i_{\sigma(k)}i_{\sigma(k+1)}}dg_{i_{\sigma(k)}i_{\sigma(k+1)}}\big);
 \end{split}
 \end{equation}
 define a \v Cech $k$-cochain $\widehat f_{k,\,\mathcal U}(E)$ by
 \begin{equation}
 \widehat f_{k,\,\mathcal U}(E):=\bigoplus\limits_{ i_1<\cdots<i_{k+1}}t_{i_1\cdots i_{k+1}}\in\bigoplus\limits_{i_1<\cdots<i_{k+1}}\Gamma(U_{i_1\cdots i_{k+1}},\Omega^k).
 \end{equation}
	Then $\widehat f_{k,\,\mathcal U}(E)$ is a \v Cech $k$-cocycle and hence defines an element $f_{k,\,\mathcal U}(E)\in\check H^k(\mathcal U,\Omega^k)$; moreover, $f_{k,\,\mathcal U}(E)$ is independent of the choice of $g$.  
	Denote by  $f_k(E)$ the image of $f_{k,\,\mathcal U}(E)$ in the \v Cech cohomology group $\check H^k(X,\Omega^k)$ under the canonical homomorphism.
	Then $f_k(E)$ is independent of the choice of $\mathcal U$. We call  $f_k(E)$ the $k$-th $T$ invariant in $\check H^k(X,\Omega^k)$.
\end{theorem}

Since $\check H^k(X,\Omega^k)\cong H^k(X,\Omega^k)$ for separated schemes,  we can compute the $T$ invariants in sheaf cohomology for separated schemes.  Moreover, we can  define the cohomological Chern character  for separated schemes as follows.
\begin{definition}\label{cdo}
	Let $X$ be a separeted scheme over a field $K$. Let $E$ be a local free sheaf over $X$. Suppose that the characteristic of $K$ is larege enough and that $X$ has a finite affine open cover. Then the cohomological Chern character is defined by
	\begin{equation}
	{\rm ch}_{coh}(E):=\sum_{k=0}^{\infty}\frac{1}{k!\cdot(2\pi\sqrt{-1})^{k}}\cdot f_{k}(E).
	\end{equation}
	Here we make the convention that when $f_k(E)=0\in H^k(X,\Omega^k)$, $\frac{1}{k!\cdot(2\pi\sqrt{-1})^{k}}\cdot f_{k}(E)=0$, even if $k$ is divisible by the characteristic of $K$.
\end{definition}
\begin{remark}
	Since $H^k(X,\Omega^k)=0$ for $k$ large enough, the above infinite sum is a finite sum by the convention.
\end{remark}
The refined first $T$ invariants and the refined $T$ invariants for vector bundles will full flag structures can be carried over to schemes over a general field as well (see Definition-Theorems \ref{adholt} and \ref{drholt}). 

Moreover, the $T$ invariants and the refined $T$ invariants can be generalized to relative schemes as well.  We leave this to the future study.
\medskip

We now briefly describe the organization of the paper and the basic ideas for the proof of Theorems \ref{tcf}-\ref{sch}.  The most technical part of the paper is the proof of Theorem \ref{tcf}.  But as long as one has a correct guess of the fomula, the remaining  is a routine algebraic computation combining the  permutation technique, the properties of the trace operator and the use of the cocycle condition of the transition functions. To compare the $T$ invariants and the Chern characters in the Dolbeault cohomology, we use the Hirsch lemma and the splitting principle. Noticing that the first $T$ invariants are in smaller sheaf $d\Omega^0$,  we can lift the first $T$ invariant from $H^1(X,\Omega^1)$ to $H^1(X,d\Omega^0)$. To refine the higher $T$ invariants,  we introduce the notion of full flag structures for holomorphic vector bundles. Since the combinotorical formula is proved by an algebraic computation, we can generalize the previous discussion easily to the locally free sheaves of schems over general field.

The organization of the paper is as follows:  In \S 2, we introduce some notations and basic facts of holomorphic vector bundles and \v Cech cohomology of various sheaves.  In \S 3,  we prove Theorem \ref{tcf} which says that the \v Cech cocycles defined by the combinatorial formula are invariant. In \S 4.1, we apply the combinatorial formula to holomorphic vector bundles and proved Theorem \ref{cfd}, Corollary \ref{comt} and Theorem \ref{compare}. In \S 4.2, we lift the first $T$ invariant from $\check H^1(X,\Omega^1)$ to $\check H^1(X,d\Omega^0)$ and prove Theorems \ref{fla}, \ref{comqf} and \ref{flat}. In \S 4.3, we prove Theorem \ref{criterion}, the strictly finer criterion. In \S 5.1, we discuss the basic properties of the cohomology ring $\Phi_X$ in which the refined $T$ invariants live and compute the ring $\Phi_X$ for complex projective spaces. In \S 5.2, we establish Theorems \ref{f} and \ref{comf} which refine the higher $T$ invariants for holomorphic vector bundles with a full flag structure. In \S 6.1, we derive Theorem \ref{sch} and define the $T$ invariants, refined $T$ invariants and the cohomological Chern characters for the locally free sheaves of schemes over a general field.  In \S 6.2, we collect an example for the \v Cech cocycle representation of the first Stiefel-Whitney class of a real vector bundle. At the end of each section, we also raise some questions for later investigation.
\medskip

{\bf Acknowledgement :} The author  appreciates greatly his advisor Prof. Xiaojun Huang for emphasizing repeatedly the importance of constructing holomorphic invariants. He would like to thank Xu Yang for helpful discussion and warm encouragement.  Also, the author would like to thank Prof. Song-Ying Li, during a hiking with whom he got the key idea of the paper.  Finally, the author would like to thank Prof. Xuguang Lu for teaching him the permutation technique when he was an undergraduate student in Tsinghua University.
\section{Preliminary}

In this section, we will  recall some notations and basic facts of holomoprhic vector bundles and \v Cech cohomoloy which will be used throughout the remaining of the paper.

Suppose that $X$ is a complex manifold and $E$ is a holomorphic vector bundle over $X$ of rank $M$. Then, there is an open cover $\mathcal U:=\{U_i\}_{i\in I}$ of $X$ and a holomorphic trivialization $(E,\{U_i\},\{\phi_i\})$ of $E$, where $I$ is an ordered set, as follows
\begin{equation}\label{tran1}
\phi_i:E\big|_{U_i}\xrightarrow{\cong}U_i\times\mathbb C^M,\,\,i\in I. 
\end{equation}
(For example, we can choose $\mathcal U$ to be any Stein open cover of $X$, for holomorphic vector bundles are holomorphically trivial over Stein sets.)
The maps 
\begin{equation}\label{tran2}
\begin{split}
\phi_i\circ\phi_j^{-1}:(U_i\cap U_j)\times\mathbb C^M&\rightarrow (U_i\cap U_j)\times\mathbb C^M\\
\end{split}
\end{equation}
are vector-space automorphisms of $\mathbb C^M$ in each fiber and hence give rise to maps 
\begin{equation}\label{tran3}
\begin{split}
g_{ij}:(U_i\cap U_j)&\rightarrow GL(M,\mathbb C)\\
g_{ij}(z)=&\phi_i\circ\phi_j^{-1}\big|_{z\times\mathbb C}
\end{split}.
\end{equation}
We call such $g:=\{g_{ij}\}$ the system of (matrix-valued) transition functions associated with trivialization $(E,\{U_i\},\{\phi_i\})$.
\begin{remark}
In the following,  we usually refer to a system of transition functions without mentioning the trivialization it associated with, for we will only use the information of the transition functions in this paper. 
\end{remark}

Notice that there is a non-abelian \v Cech cohomology interpretation of the holomorphic vector bundles as follows  (see \cite{G} or \cite{Br} for more details).  Since $g_{ij}$ is a $GL(M,\mathbb C)$-valued holomorphic function and $g_{ij}g_{jk}g_{ki}\equiv 1$ on $U_{ijk}:=U_i\cap U_j\cap U_k$, one can associate $g$ with the following \v Cech $1$-cocycle
\begin{equation}
\widehat g:=\bigoplus\limits_{i_1<i_2}g_{i_1i_2}\big|_{U_{i_1i_2}}\in\bigoplus\limits_{i_1<i_2}\Gamma_{\rm hol} (U_{i_1i_2},GL(M,\mathbb C)),
\end{equation}
where $\Gamma_{\rm hol} (U_{i_1i_2},GL(M,\mathbb C))$ is the (nonabelian) group consisting of $GL(M,\mathbb C)$-valued holomorphic functions over $U_{i_1i_2}$. 
If  $(E,\{U_i\},\{\widetilde\phi_i\})$ is another trivialization of $E$ with respect to the same cover, there exists a \v Cech $0$-cochain
\begin{equation}
h:=\bigoplus\limits_{i_1}h_{i_1}\in\bigoplus\limits_{i_1}\Gamma_{\rm hol} (U_{i_1},GL(M,\mathbb C)),
\end{equation}
such that $\widetilde g_{ij}=h_i^{-1}g_{ij}h_j$ for $i,j\in I$. Denote by $\check H_{{\rm hol}}^1(\mathcal U,GL(M,\mathbb C))$ the quotient group under the above relation.  It is easy to verify that each holomorphic line bundle determines a well-defined element in the \v Cech cohomology group 
\begin{equation}
\check H^1_{\rm hol}(X,GL(M,\mathbb C)):=\lim_{\substack{\longrightarrow\\\mathcal U}}\check H_{\rm hol}^1(\mathcal U,GL(M,\mathbb C)),
\end{equation} 
where the direct limit runs over all the open covers.

Denote by $\Omega^k$ be the sheaf of germs of holomorphic $k$-forms for $k\geq 0$; by convention, $\Omega^0=\mathcal O_X$.  Recall that $ \check H^{k}(\mathcal U,\Omega^k)$ is the $k$-th \v Cech cohomology group of the sheaf $\Omega^{k}$ with respect to $\mathcal U$. Let $I_{\mathcal U\mathcal V}$ be the natural  induced map from  $\check H^k(\mathcal U,\Omega^k)$ to $\check H^k(\mathcal V,\Omega^k)$ for open covers $\mathcal U$ and $\mathcal V$ such that $\mathcal V$ is a refinement of $\mathcal U$; denote by  $\{\check H^{k}(\mathcal U,\Omega^k),I_{\mathcal U\mathcal V}\}$ the direct system indexed by the direct set formed by the open covers under refinement. $\check H^k(X,\Omega^k)$, the $k$-th \v Cech cohomology group of $\Omega^{k}$, is the direct limit of  $\{\check H^{k}(\mathcal U,\Omega^k),I_{\mathcal U\mathcal V}\}$,
\begin{equation}
\check H^k(X,\Omega^k):=\lim_{\substack{\longrightarrow\\ \mathcal U}} \check H^{k}(\mathcal U,\Omega^k);
\end{equation}
there is a canonical homomorphism for each open cover $\mathcal U$,
\begin{equation}
L_{\mathcal U}:\check H^{k}(\mathcal U,\Omega^k)\rightarrow \check H^k(X,\Omega^k).
\end{equation}
Also, there is a natural homomorphism from \v Cech cohomology to sheaf cohomology,
\begin{equation}\label{L}
L:\check H^k(X,\Omega^k)\rightarrow H^k(X,\Omega^k).
\end{equation}
Notice that if $\mathcal U$ is a Stein cover in the sense that each $U_i$ is Stein, $H^k(X,\Omega^k)\cong\check H^{k}(\mathcal U,\Omega^k)$; moreover, if each open cover has a Stein cover as its  refinement,  then $L$ is an isomorphism.  In particular, if $X$ is a compact complex manifold, $L$ is an isomorphism.

Following \cite{HA}, for integer $q\geq 0$, we denote by   $d\Omega^q$  the sheaf of germs of closed holomorphic $(q+1)$-forms.  
Recall that for each good cover $\mathcal U$ and each $q>0$, $p>0$ and $k\geq 0$, $H^q(U_{i_1\cdots i_p},d\Omega^k)=0$ (see \S 3.3 of \cite{F}), and hence
\begin{equation}
H^{q+1}(X,d\Omega^{q})\cong\check H^{q+1}({\mathcal U},d\Omega^{q})\,\,\,{\rm for }\,\,\,q\geq 0.
\end{equation}
Moreover, if each open cover has a good cover as its refinement, then the following natural map  from \v Cech cohomology to sheaf cohomology is an isomorphism,
\begin{equation}
L:\check H^{q+1}(X,d\Omega^{q})\xrightarrow{\cong} H^{q+1}(X,d\Omega^{q})\,\,{\rm for}\,\,q\geq 0;
\end{equation}
In particular,  $L$ is an isomorphism for compact complex manifolds (see the good cover lemma in Appendix I of \cite{F} for instance).
Moreover, the natural inclusion $J_{q+1}:d\Omega^q\rightarrow \Omega^{q+1}$ of $\mathbb C$-sheaves  induces a natural homomorphism between cohomology groups
\begin{equation}
j_{q+1}:H^{q+1}(X,d\Omega^{q})\rightarrow H^{q+1}(X,\Omega^{q+1})\,\,{\rm for }\,\,q\geq 0.
\end{equation}

\section{The combinatorial formula}
Let $X$ be a complex manifold. Let $E$ be a holomorphic vector bundle over $X$; denote the rank of $E$ by $M$. In this section, we will prove Theorem \ref{tcf} by the following lemmas. 

\begin{lemma}\label{cochain}
	Let $\mathcal U:=\{U_i\}_{i\in I}$ be an open cover of $X$, where $I$ is an ordered set, and $g:=\{g_{ij}\}$ be a system of transition functions of $E$ with respect to $\mathcal U$.  For each integer $k$ and $i_1,\cdots,i_{k+1}\in I$, define $t_{i_1\cdots i_{k+1}}\in\Gamma(U_{i_1\cdots i_{k+1}},\Omega^k)$ by
\begin{equation}\label{cf}
\begin{split}
\sum_{\sigma\in S_{k+1}}\frac{sgn(\sigma)}{(k+1)!}\cdot tr\big(&g^{-1}_{i_{\sigma(1)}i_{\sigma(k+1)}}dg_{i_{\sigma(1)}i_{\sigma(k+1)}}g^{-1}_{i_{\sigma(2)}i_{\sigma(k+1)}}dg_{i_{\sigma(2)}i_{\sigma(k+1)}}g^{-1}_{i_{\sigma(3)}i_{\sigma(k+1)}}dg_{i_{\sigma(3)}i_{\sigma(k+1)}}\\
&\cdots g^{-1}_{i_{\sigma(k)}i_{\sigma(k+1)}}dg_{i_{\sigma(k)}i_{\sigma(k+1)}}\big).
\end{split}
\end{equation}
Then $t_{i_1\cdots i_ji_{j+1}\cdots i_{k+1}}=-t_{i_1\cdots i_{j+1}i_j\cdots i_{k+1}}$. Hence,  $t_{i_1\cdots i_{k+1}}=0$ if there is a repeated index in $\{i_1,\cdots,i_{k+1}\}$.
\end{lemma}
\begin{remark} We use the following convention in Lemma \ref{cochain}. $dg_{ij}$  is the matrix derived by differentiating matrix $g_{ij}$ entry by entry.  The product $\cdots dg_{i_{\sigma(1)}i_{\sigma(k)}}g^{-1}_{i_{\sigma(2)}i_{\sigma(k)}}dg_{i_{\sigma(2)}i_{\sigma(k)}}\cdots$ is the matrix product  of matrices with differential form-valued entries. For example, if 
	\begin{equation}\label{ex}
	g=\begin{pmatrix}
	g_{11}&g_{12}\\g_{21}&g_{22}
	\end{pmatrix},\,\,h=\begin{pmatrix}
	h_{11}&h_{12}\\h_{21}&h_{22}
	\end{pmatrix},\,\,\,
	\end{equation}
	then 
	\begin{equation}\label{ex1}
	dg=\begin{pmatrix}
	dg_{11}&dg_{12}\\dg_{21}&dg_{22}
	\end{pmatrix},\,\,dh=\begin{pmatrix}
	dh_{11}&dh_{12}\\dh_{21}&dh_{22}
	\end{pmatrix},\,\,\,
	\end{equation}
	and
	\begin{equation}\label{ex2}
	dgdh=\begin{pmatrix}
	dg_{11}\wedge dh_{11}+dg_{12}\wedge dh_{21}&dg_{11}\wedge dh_{12}+dg_{12}\wedge dh_{22}\\dg_{21}\wedge dh_{11}+dg_{22}\wedge dh_{21}&dg_{21}\wedge dh_{12}+dg_{22}\wedge dh_{22}
	\end{pmatrix}.
	\end{equation}
\end{remark}
\noindent{\bf Proof of Lemma \ref{cochain} :} 
Denote the permutation $(j,j+1)\in S_{k+1}$ by $\tau$; namely, $\tau(j)=j+1$, $\tau(j+1)=j$ and $\tau(q)=q$ for each $q$, $1\leq q\leq k+2$ and $q\neq j,j+1$.

Recall that
\begin{equation}
\begin{split}
t_{i_1\cdots i_ji_{j+1}\cdots i_{k+1}}&=\sum_{\sigma\in S_{k+1}}\frac{sgn(\sigma)}{(k+1)!}\cdot tr\big(g^{-1}_{i_{\sigma(1)}i_{\sigma(k+1)}}dg_{i_{\sigma(1)}i_{\sigma(k+1)}}g^{-1}_{i_{\sigma(2)}i_{\sigma(k+1)}}dg_{i_{\sigma(2)}i_{\sigma(k+1)}}\\
&\cdot g^{-1}_{i_{\sigma(3)}i_{\sigma(k+1)}}dg_{i_{\sigma(3)}i_{\sigma(k+1)}}\cdots g^{-1}_{i_{\sigma(k)}i_{\sigma(k+1)}}dg_{i_{\sigma(k)}i_{\sigma(k+1)}}\big).
\end{split}
\end{equation}
Also, we have
\begin{equation}\label{tij}
\begin{split}
t_{i_1\cdots i_{j+1}i_j\cdots i_{k+1}}=\sum_{\sigma\in S_{k+1}}&\frac{sgn(\sigma)}{(k+1)!}\cdot tr\big(g^{-1}_{i_{(\tau\circ\sigma)(1)}i_{(\tau\circ\sigma)(k+1)}}dg_{i_{(\tau\circ\sigma)(1)}i_{(\tau\circ\sigma)(k+1)}}g^{-1}_{i_{(\tau\circ\sigma)(2)}i_{(\tau\circ\sigma)(k+1)}}dg_{i_{(\tau\circ\sigma)(2)}i_{(\tau\circ\sigma)(k+1)}}\\
&\cdot g^{-1}_{i_{(\tau\circ\sigma)(3)}i_{(\tau\circ\sigma)(k+1)}}dg_{i_{(\tau\circ\sigma)(3)}i_{(\tau\circ\sigma)(k+1)}}\cdots g^{-1}_{i_{(\tau\circ\sigma)(k)}i_{(\tau\circ\sigma)(k+1)}}dg_{i_{(\tau\circ\sigma)(k)}i_{(\tau\circ\sigma)(k+1)}}\big)\\
=sgn(\tau)\cdot\sum_{\sigma\in S_{k+1}}&\frac{sgn(\tau\circ\sigma)}{(k+1)!}\cdot tr\big(g^{-1}_{i_{(\tau\circ\sigma)(1)}i_{(\tau\circ\sigma)(k+1)}}dg_{i_{(\tau\circ\sigma)(1)}i_{(\tau\circ\sigma)(k+1)}}g^{-1}_{i_{(\tau\circ\sigma)(2)}i_{(\tau\circ\sigma)(k+1)}}dg_{i_{(\tau\circ\sigma)(2)}i_{(\tau\circ\sigma)(k+1)}}\\
\cdot& g^{-1}_{i_{(\tau\circ\sigma)(3)}i_{(\tau\circ\sigma)(k+1)}}dg_{i_{(\tau\circ\sigma)(3)}i_{(\tau\circ\sigma)(k+1)}}\cdots g^{-1}_{i_{(\tau\circ\sigma)(k)}i_{(\tau\circ\sigma)(k+1)}}dg_{i_{(\tau\circ\sigma)(k)}i_{(\tau\circ\sigma)(k+1)}}\big)\\
=sgn(\tau)\cdot\sum_{\widetilde\sigma\in S_{k+1}}&\frac{sgn(\widetilde\sigma)}{(k+1)!}\cdot tr\big(g^{-1}_{i_{\widetilde\sigma(1)}i_{\widetilde\sigma(k+1)}}dg_{i_{\widetilde\sigma(1)}i_{(\widetilde\sigma)(k+1)}}g^{-1}_{i_{\widetilde\sigma(2)}i_{\widetilde\sigma(k+1)}}dg_{i_{\widetilde\sigma(2)}i_{\widetilde\sigma(k+1)}}\\
\cdot& g^{-1}_{i_{\widetilde\sigma(3)}i_{\widetilde\sigma(k+1)}}dg_{i_{\widetilde\sigma(3)}i_{\widetilde\sigma(k+1)}}\cdots g^{-1}_{i_{\widetilde\sigma(k)}i_{\widetilde\sigma(k+1)}}dg_{i_{\widetilde\sigma(k)}i_{\widetilde\sigma(k+1)}}\big)\\
=-t_{i_1\cdots i_ji_{j+1}\cdots i_{k+1}}&.
\end{split}
\end{equation}

Clearly, $t_{i_1\cdots i_{k+1}}=0$ if there is a repeated index in $\{i_1,\cdots,i_{k+1}\}$ by equality (\ref{tij}). 

We complete the proof of Lemma \ref{cochain}.
\,\,\,\,$\endpf$
\medskip

We define a \v Cech $k$-cochain $\widehat f_k(E,g)$ by
\begin{equation}\label{gcf}
\widehat f_k(E,g):=\bigoplus\limits_{i_1<\cdots<i_{k+1}}t_{i_1\cdots i_{k+1}}\in\bigoplus\limits_{i_1<\cdots<i_{k+1}}\Gamma(U_{i_1\cdots i_{k+1}},\Omega^k).
\end{equation}
\begin{remark}\label{all}
	For convenience, we extend the components of the \v Cech $k$-cochain $\widehat f_k(E,g)$ in formula $(\ref{gcf})$ to all $(k+1)$-tuples of elements in $I$ as follows. For any elements $j_1,\cdots,j_{k+1}\in I$, we define $t_{j_1\cdots j_{k+1}}=0$, if there is a repeated index in the set $\{j_1,\cdots, j_{k+1}\}$;  we define $t_{j_1\cdots j_{k+1}}=sgn(\sigma)\cdot\alpha_{j_{\sigma(1)}\cdots j_{\sigma(k+1)}}$, if the indices are all distinct, where $\sigma\in S_{k+1}$ and $j_{\sigma(1)}<\cdots<j_{\sigma(k+1)}$. Then by Lemma \ref{cochain}, formula (\ref{cf}) is compatible with the extension of the components of $\widehat f_k(E,g)$ to all $(k+1)$-tuples of elements in $I$.
\end{remark}

Next, we will prove that $\widehat f_k(E,g)$ is closed.
\begin{lemma}\label{cocycle}
	The \v Cech $k$-cochain $\widehat f_k(E,g)$ defined by formula (\ref{gcf}) is a \v Cech $k$-cocycle.
\end{lemma}

\noindent{\bf Proof of Lemma \ref{cocycle} :} It suffices to prove that for any $i_1,\cdots,i_{k+2}\in I$ the following equality holds:
\begin{equation}\label{cocy}
\sum_{j=1}^{k+2}(-1)^{j+1}t_{i_1\cdots\widehat {i_j}\cdots i_{k+2}}=0,
\end{equation}
where $\widehat {i_j}$ is the usual notation for omitting the $j$-th index $i_j$.

Let $\sigma\in S_{k+1}$; for $j=1,\cdots,k+2$, define the $j$-th lift $\widehat\sigma^j\in S_{k+2}$ of $\sigma$ as follows,
\begin{equation}\label{lift}
\widehat\sigma^j(q)=\left\{\begin{matrix}
\sigma(q)&\,\,\,\,\,\,\,\,\,\,\,{\rm if}\,\,\,\,\,\,\,\,\,\,\,\sigma(q)< j\,\,\,\,{\rm and}\,\,\,\,1\leq q\leq k+1;\\\sigma(q)+1&\,\,\,\,\,\,\,\,\,\,\,{\rm if}\,\,\,\,\,\,\,\,\,\,\,\sigma(q)\geq j\,\,\,\,{\rm and}\,\,\,\,1\leq q\leq k+1;\\j&\,\,\,\,\,\,\,\,\,\,\,\,\,\,\,\,\,\,\,\,\,\,\,\,{\rm if}\,\,\,\,\,\,\,\,\,\,\,\,\,\,\,\,\,\,\,\,\,\,\,\,\,\,\,\,\,\,\,\,q=k+2.\,\,\,\,\,\,\,\,\,\,\,\,\,\,\,\,\,\,\,\,\,\,\,\,\,\,\,\,\,\,\,\,\,\,\,\,\,\,\,\,\,
\end{matrix}\right  .
\end{equation}
Denote by $S_{k+2}^j$ the  set $\big\{\tau\in S_{k+2}\big|\tau(k+2)=j\big\}$. Notice  that \begin{equation}
sgn(\widehat\sigma^j)=sgn(\sigma)\cdot(-1)^{k+2-j},
\end{equation}
and 
\begin{equation}
S_{k+2}=\bigcup_{j=1}^{k+2}S_{k+2}^j\,\,{\rm where \,\,}S_{k+2}^i\cap S_{k+2}^j=\emptyset\,\,{\rm for\,\,} i\neq j.
\end{equation}

We then have the following formula under the above notation for $j=1,\cdots,k+2$.
\begin{equation}\label{ij}
\begin{split}
(-1)^{j+1}t_{i_1\cdots\widehat {i_j}\cdots i_{k+2}}=&\sum_{\sigma\in S_{k+1}}(-1)^{j+1}\cdot\frac{sgn(\sigma)}{(k+1)!}\cdot tr\big(g^{-1}_{i_{\widehat\sigma^j(1)}i_{\widehat\sigma^j(k+1)}}dg_{i_{\widehat\sigma^j(1)}i_{\widehat\sigma^j(k+1)}}g^{-1}_{i_{\widehat\sigma^j(2)}i_{\widehat\sigma^j(k+1)}}dg_{i_{\widehat\sigma^j(2)}i_{\widehat\sigma^j(k+1)}}\\
&\,\,\,\,\,\,\,\,g^{-1}_{i_{\widehat\sigma^j(3)}i_{\widehat\sigma^j(k+1)}}dg_{i_{\widehat\sigma^j(3)}i_{\widehat\sigma^j(k+1)}}\cdots g^{-1}_{i_{\widehat\sigma^j(k)}i_{\widehat\sigma^j(k+1)}}dg_{i_{\widehat\sigma^j(k)}i_{\widehat\sigma^j(k+1)}}\big)\\
=&\sum_{\sigma\in S_{k+1}}(-1)^{k+1}\cdot\frac{sgn(\widehat\sigma^j)}{(k+1)!}\cdot tr\big(g^{-1}_{i_{\widehat\sigma^j(1)}i_{\widehat\sigma^j(k+1)}}dg_{i_{\widehat\sigma^j(1)}i_{\widehat\sigma^j(k+1)}}g^{-1}_{i_{\widehat\sigma^j(2)}i_{\widehat\sigma^j(k+1)}}dg_{i_{\widehat\sigma^j(2)}i_{\widehat\sigma^j(k+1)}}\\
&\,\,\,\,\,\,\,\,g^{-1}_{i_{\widehat\sigma^j(3)}i_{\widehat\sigma^j(k+1)}}dg_{i_{\widehat\sigma^j(3)}i_{\widehat\sigma^j(k+1)}}\cdots g^{-1}_{i_{\widehat\sigma^j(k)}i_{\widehat\sigma^j(k+1)}}dg_{i_{\widehat\sigma^j(k)}i_{\widehat\sigma^j(k+1)}}\big)\\
=&\sum_{\widehat\sigma\in S^j_{k+2}}(-1)^{k+1}\cdot\frac{sgn(\widehat\sigma)}{(k+1)!}\cdot tr\big(g^{-1}_{i_{\widehat\sigma(1)}i_{\widehat\sigma(k+1)}}dg_{i_{\widehat\sigma(1)}i_{\widehat\sigma(k+1)}}g^{-1}_{i_{\widehat\sigma(2)}i_{\widehat\sigma(k+1)}}dg_{i_{\widehat\sigma(2)}i_{\widehat\sigma(k+1)}}\\
&\,\,\,\,\,\,\,\,g^{-1}_{i_{\widehat\sigma(3)}i_{\widehat\sigma(k+1)}}dg_{i_{\widehat\sigma(3)}i_{\widehat\sigma(k+1)}}\cdots g^{-1}_{i_{\widehat\sigma(k)}i_{\widehat\sigma(k+1)}}dg_{i_{\widehat\sigma(k)}i_{\widehat\sigma(k+1)}}\big).
\end{split}
\end{equation}
Substituting into formula (\ref{cocy}), we derive that
\begin{equation}\label{I_1}
\begin{split}
\sum_{j=1}^{k+2}(-1)^{j+1}&t_{i_1\cdots\widehat {i_j}\cdots i_{k+2}}=\sum_{j=1}^{k+2}\sum_{\widehat\sigma\in S^j_{k+2}}(-1)^{k+1}\cdot\frac{sgn(\widehat\sigma)}{(k+1)!}\cdot tr\big(g^{-1}_{i_{\widehat\sigma(1)}i_{\widehat\sigma(k+1)}}dg_{i_{\widehat\sigma(1)}i_{\widehat\sigma(k+1)}}\\
&g^{-1}_{i_{\widehat\sigma(2)}i_{\widehat\sigma(k+1)}}dg_{i_{\widehat\sigma(2)}i_{\widehat\sigma^j(k+1)}}g^{-1}_{i_{\widehat\sigma(3)}i_{\widehat\sigma(k+1)}}dg_{i_{\widehat\sigma(3)}i_{\widehat\sigma(k+1)}}\cdots g^{-1}_{i_{\widehat\sigma(k)}i_{\widehat\sigma(k+1)}}dg_{i_{\widehat\sigma(k)}i_{\widehat\sigma(k+1)}}\big)\\
&=\sum_{\widehat\sigma\in S_{k+2}}(-1)^{k+1}\cdot\frac{sgn(\widehat\sigma)}{(k+1)!}\cdot tr\big(g^{-1}_{i_{\widehat\sigma(1)}i_{\widehat\sigma(k+1)}}dg_{i_{\widehat\sigma(1)}i_{\widehat\sigma(k+1)}}\\
&g^{-1}_{i_{\widehat\sigma(2)}i_{\widehat\sigma(k+1)}}dg_{i_{\widehat\sigma(2)}i_{\widehat\sigma(k+1)}}g^{-1}_{i_{\widehat\sigma(3)}i_{\widehat\sigma(k+1)}}dg_{i_{\widehat\sigma(3)}i_{\widehat\sigma(k+1)}}\cdots g^{-1}_{i_{\widehat\sigma(k)}i_{\widehat\sigma(k+1)}}dg_{i_{\widehat\sigma(k)}i_{\widehat\sigma(k+1)}}\big).
\end{split}
\end{equation}

Notice that $g_{\alpha\beta}=g_{\alpha\gamma}\cdot g_{\gamma\beta}$ and hence $dg_{\alpha\beta}=dg_{\alpha\gamma}\cdot g_{\gamma\beta}+g_{\alpha\gamma}\cdot dg_{\gamma\beta}$. Then,  we have 
\begin{equation}\label{tra}
g^{-1}_{\alpha\beta}dg_{\alpha\beta}=\big(g^{-1}_{\gamma\beta}g^{-1}_{\alpha\gamma}dg_{\alpha\gamma}g_{\gamma\beta}+g^{-1}_{\gamma\beta} dg_{\gamma\beta}\big)=g^{-1}_{\gamma\beta}\big(g^{-1}_{\alpha\gamma}dg_{\alpha\gamma}+(-1)g^{-1}_{\beta\gamma} dg_{\beta\gamma}\big)g_{\gamma\beta}.
\end{equation}
Substituting formula (\ref{tra}) into formula (\ref{ij}) by taking $\gamma=i_{\widehat\sigma(k+2)}$, $\beta=i_{\widehat\sigma(k+1)}$ and $\alpha=i_{\widehat\sigma(1)},i_{\widehat\sigma(2)},i_{\widehat\sigma(3)},\cdots,i_{\widehat\sigma(k)}$, we get
\begin{equation}\label{traij}
\begin{split}
&(-1)^{j+1}t_{i_1\cdots\widehat {i_j}\cdots i_{k+2}}=\sum_{\widehat\sigma\in S^j_{k+2}}(-1)^{k+1}\cdot\frac{sgn(\widehat\sigma)}{(k+1)!}\cdot tr\Big\{\prod_{l=1}^kg^{-1}_{i_{\widehat\sigma(k+2)}i_{\widehat\sigma(k+1)}}\big(g^{-1}_{i_{\widehat\sigma(l)}i_{\widehat\sigma(k+2)}}dg_{i_{\widehat\sigma(l)}i_{\widehat\sigma(k+2)}}\\ &\,\,\,\,\,\,\,\,\,\,\,\,\,\,\,\,\,\,\,\,\,\,\,\,\,\,\,\,\,\,\,\,\,\,\,\,\,\,\,\,\,\,\,\,\,\,\,\,\,\,\,\,\,\,\,\,\,\,\,\,+(-1)g^{-1}_{i_{\widehat\sigma(k+1)}i_{\widehat\sigma(k+2)}} dg_{i_{\widehat\sigma(k+1)}i_{\widehat\sigma(k+2)}}\big)g_{i_{\widehat\sigma(k+2)}i_{\widehat\sigma(k+1)}}\Big\}\\
&=\sum_{\widehat\sigma\in S^j_{k+2}}(-1)^{k+1}\cdot\frac{sgn(\widehat\sigma)}{(k+1)!}\cdot tr\big\{\prod_{l=1}^k(g^{-1}_{i_{\widehat\sigma(l)}i_{\widehat\sigma(k+2)}}dg_{i_{\widehat\sigma(l)}i_{\widehat\sigma(k+2)}}-g^{-1}_{i_{\widehat\sigma(k+1)}i_{\widehat\sigma(k+2)}} dg_{i_{\widehat\sigma(k+1)}i_{\widehat\sigma(k+2)}}\big)\big\}.
\end{split}
\end{equation}
Notice that in the last step we use the fact that $tr(ABC)=tr(BCA)$.

To simplify the notation, for elements $p_1,\cdots,p_{k+1}\in I$ (not necessarily distinct), we define
\begin{equation}
\Delta(p_1,\cdots,p_{k+1}):=tr\big(g^{-1}_{p_1p_{k+1}}dg_{p_1p_{k+1}}g^{-1}_{p_2p_{k+1}}dg_{p_2p_{k+1}}g^{-1}_{p_3p_{k+1}}dg_{p_3p_{k+1}}\cdots g^{-1}_{p_kp_{k+1}}dg_{p_kp_{k+1}}\big);
\end{equation}
for positive integer $j$ and $s_1,\cdots,s_j$ such that  $1\leq j\leq k$ and $1\leq s_1<s_2<\cdots<s_j\leq k$, and an element $p_{k+2}\in I$, define 
\begin{equation}\label{expan}
\begin{split}
\Delta_{s_1\cdots s_j,\,p_{k+2}}&(p_1,\cdots,p_{k+1}):=\Delta(p_1,\cdots,\widehat{p_{s_1}},p_{k+2},p_{s_1+1},\cdots,\widehat{p_{s_j}},p_{k+2},p_{s_j+1},\cdots, p_{k+1})\,\,, 
\end{split}
\end{equation}
namely, one replaces the indices $p_{s_1},\cdots,p_{s_j}$ by $p_{k+2}$ in $\Delta(p_1,\cdots,p_{k+1})$.

Under the above notation, we can expand formula (\ref{traij}) and derive that
\begin{equation}\label{expansion}
\begin{split}
&(-1)^{j+1}t_{i_1\cdots\widehat {i_j}\cdots i_{k+2}}=\frac{(-1)^{k+1}}{(k+1)!}\sum_{\widehat\sigma\in S^j_{k+2}}sgn(\widehat\sigma)\cdot\Big\{ \Delta(i_{\widehat\sigma(1)},i_{\widehat\sigma(2)},\cdots,i_{\widehat\sigma(k)},i_{\widehat\sigma(k+2)})\\
&\,\,\,\,-\Delta(i_{\widehat\sigma(k+1)},i_{\widehat\sigma(2)},i_{\widehat\sigma(3)},\cdots,i_{\widehat\sigma(k)},i_{\widehat\sigma(k+2)})-\Delta(i_{\widehat\sigma(1)},i_{\widehat\sigma(k+1)},i_{\widehat\sigma(3)},\cdots,i_{\widehat\sigma(k)},i_{\widehat\sigma(k+2)})\\
&\,\,\,\,\,\,\,\,\,-\cdots-\Delta(i_{\widehat\sigma(1)},i_{\widehat\sigma(2)},i_{\widehat\sigma(3)},\cdots,i_{\widehat\sigma(k-1)},i_{\widehat\sigma(k+1)},i_{\widehat\sigma(k+2)})
\Big\}\\
&+\frac{(-1)^{k+1}}{(k+1)!}\sum_{\widehat\sigma\in S^j_{k+2}}sgn(\widehat\sigma) 
\sum_{\substack{2\leq m\leq k \\ 1\leq s_1<\cdots<s_m\leq k}}(-1)^m\Delta_{s_1\cdots s_m,\,i_{\widehat\sigma(k+1)}}(i_{\widehat\sigma(1)},i_{\widehat\sigma(2)},\cdots,i_{\widehat\sigma(k)},i_{\widehat\sigma(k+2)}).
\end{split}
\end{equation}

\noindent{\bf Claim I :} For $1\leq j\leq k+2$, $2\leq m\leq k$ and $1\leq s_1<\cdots<s_m\leq k$, we have that
\begin{equation}
\sum_{\widehat\sigma\in S^j_{k+2}}sgn(\widehat\sigma)\cdot \Delta_{s_1\cdots s_m,\,i_{\widehat\sigma(k+1)}}(i_{\widehat\sigma(1)},i_{\widehat\sigma(2)},\cdots,i_{\widehat\sigma(k)},i_{\widehat\sigma(k+2)})=0\,\,.
\end{equation}
\noindent{\bf Proof of Claim I :} Denote by $\tau$ the permutation $(s_1,s_2)\in S_{k+2}$. Since $1\leq s_1<s_2\leq k$, $\widehat\sigma\circ\tau\in S_{k+2}^j$ for each $\widehat\sigma\in S_{k+2}^j$. It is easy to verify that the following composition map is a bijection for $1\leq j\leq k+2$,
\begin{equation}\label{iso}
\begin{split}
T_{\tau}:S_{k+2}^j\xrightarrow{\cong}S_{k+2}^j\\
\widehat\sigma\mapsto\widehat\sigma\circ\tau
\end{split}\,\,\,\,\,\,.
\end{equation}
Moreover, we have that
\begin{equation}\label{change}
\begin{split}
\Delta_{s_1\cdots s_m,i_{\widehat\sigma(k+1)}}&(i_{\widehat\sigma(1)},i_{\widehat\sigma(2)},\cdots,i_{\widehat\sigma(k)},i_{\widehat\sigma(k+2)})\\
&=\Delta_{s_1\cdots s_m,i_{(\widehat\sigma\circ\tau)(k+1)}}(i_{(\widehat\sigma\circ\tau)(1)},i_{(\widehat\sigma\circ\tau)(2)},\cdots,i_{(\widehat\sigma\circ\tau)(k)},i_{(\widehat\sigma\circ\tau)(k+2)}).
\end{split}
\end{equation}
Combining bijection map (\ref{iso}) and equality (\ref{change}), we derive that
\begin{equation}\label{t1j}
\begin{split}
\sum_{\widehat\sigma\in S^j_{k+2}}&sgn(\widehat\sigma)\cdot \Delta_{s_1\cdots s_m,i_{\widehat\sigma(k+1)}}(i_{\widehat\sigma(1)},i_{\widehat\sigma(2)},\cdots,i_{\widehat\sigma(k)},i_{\widehat\sigma(k+2)})\\
=&\sum_{\widehat\sigma\circ\tau\in S^j_{k+2}}sgn(\widehat\sigma\circ\tau)\cdot\Delta_{s_1\cdots s_m,i_{(\widehat\sigma\circ\tau)(k+1)}}(i_{(\widehat\sigma\circ\tau)(1)},i_{(\widehat\sigma\circ\tau)(2)},\cdots,i_{(\widehat\sigma\circ\tau)(k)},i_{(\widehat\sigma\circ\tau)(k+2)})\\
=&sgn(\tau)\cdot\sum_{\widehat\sigma\circ\tau\in S^j_{k+2}}sgn(\widehat\sigma)\cdot \Delta_{s_1\cdots s_m,i_{\widehat\sigma(k+1)}}(i_{\widehat\sigma(1)},i_{\widehat\sigma(2)},\cdots,i_{\widehat\sigma(k)},i_{\widehat\sigma(k+2)})\\
=&-\sum_{\widehat\sigma\in S^j_{k+2}}sgn(\widehat\sigma)\cdot \Delta_{s_1\cdots s_m,i_{\widehat\sigma(k+1)}}(i_{\widehat\sigma(1)},i_{\widehat\sigma(2)},\cdots,i_{\widehat\sigma(k)},i_{\widehat\sigma(k+2)}).
\end{split}
\end{equation}
Therefore, we have
\begin{equation}
\sum_{\widehat\sigma\in S^j_{k+2}}sgn(\widehat\sigma)\cdot \Delta_{s_1\cdots s_m,i_{\widehat\sigma(k+1)}}(i_{\widehat\sigma(1)},i_{\widehat\sigma(2)},\cdots,i_{\widehat\sigma(k)},i_{\widehat\sigma(k+2)})=0\,\,.
\end{equation}
We complete the proof of Claim I. 
\,\,\,\,$\endpf$
\medskip

Denote the permutation $(q,k+1)\in S_{k+2}$ by $\tau_q$ for $1\leq q\leq k+2$. Then, the following map is a bijection for $1\leq q\leq k+2$,
\begin{equation}\label{tiso}
\begin{split}
T_{\tau_q}:S_{k+2}\xrightarrow{\cong}S_{k+2}\\
\widehat\sigma\mapsto\widehat\sigma\circ\tau_q
\end{split}\,\,\,\,\,\,;
\end{equation}
$sgn(\tau_q)=-1$ for $q\neq k+1$  and $sgn(\tau_q)=1$ for $q=k+1$.

By Claim I and formula (\ref{expansion}), we have that
\begin{equation}\label{expansion1}
\begin{split}
&\sum_{j=1}^{k+2}(-1)^{j+1}t_{i_1\cdots\widehat {i_j}\cdots i_{k+2}}=\frac{(-1)^{k+1}}{(k+1)!}\sum_{j=1}^{k+2}\sum_{\widehat\sigma\in S^j_{k+2}}sgn(\widehat\sigma)\cdot\Big\{ \Delta(i_{\widehat\sigma(1)},i_{\widehat\sigma(2)},\cdots,i_{\widehat\sigma(k)},i_{\widehat\sigma(k+2)})\\
&\,\,\,\,-\Delta(i_{\widehat\sigma(k+1)},i_{\widehat\sigma(2)},i_{\widehat\sigma(3)},\cdots,i_{\widehat\sigma(k)},i_{\widehat\sigma(k+2)})-\Delta(i_{\widehat\sigma(1)},i_{\widehat\sigma(k+1)},i_{\widehat\sigma(3)},\cdots,i_{\widehat\sigma(k)},i_{\widehat\sigma(k+2)})\\
&\,\,\,\,\,\,\,\,\,-\cdots-\Delta(i_{\widehat\sigma(1)},i_{\widehat\sigma(2)},i_{\widehat\sigma(3)},\cdots,i_{\widehat\sigma(k-1)},i_{\widehat\sigma(k+1)},i_{\widehat\sigma(k+2)})\Big\}\\
&\,\,\,\,=\frac{(-1)^{k+1}}{(k+1)!}\sum_{\widehat\sigma\in S_{k+2}}sgn(\widehat\sigma)\cdot\Big\{+\Delta(i_{\widehat\sigma(1)},i_{\widehat\sigma(2)},\cdots,i_{\widehat\sigma(k)},i_{\widehat\sigma(k+2)}) \\
&\,\,\,\,\,\,\,\,\,\,\,\,\,\,\,-\sum_{q=1}^k\Delta(i_{(\widehat\sigma\circ\tau_q)(1)},i_{(\widehat\sigma\circ\tau_q)(2)},i_{(\widehat\sigma\circ\tau_q)(3)},\cdots,i_{(\widehat\sigma\circ\tau_q)(k)},i_{(\widehat\sigma\circ\tau_q)(k+2)})\Big\}\\
&\,\,\,=\frac{(-1)^{k+1}}{(k+1)!}\sum_{\widehat\sigma\in S_{k+2}}sgn(\widehat\sigma)\cdot\Big\{ (k+1)\Delta(i_{\widehat\sigma(1)},i_{\widehat\sigma(2)},\cdots,i_{\widehat\sigma(k)},i_{\widehat\sigma(k+2)})\Big\}.
\end{split}
\end{equation}

Together with formula (\ref{I_1}), we have that
\begin{equation}\label{tkj}
\begin{split}
&\sum_{j=1}^{k+2}(-1)^{j+1}t_{i_1\cdots\widehat {i_j}\cdots i_{k+2}}=\frac{(-1)^{k+1}}{(k+1)!}\sum_{\widehat\sigma\in S_{k+2}}sgn(\widehat\sigma)\cdot\Big\{ (k+1)\Delta(i_{\widehat\sigma(1)},i_{\widehat\sigma(2)},\cdots,i_{\widehat\sigma(k)},i_{\widehat\sigma(k+2)})\Big\}\\
&=\frac{(-1)^{k+1}(k+1)\cdot sgn(\tau_{k+2})}{(k+1)!}\sum_{\widehat\sigma\circ\tau_{k+2}\in S_{k+2}}sgn(\widehat\sigma\circ\tau_{k+2})\cdot \Delta(i_{(\widehat\sigma\circ\tau_{k+2})(1)},\cdots,i_{(\widehat\sigma\circ\tau_{k+2})(k)},i_{(\widehat\sigma\circ\tau_{k+2})(k+1)})\\
&=-(k+1)\cdot\sum_{j=1}^{k+2}(-1)^{j+1}t_{i_1\cdots\widehat {i_j}\cdots i_{k+2}}.
\end{split}
\end{equation}
Therefore, we complete the proof of Lemma \ref{cocycle}.
\,\,\,\,$\endpf$
\medskip

By Lemma \ref{cocycle}, $\widehat f_k(E,g)$ determines an element in $\check H^k(\mathcal U,\Omega^k)$; we denote this element by $f_k(E,g)$. Next, we will prove that $f_k(E,g)$ is independent of the choice of $g$.
 
\begin{lemma}\label{invariant}
	The element $f_k(E,g)\in \check H^k(\mathcal U,\Omega^k)$ is independent of the choice of the system of transition functions of $E$ with respect to the open cover $\mathcal U$; namely, for $g$ and $\widetilde g$  systems of transition functions   with respect to $\mathcal U$, $f_k(E,g)=f_k(E,\widetilde g)$ in $ \check H^k(\mathcal U,\Omega^k)$.
\end{lemma}

\noindent{\bf Proof of Lemma \ref{invariant} :} 
Since $\{g_{ij}\}$ and $\{\widetilde g_{ij}\}$ are two systems of transition functions of $E$ with respect to $\mathcal U$,  there exists a \v Cech $0$-cochain
\begin{equation}
h:=\bigoplus\limits_{i_1}h_{i_1}\in\bigoplus\limits_{i_1}\Gamma_{\rm hol} (U_{i_1},GL(M,\mathbb C)),
\end{equation}
such that $\widetilde g_{ij}=h_i^{-1}g_{ij}h_j$ for $i,j\in I$. Let $\widehat f_k(E,g)$, $\widehat f_k(E,\widetilde g)$ be the \v Cech $k$-cocycles associated with $f_k(E,g)$, $f_k(E,g)$, respectively, as follows,
\begin{equation}
\widehat f_k(E,g)=\bigoplus\limits_{i_1<\cdots<i_{k+1}}t_{i_1\cdots i_{k+1}}\in\bigoplus\limits_{i_1<\cdots<i_{k+1}}\Gamma(U_{i_1\cdots i_{k+1}},\Omega^k),
\end{equation}
\begin{equation}
\widehat f_k(E,\widetilde g)=\bigoplus\limits_{i_1<\cdots<i_{k+1}}\widetilde t_{i_1\cdots i_{k+1}}\in\bigoplus\limits_{i_1<\cdots<i_{k+1}}\Gamma(U_{i_1\cdots i_{k+1}},\Omega^k).
\end{equation}
Here $t_{i_1\cdots i_{k+1}}$ and $\widetilde t_{i_1\cdots i_{k+1}}$ are defined by formula (\ref{gcf}) with respect to $g$ and $\widetilde g$, respectively.

To prove Lemma \ref{invariant}, it suffices to prove that there is a \v Cech $(k-1)$-cochain
\begin{equation}
h_{k-1}(E,g,\widetilde g)=\bigoplus\limits_{j_1<\cdots<j_{k}}s_{j_1\cdots j_{k}}\in\bigoplus\limits_{j_1<\cdots<j_{k}}\Gamma(U_{j_1\cdots j_{k}},\Omega^k),
\end{equation}
such that for any $i_1,\cdots,i_{k+1}\in I$, we have
\begin{equation}\label{cobound}
\widetilde t_{i_1\cdots i_{k+1}}- t_{i_1\cdots i_{k+1}}=\sum_{j=1}^{k+1}(-1)^{j-1}s_{i_1\cdots\widehat i_j\cdots i_{k+1}}\big|_{U_{i_1\cdots i_{k+1}}}.
\end{equation}

Notice that 
\begin{equation}\label{tra2}
\widetilde g_{\alpha\beta}^{-1}d\widetilde g_{\alpha\beta}=h_{\beta}^{-1}\big(g_{\alpha\beta}^{-1}dg_{\alpha\beta}+g_{\alpha\beta}^{-1}h_{\alpha}dh_{\alpha}^{-1}g_{\alpha\beta}+(-1)\cdot h_{\beta}dh_{\beta}^{-1}\big)h_{\beta}.
\end{equation}
Applying formula (\ref{tra2}) for $\beta=i_{\sigma(k+1)}$ and $\alpha=i_{\sigma(1)},i_{\sigma(2)},i_{\sigma(3)},\cdots,i_{\sigma(k)}$, we get
\begin{equation}\label{tildet}
\begin{split}
\widetilde t_{i_1\cdots i_{k+1}}&=\sum_{\sigma\in S_{k+1}}\frac{sgn(\sigma)}{(k+1)!}\cdot tr\big(\widetilde g^{-1}_{i_{\sigma(1)}i_{\sigma(k+1)}}d\widetilde g_{i_{\sigma(1)}i_{\sigma(k+1)}}\widetilde g^{-1}_{i_{\sigma(2)}i_{\sigma(k+1)}}d\widetilde g_{i_{\sigma(2)}i_{\sigma(k+1)}}\widetilde g^{-1}_{i_{\sigma(3)}i_{\sigma(k+1)}}d\widetilde g_{i_{\sigma(3)}i_{\sigma(k+1)}}\\
&\,\,\,\,\,\,\,\,\,\,\,\,\,\,\,\,\,\,\,\,\,\,\,\,\,\,\,\,\,\,\,\,\,\,\,\,\,\,\,\,\,\,\,\,\,\,\cdots \widetilde g^{-1}_{i_{\sigma(k)}i_{\sigma(k+1)}}d\widetilde g_{i_{\sigma(k)}i_{\sigma(k+1)}}\big)\\
&=\sum_{\sigma\in S_{k+1}}\frac{sgn(\sigma)}{(k+1)!}\cdot tr\big\{\prod_{l=1}^kh_{i_{\sigma(k+1)}}^{-1}\big(g_{i_{\sigma(l)} i_{\sigma(k+1)}}^{-1}dg_{i_{\sigma(l)} i_{\sigma(k+1)}}+g_{i_{\sigma(l)} i_{\sigma(k+1)}}^{-1}h_{i_{\sigma(l)}}dh_{i_{\sigma(l)}}^{-1}g_{i_{\sigma(l)}i_{\sigma(k+1)}}\\
&\,\,\,\,\,\,\,\,\,\,\,\,\,\,\,\,\,\,\,\,\,\,\,\,\,\,\,\,\,\,\,\,\,\,\,\,\,\,\,\,\,\,\,\,\,\,\,\,\,\,\,\,\,\,\,\,\,\,\,\,\,\,\,\,\,\,\,+(-1)\cdot h_{i_{\sigma(k+1)}}dh_{i_{\sigma(k+1)}}^{-1}\big)h_{i_{\sigma(k+1)}}\big\}\\
&=\sum_{\sigma\in S_{k+1}}\frac{sgn(\sigma)}{(k+1)!}\cdot tr\big\{\prod_{l=1}^k\big(g_{i_{\sigma(l)} i_{\sigma(k+1)}}^{-1}dg_{i_{\sigma(l)} i_{\sigma(k+1)}}+g_{i_{\sigma(l)}i_{\sigma(k+1)}}^{-1}h_{i_{\sigma(l)}}dh_{i_{\sigma(l)}}^{-1}g_{i_{\sigma(l)}i_{\sigma(k+1)}}\\
&\,\,\,\,\,\,\,\,\,\,\,\,\,\,\,\,\,\,\,\,\,\,\,\,\,\,\,\,\,\,\,\,\,\,\,\,\,\,\,\,\,\,\,\,\,\,\,\,\,\,\,\,\,\,\,\,\,\,\,\,\,\,\,\,\,\,\,+(-1)\cdot h_{i_{\sigma(k+1)}}dh_{i_{\sigma(k+1)}}^{-1}\big)\big\}.
\end{split}
\end{equation}

In order to compute the monomials appearing in the expansion of formula (\ref{tildet}) by distribution law, we introduce the following notation. For integer $t\geq 0$ and integers $a_1,\cdots,a_t$, we define the counting function $\delta_{a_1\cdots a_t}$ from $\mathbb Z$ to the set $\{0,1\}$ as follows,
\begin{equation}
\delta_{a_1\cdots a_t}(l)=\left\{\begin{matrix}
1&\,\,\,\,\,\,\,\,\,\,\,{\rm if}\,\,\,\,\,\,\,\,\,\,\,l\in\{a_1,\cdots,a_t\};\\0&\,\,\,\,\,\,\,\,\,\,\,{\rm if}\,\,\,\,\,\,\,\,\,\,\,l\notin\{a_1,\cdots,a_t\}.
\end{matrix}\right  .
\end{equation}
Notice that when $t=0$, $\delta_{a_1\cdots a_t}\equiv 0$; we denote it by $\delta_{\emptyset}$.

Fix a permutation $\sigma\in S_{k+1}$. For integers $j, n$ such that $j\geq0$, $n\geq0$ and $j+n\leq k$, and distinct integers $u_1,\cdots,u_j,v_1,\cdots,v_n\in\{1,\cdots,k\}$ (by convention if $j=0$, there is no $u$; if $n=0$, there is no $v$), we define 
\begin{equation}\label{dis}
\begin{split}
&\Delta^{\sigma}_{u_1\cdots u_j,v_{1}\cdots v_{n}}(\widetilde t_{i_1\cdots i_{k+1}}):=tr\big\{\prod_{l=1}^k\big(g_{i_{\sigma(l)} i_{\sigma(k+1)}}^{-1}dg_{i_{\sigma(l)} i_{\sigma(k+1)}}\big)^{1-\delta_{u_1\cdots u_j}(l)-\delta_{v_1\cdots v_n}(l)}\\
&\cdot\big(g_{i_{\sigma(l)}i_{\sigma(k+1)}}^{-1}h_{i_{\sigma(l)}}dh_{i_{\sigma(l)}}^{-1}g_{i_{\sigma(l)}i_{\sigma(k+1)}}\big)^{\delta_{u_1\cdots u_j}(l)}\cdot\big(-h_{i_{\sigma(k+1)}}dh_{i_{\sigma(k+1)}}^{-1}\big)^{\delta_{v_1\cdots v_n}(l)}\big\}.
\end{split}
\end{equation}
Note that here we use the convention that 
\begin{equation}
\big(g_{i_{\sigma(l)} i_{\sigma(k+1)}}^{-1}dg_{i_{\sigma(l)} i_{\sigma(k+1)}}\big)^{0}=
\big(g_{i_{\sigma(l)}i_{\sigma(k+1)}}^{-1}h_{i_{\sigma(l)}}dh_{i_{\sigma(l)}}^{-1}g_{i_{\sigma(l)}i_{\sigma(k+1)}}\big)^{0}=\big(-h_{i_{\sigma(k+1)}}dh_{i_{\sigma(k+1)}}^{-1}\big)^{0}=I_{M\times M},
\end{equation}
where $I_{M\times M}$ is the identity matrix of rank $M$.  We also denote $\Delta^{\sigma}_{u_1\cdots u_j,v_{1}\cdots v_{n}}(\widetilde t_{i_1\cdots i_{k+1}})$ by $\Delta^{\sigma}_{\emptyset,v_{1}\cdots v_{n}}(\widetilde t_{i_1\cdots i_{k+1}})$ when $j=0$ and $n\geq 1$;  denote $\Delta^{\sigma}_{u_1\cdots u_j,v_{1}\cdots v_{n}}(\widetilde t_{i_1\cdots i_{k+1}})$ by $\Delta^{\sigma}_{u_{1}\cdots u_{j},\emptyset}(\widetilde t_{i_1\cdots i_{k+1}})$  when $n=0$ and $j\geq 1$; denote $\Delta^{\sigma}_{u_1\cdots u_j,v_{1}\cdots v_{n}}(\widetilde t_{i_1\cdots i_{k+1}})$ by $\Delta^{\sigma}_{\emptyset,\emptyset}(\widetilde t_{i_1\cdots i_{k+1}})$ when $j=n=0$.
\medskip

Applying distribution law to formula (\ref{tildet}), we have
\begin{equation}\label{tildet2}
\begin{split}
&\widetilde t_{i_1\cdots i_{k+1}}=\sum_{\sigma\in S_{k+1}}\frac{sgn(\sigma)}{(k+1)!} \sum_{\substack{n\geq 2,\,0\leq j\leq k-n\,; \\  1\leq u_1<u_2<\cdots<u_j\leq k\,;\\ 1\leq v_1<v_2<\cdots<v_n\leq k\,;\\u_1,\cdots,u_j,v_1,\cdots,v_n\, {\text{are distinct}}}}\Delta^{\sigma}_{u_1\cdots u_j,v_{1}\cdots v_{n}}(\widetilde t_{i_1\cdots i_{k+1}})\\
&+\sum_{\sigma\in S_{k+1}}\frac{sgn(\sigma)}{(k+1)!} \sum_{\substack{n=1,\,1\leq j\leq k-1\,; \\  1\leq u_1<\cdots<u_j\leq k\,;\\ 1\leq v_1\leq k\,;\\u_1,\cdots,u_j,v_1\, {\text{are distinct}}}}\Delta^{\sigma}_{u_1\cdots u_j,v_{1}}(\widetilde t_{i_1\cdots i_{k+1}})+\sum_{\sigma\in S_{k+1}}\frac{sgn(\sigma)}{(k+1)!} \sum_{\substack{n=1,\, j=0\,; \\ 1\leq v_1\leq k\,;}}\Delta^{\sigma}_{\emptyset,v_{1}}(\widetilde t_{i_1\cdots i_{k+1}})\\
&+\sum_{\sigma\in S_{k+1}}\frac{sgn(\sigma)}{(k+1)!} \sum_{\substack{n=0,\,1\leq j\leq k\,; \\  1\leq u_1<u_2<\cdots<u_j\leq k\,;}}\Delta^{\sigma}_{u_1\cdots u_j,\emptyset}(\widetilde t_{i_1\cdots i_{k+1}})+\sum_{\sigma\in S_{k+1}}\frac{sgn(\sigma)}{(k+1)!}\cdot \Delta^{\sigma}_{\emptyset,\emptyset}(\widetilde t_{i_1\cdots i_{k+1}}).\\
&=:I_1+I_2+I_3+I_4+I_5.
\end{split}
\end{equation}

We first compute $I_1$ as follows.
\medskip

\noindent{\bf Claim I : } 
For integers $j, n$ such that $n\geq2$ and $0\leq j\leq k-n$, and distinct integers $u_1,\cdots,u_j,v_1,\cdots,v_n\in\{1,\cdots,k\}$,  the following identity holds,
\begin{equation}
\sum_{\sigma\in S_{k+1}}\frac{sgn(\sigma)}{(k+1)!}\cdot \Delta^{\sigma}_{u_1\cdots u_j,v_{1}\cdots v_{n}}(\widetilde t_{i_1\cdots i_{k+1}})=0.
\end{equation}
Therefore, $I_1\equiv 0$.
\medskip

\noindent{\bf Proof of Claim I :} The proof is similar to the proof of Claim I of Lemma \ref{cocycle}.  Denote by $\tau$ the permutation $(v_1,v_2)\in S_{k+1}$. Notice that
\begin{equation}\label{changeh}
\begin{split}
&\Delta^{\sigma\circ\tau}_{u_1\cdots u_j,v_{1}\cdots v_{n}}(\widetilde t_{i_1\cdots i_{k+1}})=tr\big\{\prod_{l=1}^k\big(g_{i_{(\sigma\circ\tau)(l)} i_{(\sigma\circ\tau)(k+1)}}^{-1}dg_{i_{(\sigma\circ\tau)(l)} i_{(\sigma\circ\tau)(k+1)}}\big)^{1-\delta_{u_1\cdots u_j}(l)-\delta_{v_1\cdots v_n}(l)}\\
&\cdot\big(g_{i_{(\sigma\circ\tau)(l)}i_{(\sigma\circ\tau)(k+1)}}^{-1}h_{i_{(\sigma\circ\tau)(l)}}dh_{i_{(\sigma\circ\tau)(l)}}^{-1}g_{i_{(\sigma\circ\tau)(l)}i_{(\sigma\circ\tau)(k+1)}}\big)^{\delta_{u_1\cdots u_j}(l)}\cdot\big(-h_{i_{(\sigma\circ\tau)(k+1)}}dh_{i_{(\sigma\circ\tau)(k+1)}}^{-1}\big)^{\delta_{v_1\cdots v_n}(l)}\big\}\\
&=tr\big\{\prod_{l=1}^k\big(g_{i_{\sigma(l)} i_{\sigma(k+1)}}^{-1}dg_{i_{\sigma(l)} i_{\sigma(k+1)}}\big)^{1-\delta_{u_1\cdots u_j}(l)-\delta_{v_1\cdots v_n}(l)}\cdot\big(g_{i_{\sigma(l)}i_{\sigma(k+1)}}^{-1}h_{i_{\sigma(l)}}dh_{i_{\sigma(l)}}^{-1}g_{i_{\sigma(l)}i_{\sigma(k+1)}}\big)^{\delta_{u_1\cdots u_j}(l)}\\
&\cdot\big(-h_{i_{\sigma(k+1)}}dh_{i_{\sigma(k+1)}}^{-1}\big)^{\delta_{v_1\cdots v_n}(l)}\big\}=\Delta^{\sigma}_{u_1\cdots u_j,v_1v_2\cdots v_n}(\widetilde t_{i_1\cdots i_{k+1}}).
\end{split}
\end{equation}
Hence,
\begin{equation}\label{t2j}
\begin{split}
&\sum_{\sigma\in S_{k+1}}\frac{sgn(\sigma)}{(k+1)!}\cdot\Delta^{\sigma}_{u_1\cdots u_j,v_1v_2\cdots v_n}(\widetilde t_{i_1\cdots i_{k+1}})
=\sum_{\sigma\circ\tau\in S_{k+1}}\frac{sgn(\sigma\circ\tau)}{(k+1)!}\cdot\Delta^{\sigma\circ\tau}_{u_1\cdots u_j,v_{1}\cdots v_{n}}(\widetilde t_{i_1\cdots i_{k+1}})\\
&=sgn(\tau)\cdot\sum_{\sigma\in S_{k+1}}\frac{sgn(\sigma)}{(k+1)!}\cdot \Delta^{\sigma\circ\tau}_{u_1\cdots u_j,v_{1}\cdots v_{n}}(\widetilde t_{i_1\cdots i_{k+1}})=-\sum_{\sigma\in S_{k+1}}\frac{sgn(\sigma)}{(k+1)!}\cdot\Delta^{\sigma}_{u_1\cdots u_j,v_1v_2\cdots v_n}(\widetilde t_{i_1\cdots i_{k+1}}).
\end{split}
\end{equation}
We complete the proof of Claim I. \,\,\,\,$\endpf$
\medskip

In order to compute $I_4$, we introduce the following notation.

Fix $\sigma\in S_{k+1}$. For integers $j, l$ such that $j\geq1$ and $0\leq l\leq k-j$, and distinct integers $u_1,\cdots,u_j,x_1,\cdots,x_l\in\{1,\cdots,k\}$ (by convention when $l=0$, there is no $x$), we define 
\begin{equation}\label{;empty}
\begin{split}
&\Delta^{\sigma}_{(u_1\cdots u_j;\,x_1\cdots x_{l})}(\widetilde t_{i_1\cdots i_{k+1}})=tr\big\{\prod_{l=1}^k\big(g_{i_{\sigma(l)} i_{\sigma(u_1)}}^{-1}dg_{i_{\sigma(l)} i_{\sigma(u_1)}}\big)^{1-\delta_{x_1\cdots x_l}(l)-\delta_{u_1\cdots u_j}(l)}\\
&\,\,\,\,\,\,\,\cdot\big(-g_{i_{\sigma(u_1)} i_{\sigma(k+1)}}dg_{i_{\sigma(u_1)} i_{\sigma(k+1)}}^{-1}\big)^{\delta_{x_1\cdots x_l}(l)}\cdot\big(g_{i_{\sigma(l)}i_{\sigma(u_1)}}^{-1}h_{i_{\sigma(l)}}dh_{i_{\sigma(l)}}^{-1}g_{i_{\sigma(l)}i_{\sigma(u_1)}}\big)^{\delta_{u_1\cdots u_j}(l)}\big\}.
\end{split}
\end{equation}
We also denote $\Delta^{\sigma}_{(u_1\cdots u_j;\,x_1\cdots x_{l})}(\widetilde t_{i_1\cdots i_{k+1}})$ by $\Delta^{\sigma}_{(u_1\cdots u_j;\,\emptyset)}(\widetilde t_{i_1\cdots i_{k+1}})$ when $l=0$. 
\medskip

\noindent{\bf Claim II : } 
\begin{equation}\label{I_4}
\begin{split}
I_4=\sum_{\sigma\in S_{k+1}}&\frac{sgn(\sigma)}{(k+1)!}\sum_{\substack{1\leq j\leq k\,; \\  1\leq u_1<u_2<\cdots<u_j\leq k;}}\Delta^{\sigma}_{(u_1\cdots u_j;\,\emptyset)}(\widetilde t_{i_1\cdots i_{k+1}})\\
&+\sum_{\sigma\in S_{k+1}}\frac{sgn(\sigma)}{(k+1)!}\sum_{\substack{1\leq j\leq k-1\,; \\  1\leq u_1<u_2<\cdots<u_j\leq k;\\   1\leq x_1\leq k,\, x_1\notin\{u_1,\cdots,u_j\}}}\Delta^{\sigma}_{(u_1\cdots u_j;\,x_1)}(\widetilde t_{i_1\cdots i_{k+1}}).\\
\end{split}
\end{equation}


\noindent{\bf Proof of Claim II :}
Recall the following identities  
\begin{equation}\label{tran}
g_{\gamma\beta}^{-1}h_{\gamma}^{-1}dh_{\gamma}g_{\gamma\beta}=g_{\alpha\beta}^{-1}g_{\gamma\alpha}^{-1}h_{\gamma}^{-1}dh_{\gamma}g_{\gamma\alpha}g_{\alpha\beta}\,\,\,\,\,{\rm and}\,\,\,g_{\delta\beta}^{-1}dg_{\delta\beta}=g_{\alpha\beta}^{-1}g_{\delta\alpha}^{-1}dg_{\delta\alpha}g_{\alpha\beta}+g_{\alpha\beta}^{-1}dg_{\alpha\beta}.
\end{equation}
Applying formula (\ref{tran}), we  have
\begin{equation}\label{j+0}
\begin{split}
\Delta^{\sigma}_{u_1\cdots u_j,\emptyset}&(\widetilde t_{i_1\cdots i_{k+1}})=tr\big\{\prod_{l=1}^k\big(g_{i_{\sigma(l)} i_{\sigma(k+1)}}^{-1}dg_{i_{\sigma(l)} i_{\sigma(k+1)}}\big)^{1-\delta_{u_1\cdots u_j}(l)}\cdot\big(g_{i_{\sigma(l)}i_{\sigma(k+1)}}^{-1}h_{i_{\sigma(l)}}dh_{i_{\sigma(l)}}^{-1}g_{i_{\sigma(l)}i_{\sigma(k+1)}}\big)^{\delta_{u_1\cdots u_j}(l)}\big\}\\
&=tr\big\{\prod_{l=1}^k\big(g_{i_{\sigma(u_1)}i_{\sigma(k+1)}}^{-1}g_{i_{\sigma(l)} i_{\sigma(u_1)}}^{-1}dg_{i_{\sigma(l)} i_{\sigma(u_1)}}g_{i_{\sigma(u_1)}i_{\sigma(k+1)}}+g_{i_{\sigma(u_1)} i_{\sigma(k+1)}}^{-1}dg_{i_{\sigma(u_1)} i_{\sigma(k+1)}}\big)^{1-\delta_{u_1\cdots u_j}(l)}\\
&\,\,\,\,\,\,\,\,\,\,\,\,\,\,\,\,\,\,\,\,\,\,\cdot\big(g_{i_{\sigma(u_1)}i_{\sigma(k+1)}}^{-1}g_{i_{\sigma(l)}i_{\sigma(u_1)}}^{-1}h_{i_{\sigma(l)}}dh_{i_{\sigma(l)}}^{-1}g_{i_{\sigma(l)}i_{\sigma(u_1)}}g_{i_{\sigma(u_1)}i_{\sigma(k+1)}}\big)^{\delta_{u_1\cdots u_j}(l)}\big\}\\
&=tr\big\{\prod_{l=1}^k\big(g_{i_{\sigma(l)} i_{\sigma(u_1)}}^{-1}dg_{i_{\sigma(l)} i_{\sigma(u_1)}}+(-1)\cdot g_{i_{\sigma(u_1)} i_{\sigma(k+1)}}dg_{i_{\sigma(u_1)} i_{\sigma(k+1)}}^{-1}\big)^{1-\delta_{u_1\cdots u_j}(l)}\\
&\,\,\,\,\,\,\,\,\,\,\,\,\,\,\,\,\,\,\,\,\,\,\cdot\big(g_{i_{\sigma(l)}i_{\sigma(u_1)}}^{-1}h_{i_{\sigma(l)}}dh_{i_{\sigma(l)}}^{-1}g_{i_{\sigma(l)}i_{\sigma(u_1)}}\big)^{\delta_{u_1\cdots u_j}(l)}\big\}\\
&=\sum_{\substack{0\leq l\leq k-j\,; \\  1\leq x_1<\cdots<x_l\leq k\,;\\ u_1,\cdots,u_j,x_1,\cdots,x_l\, {\text{are distinct}}}}\Delta^{\sigma}_{(u_1\cdots u_j;\,x_1\cdots x_{l})}(\widetilde t_{i_1\cdots i_{k+1}}).
\end{split}
\end{equation}
Then,
\begin{equation}
\begin{split}
I_4&=\sum_{\sigma\in S_{k+1}}\frac{sgn(\sigma)}{(k+1)!}\sum_{\substack{j\geq 1; \\  1\leq u_1<u_2<\cdots<u_j\leq k;\,}}\sum_{\substack{0\leq l\leq k-j\,; \\  1\leq x_1<\cdots<x_l\leq k\,;\\ u_1,\cdots,u_j,x_1,\cdots,x_l\, {\text{are distinct}}}}\Delta^{\sigma}_{(u_1\cdots u_j;\,x_1\cdots x_{l})}(\widetilde t_{i_1\cdots i_{k+1}})\\
&=\sum_{\sigma\in S_{k+1}}\frac{sgn(\sigma)}{(k+1)!}\sum_{\substack{ 1\leq j\leq k\,; \\  1\leq u_1<u_2<\cdots<u_j\leq k\,}}\Delta^{\sigma}_{(u_1\cdots u_j;\,\emptyset)}(\widetilde t_{i_1\cdots i_{k+1}})\\
&+\sum_{\sigma\in S_{k+1}}\frac{sgn(\sigma)}{(k+1)!}\sum_{\substack{1\leq j\leq k-1\,; \\  1\leq u_1<u_2<\cdots<u_j\leq k;\\   1\leq x_1\leq k,\,x_1\notin\{u_1,\cdots,u_j\}}}\Delta^{\sigma}_{(u_1\cdots u_j;\,x_1)}(\widetilde t_{i_1\cdots i_{k+1}})\\
&+\sum_{\sigma\in S_{k+1}}\frac{sgn(\sigma)}{(k+1)!}\sum_{\substack{1\leq j\leq k-2\,; \\  1\leq u_1<u_2<\cdots<u_j\leq k\,}}\sum_{\substack{2\leq l\leq k-j\,; \\  1\leq x_1<\cdots<x_l\leq k\,;\\ u_1,\cdots,u_j,x_1,\cdots,x_l\, {\text{are distinct}}}}\Delta^{\sigma}_{(u_1\cdots u_j;\,x_1\cdots x_{l})}(\widetilde t_{i_1\cdots i_{k+1}})\\
&=:J_1+J_2+J_3.
\end{split}
\end{equation}
Interchanging the order of summation, we have that
\begin{equation}
\begin{split}
J_3&=\sum_{\substack{1\leq j\leq k-2\,; \\  1\leq u_1<u_2<\cdots<u_j\leq k\,}}\sum_{\substack{2\leq l\leq k-j\,; \\  1\leq x_1<\cdots<x_l\leq k\,;\\ u_1,\cdots,u_j,x_1,\cdots,x_l\, {\text{are distinct}}}}\sum_{\sigma\in S_{k+1}}\frac{sgn(\sigma)}{(k+1)!}\cdot\Delta^{\sigma}_{(u_1\cdots u_j;\,x_1\cdots x_{l})}(\widetilde t_{i_1\cdots i_{k+1}})\\
\end{split}
\end{equation}
Hence, it suffices to prove that for integers $j, l$ such that $1\leq j$, $2\leq l\leq k-j$, and distinct integers $u_1,\cdots,u_j,x_1,\cdots,x_l\in\{1,\cdots,k\}$, the following equality holds,
\begin{equation}\label{j3}
\sum_{\sigma\in S_{k+1}}sgn(\sigma)\cdot\Delta^{\sigma}_{(u_1\cdots u_j;\,x_1\cdots x_{l})}(\widetilde t_{i_1\cdots i_{k+1}})\equiv 0.
\end{equation}

We shall prove equality (\ref{j3}) in the same way as the proof of Claim I of Lemma \ref{cocycle}.   Denote by $\tau$ the permutation $(x_1,x_2)\in S_{k+1}$. Then,
\begin{equation}
\begin{split}
&\Delta^{\sigma\circ\tau}_{(u_1\cdots u_j;\,x_1\cdots x_{l})}(\widetilde t_{i_1\cdots i_{k+1}})=tr\big\{\prod_{l=1}^k\big(g_{i_{(\sigma\circ\tau)(l)} i_{(\sigma\circ\tau)(u_1)}}^{-1}dg_{i_{(\sigma\circ\tau)(l)} i_{(\sigma\circ\tau)(u_1)}}\big)^{1-\delta_{x_1\cdots x_{l}}(l)-\delta_{u_1\cdots u_j}(l)}\\
&\,\,\,\,\,\,\,\,\,\,\,\,\,\,\,\,\,\,\,\,\cdot\big(-g_{i_{(\sigma\circ\tau)(u_1)} i_{(\sigma\circ\tau)(k+1)}}dg_{i_{(\sigma\circ\tau)(u_1)} i_{(\sigma\circ\tau)(k+1)}}^{-1}\big)^{\delta_{x_1\cdots x_{l}}(l)}\cdot\\
&\,\,\,\,\,\,\,\,\,\,\,\,\,\,\,\,\,\,\,\,\cdot\big(g_{i_{(\sigma\circ\tau)(l)}i_{(\sigma\circ\tau)(u_1)}}^{-1}h_{i_{(\sigma\circ\tau)(l)}}dh_{i_{(\sigma\circ\tau)(l)}}^{-1}g_{i_{(\sigma\circ\tau)(l)}i_{(\sigma\circ\tau)(u_1)}}\big)^{\delta_{u_1\cdots u_j}(l)}\big\}\\
&\,\,\,\,=tr\big\{\prod_{l=1}^k\big(g_{i_{\sigma(l)} i_{\sigma(u_1)}}^{-1}dg_{i_{\sigma(l)} i_{\sigma(u_1)}}\big)^{1-\delta_{x_1\cdots x_l}(l)-\delta_{u_1\cdots u_j}(l)}\cdot\\
&\,\,\,\,\,\,\,\,\,\,\cdot\big(-g_{i_{\sigma(u_1)} i_{\sigma(k+1)}}dg_{i_{\sigma(u_1)} i_{\sigma(k+1)}}^{-1}\big)^{\delta_{x_1\cdots x_l}(l)}\big(g_{i_{\sigma(l)}i_{\sigma(u_1)}}^{-1}h_{i_{\sigma(l)}}dh_{i_{\sigma(l)}}^{-1}g_{i_{\sigma(l)}i_{\sigma(u_1)}}\big)^{\delta_{u_1\cdots u_j}(l)}\big\}\\
&\,\,\,\,=\Delta^{\sigma}_{(u_1\cdots u_j;\,x_1\cdots x_{l})}(\widetilde t_{i_1\cdots i_{k+1}}).
\end{split}
\end{equation}
Therefore, we have that
\begin{equation}\label{t3j}
\begin{split}
\sum_{\sigma\in S_{k+1}}&sgn(\sigma)\cdot\Delta^{\sigma}_{(u_1\cdots u_j;\,x_1\cdots x_{l})}(\widetilde t_{i_1\cdots i_{k+1}})=\sum_{\sigma\circ\tau\in S_{k+1}}sgn(\sigma\circ\tau)\cdot\Delta^{\sigma\circ\tau}_{(u_1\cdots u_j;\,x_1\cdots x_{l})}(\widetilde t_{i_1\cdots i_{k+1}})\\
&=sgn(\tau)\cdot\sum_{\sigma\in S_{k+1}}sgn(\sigma)\cdot\Delta^{\sigma\circ\tau}_{(u_1\cdots u_j;\,x_1\cdots x_{l})}(\widetilde t_{i_1\cdots i_{k+1}})\\
&=-\sum_{\sigma\in S_{k+1}}sgn(\sigma)\cdot\Delta^{\sigma}_{(u_1\cdots u_j;\,x_1\cdots x_{l})}(\widetilde t_{i_1\cdots i_{k+1}}).\\
\end{split}
\end{equation}
We complete the proof of Claim II.\,\,\,\,$\endpf$
\medskip

We now construct a \v Cech $(k-1)$-cochain $s$ as follows. For $i_1,i_2,\cdots,i_{k+1}\in I$ and $1\leq\alpha\leq k+1$, define 
\begin{equation}\label{s0}
\begin{split}
s_{i_1\cdots\widehat{i_\alpha}\cdots i_{k+1};\,i_1\cdots i_{k+1}}:&=\frac{(-1)^{\alpha-1}}{(k+1)!} \sum_{\substack{1\leq j\leq k-1\,; \\  1\leq u_1<\cdots<u_j\leq k\,;\\ 1\leq v_1\leq k\,;\\u_1,\cdots,u_j,v_1\, {\text{are distinct}}}}\sum_{\substack{\sigma\in S_{k+1},\\\sigma(v_1)=\alpha}}sgn(\sigma)\cdot\Delta^{\sigma}_{u_1\cdots u_j,v_{1}}(\widetilde t_{i_1\cdots i_{k+1}})\\
&+\frac{(-1)^{\alpha-1}}{(k+1)!} \sum_{\substack{ 1\leq v_1\leq k}}\,\,\,\sum_{\substack{\sigma\in S_{k+1},\\\sigma(v_1)=\alpha}}{sgn(\sigma)}\cdot\Delta^{\sigma}_{\emptyset,v_{1}}(\widetilde t_{i_1\cdots i_{k+1}})\\
&+\frac{(-1)^{\alpha-1}}{(k+1)!}\sum_{\substack{1\leq j\leq k\,; \\  1\leq u_1<u_2<\cdots<u_j\leq k;}}\sum_{\substack{\sigma\in S_{k+1},\\\sigma(k+1)=\alpha}}{sgn(\sigma)}\cdot\Delta^{\sigma}_{(u_1\cdots u_j;\,\emptyset)}(\widetilde t_{i_1\cdots i_{k+1}})\\
&+\frac{(-1)^{\alpha-1}}{(k+1)!}\sum_{\substack{1\leq j\leq k-1\,; \\  1\leq u_1<u_2<\cdots<u_j\leq k;\\   1\leq x_1\leq k,\, x_1\notin\{u_1,\cdots,u_j\}}}\sum_{\substack{\sigma\in S_{k+1},\\\sigma(x_1)=\alpha}}{sgn(\sigma)}\cdot\Delta^{\sigma}_{(u_1\cdots u_j;\,x_1)}(\widetilde t_{i_1\cdots i_{k+1}}).\\
\end{split}
\end{equation}

\noindent{\bf Claim III : }  The above definition depends on   $(i_1,\cdots,\widehat{i_{\alpha}},\cdots,i_{k+1})$ but  not  on $(i_1,\cdots,i_{k+1})$; namely, if $(i_1,\cdots,\widehat{i_{\alpha}},\cdots,i_{k+1})=(j_1,\cdots,\widehat {j_{\beta}},\cdots,j_{k+1})$,   $s_{i_1\cdots\widehat{i_\alpha}\cdots i_{k+1};\,i_1\cdots i_{k+1}}=s_{j_1\cdots\widehat{j_\beta}\cdots j_{k+1};\,j_1\cdots j_{k+1}}$. 
\smallskip

\noindent{\bf Proof of Claim III :} 
For any integers $\alpha$, $\beta$ and $i_1,\cdots, i_{k+2}$ such that $1\leq\alpha<\beta\leq k+2$ and  $i_1,i_2,\cdots,i_{k+1}\in I$, we can define 
$s_{i_1\cdots\widehat{i_\alpha}\cdots\widehat{i_\beta}\cdots i_{k+1};\,i_1\cdots\widehat{i_\beta}\cdots i_{k+2}}$  and
$s_{i_1\cdots\widehat{i_\alpha}\cdots\widehat{i_\beta}\cdots i_{k+1};\,i_1\cdots\widehat {i_{\alpha}}\cdots i_{k+2}}$ by formula (\ref{s0}).  Firstly, we have
\begin{equation}\label{s1}
\begin{split}
s_{i_1\cdots\widehat{i_\alpha}\cdots\widehat{i_\beta}\cdots i_{k+1};\,i_1\cdots\widehat{i_\beta}\cdots i_{k+2}}&=\frac{(-1)^{\alpha-1}}{(k+1)!} \sum_{\substack{1\leq j\leq k-1\,; \\  1\leq u_1<\cdots<u_j\leq k\,;\\ 1\leq v_1\leq k\,;\\u_1,\cdots,u_j,v_1\, {\text{are distinct}}}}\sum_{\substack{\sigma\in S_{k+1},\\\sigma(v_1)=\alpha}}sgn(\sigma)\cdot\Delta^{\sigma}_{u_1\cdots u_j,v_{1}}(\widetilde t_{i_1\cdots\widehat{i_\beta}\cdots i_{k+2}})\\
&+\frac{(-1)^{\alpha-1}}{(k+1)!} \sum_{\substack{j=0\,; \\ 1\leq v_1\leq k\,;}}\sum_{\substack{\sigma\in S_{k+1},\\\sigma(v_1)=\alpha}}{sgn(\sigma)}\cdot\Delta^{\sigma}_{\emptyset,v_{1}}(\widetilde t_{i_1\cdots\widehat{i_\beta}\cdots i_{k+2}})\\
&+\frac{(-1)^{\alpha-1}}{(k+1)!}\sum_{\substack{1\leq j\leq k\,; \\  1\leq u_1<u_2<\cdots<u_j\leq k;}}\sum_{\substack{\sigma\in S_{k+1},\\\sigma(k+1)=\alpha}}{sgn(\sigma)}\cdot\Delta^{\sigma}_{(u_1\cdots u_j;\,\emptyset)}(\widetilde t_{i_1\cdots\widehat{i_\beta}\cdots i_{k+2}})\\
&+\frac{(-1)^{\alpha-1}}{(k+1)!}\sum_{\substack{1\leq j\leq k-1\,; \\  1\leq u_1<u_2<\cdots<u_j\leq k;\\   1\leq x_1\leq k,\, x_1\notin\{u_1,\cdots,u_j\}}}\sum_{\substack{\sigma\in S_{k+1},\\\sigma(x_1)=\alpha}}{sgn(\sigma)}\cdot\Delta^{\sigma}_{(u_1\cdots u_j;\,x_1)}(\widetilde t_{i_1\cdots\widehat{i_\beta}\cdots i_{k+2}})\\
&=:\kappa^{\alpha}_1+\kappa^{\alpha}_2+\kappa^{\alpha}_3+\kappa^{\alpha}_4.
\end{split}
\end{equation}

Notice that in the definition of $s_{i_1\cdots\widehat{i_\alpha}\cdots i_{k+1};\,i_1\cdots i_{k+1}}$, the number $\alpha$ on the right hand side is referring to the position of the omitted index on the left hand side. Therefore, 
$s_{i_1\cdots\widehat{i_\alpha}\cdots\widehat{i_\beta}\cdots i_{k+1};\,i_1\cdots\widehat{i_\alpha}\cdots i_{k+2}}$  takes the following form,
\begin{equation}\label{s2}
\begin{split}
s_{i_1\cdots\widehat{i_\alpha}\cdots\widehat{i_\beta}\cdots i_{k+1};\,i_1\cdots\widehat{i_{\alpha}}\cdots i_{k+2}}&=\frac{(-1)^{\beta-2}}{(k+1)!} \sum_{\substack{1\leq j\leq k-1\,; \\  1\leq u_1<\cdots<u_j\leq k\,;\\ 1\leq v_1\leq k\,;\\u_1,\cdots,u_j,v_1\, {\text{are distinct}}}}\sum_{\substack{\tau\in S_{k+1},\\\tau(v_1)=\beta-1}}sgn(\tau)\cdot\Delta^{\tau}_{u_1\cdots u_j,v_{1}}(\widetilde t_{i_1\cdots\widehat{i_\alpha}\cdots i_{k+2}})\\
&+\frac{(-1)^{\beta-2}}{(k+1)!} \sum_{\substack{j=0\,; \\ 1\leq v_1\leq k\,;}}\sum_{\substack{\tau\in S_{k+1},\\\tau(v_1)=\beta-1}}{sgn(\tau)}\cdot\Delta^{\tau}_{\emptyset,v_{1}}(\widetilde t_{i_1\cdots\widehat{i_\alpha}\cdots i_{k+2}})\\
&+\frac{(-1)^{\beta-2}}{(k+1)!}\sum_{\substack{1\leq j\leq k\,; \\  1\leq u_1<u_2<\cdots<u_j\leq k;}}\sum_{\substack{\sigma\in S_{k+1},\\\sigma(k+1)=\beta-1}}{sgn(\sigma)}\cdot\Delta^{\sigma}_{(u_1\cdots u_j;\,\emptyset)}(\widetilde t_{i_1\cdots\widehat{i_\alpha}\cdots i_{k+2}})\\
&+\frac{(-1)^{\beta-2}}{(k+1)!}\sum_{\substack{1\leq j\leq k-1\,; \\  1\leq u_1<u_2<\cdots<u_j\leq k;\\   1\leq x_1\leq k,\, x_1\notin\{u_1,\cdots,u_j\}}}\sum_{\substack{\tau\in S_{k+1},\\\tau(x_1)=\beta-1}}{sgn(\tau)}\cdot\Delta^{\tau}_{(u_1\cdots u_j;\,x_1)}(\widetilde t_{i_1\cdots\widehat{i_\alpha}\cdots i_{k+2}})\\
&=:
\kappa^{\beta}_1+\kappa^{\beta}_2+\kappa^{\beta}_3+\kappa^{\beta}_4.
\end{split}
\end{equation}

In order to prove Claim III, it suffices to prove that  
\begin{equation}\label{s3}
s_{i_1\cdots\widehat{i_\alpha}\cdots\widehat{i_\beta}\cdots i_{k+1};\,i_1\cdots\widehat{i_\beta}\cdots i_{k+2}}=s_{i_1\cdots\widehat{i_\alpha}\cdots\widehat{i_\beta}\cdots i_{k+1};\,i_1\cdots\widehat {i_{\alpha}}\cdots i_{k+2}}.
\end{equation}
By  Claim IV to be proved in the following,  formula (\ref{s3}) holds. We complete the proof of Claim III.
\,\,\,\,$\endpf$
\medskip

Next we will prove  Claim IV used above.
\medskip

\noindent{\bf Claim IV :}
\begin{equation}
\kappa^{\alpha}_i=(-1)^{\alpha-\beta+1}\cdot\kappa^{\beta}_i\,\,\,\,\,{\rm for}\,\,\,\, i=1,2,3,4.
\end{equation}
\medskip

\noindent{\bf Proof of Claim IV :}  We will prove Claim IV based on a case by case argument.
\medskip

\noindent{\bf Case I ($\kappa^{\alpha}_1=(-1)^{\alpha-\beta+1}\cdot\kappa^{\beta}_1$) :} It suffices to prove that for fixed integer $j$ such that $1\leq j\leq k-1$, and distinct integers $u_1,\cdots,u_j,v_1\in\{1,\cdots,k\}$, the following equality holds,
\begin{equation}\label{k1}
\sum_{\substack{\sigma\in S_{k+1},\\\sigma(v_1)=\alpha}}sgn(\sigma)\cdot\Delta^{\sigma}_{u_1\cdots u_j,v_{1}}(\widetilde t_{i_1\cdots\widehat{i_\beta}\cdots i_{k+2}})=(-1)^{\alpha-\beta+1}\cdot\sum_{\substack{\tau\in S_{k+1},\\\tau(v_1)=\beta-1}}sgn(\tau)\cdot\Delta^{\tau}_{u_1\cdots u_j,v_{1}}(\widetilde t_{i_1\cdots\widehat{i_\alpha}\cdots i_{k+2}}).
\end{equation}

In order to compare the terms on both sides, we lift each element of $S_{k+1}$ to  an element of $S_{k+2}$ as follows. Denote by $S_{k+1}^{\beta}$ the set consisting of all bijections from $\{1,2,\cdots,k+1\}$ to $\{1,\cdots,\widehat\beta, \cdots,k+2\}$; denote by $S_{k+1}^{\alpha}$ the set consisting of all bijections from $\{1,2,\cdots,k+1\}$ to $\{1,\cdots,\widehat\alpha, \cdots,k+2\}$.  Define a bijection 
\begin{equation}
\begin{split}
\widehat\bullet\,\,\,\,\,:\,\,\,\,\,\,S_{k+1}&\longrightarrow S_{k+1}^{\beta} \\
\sigma&\longmapsto\widehat{\sigma}
\end{split}
\end{equation}
 by
\begin{equation}\label{hlift1}
\widehat\sigma(l)=\left\{\begin{matrix}
\sigma(l)&\,\,\,\,\,\,\,\,\,\,\,{\rm if}\,\,\,\,\,\,\,\,\,\,\,1\leq l\leq k+1\,\,\,\,{\rm and}\,\,\,\,\sigma(l)< \beta,\\\sigma(l)+1&\,\,\,\,\,\,\,\,\,\,\,{\rm if}\,\,\,\,\,\,\,\,\,\,\,1\leq l\leq k+1\,\,\,\,{\rm and}\,\,\,\,\sigma(l)\geq \beta.\\
\end{matrix}\right.
\end{equation}
Moreover, to compare the signature,  we define a bijection 
\begin{equation}
\begin{split}
[\bullet]_{\beta}\,\,\,\,\,:\,\,\,\,\,\,S^{\beta}_{k+1}&\longrightarrow \big\{\eta\big|\eta\in S_{k+2},\eta(k+2)=\beta\big\} \\
\widehat\sigma&\longmapsto[\widehat\sigma]_{\beta}
\end{split}
\end{equation}
by
\begin{equation}\label{lift1}
[\widehat\sigma]_{\beta}(l)=\left\{\begin{matrix}
\widehat\sigma(l)&\,\,\,\,\,\,\,\,\,{\rm if}\,\,\,\,\,\,\,\,\,\,\,\,\,\,\,\,\,\,\,\,\,\,\,\,\,\,1\leq l\leq k+1,\\\beta&\,\,\,\,\,\,\,\,\,\,\,\,\,\,\,\,\,\,\,\,\,\,\,\,\,\,\,\,\,\,\,\,\,\,\,\,\,\,\,\,\,\,{\rm if}\,\,\,\,\,\,\,\,\,\,\,\,\,\,\,\,\,\,\,\,\,\,\,\,\,\,\,\,\,l=k+2\,;\,\,\,\,\,\,\,\,\,\,\,\,\,\,\,\,\,\,\,\,\,\,\,\,\,\,\,\,\,\,\,\,\,\,\,\,\,\,\,\,\\
\end{matrix}\right.
\end{equation}
then, the following equality holds,
\begin{equation}
sgn(\sigma)=(-1)^{k+2-\beta}\cdot sgn([\widehat\sigma]_{\beta}).
\end{equation}
Similarly, we can define a bijection 
\begin{equation}
\begin{split}
\overline\bullet\,\,\,\,\,:\,\,\,\,\,\,S_{k+1}&\longrightarrow S_{k+1}^{\alpha} \\
\tau&\longmapsto\overline{\tau}
\end{split}
\end{equation}
by
\begin{equation}\label{hlift2}
\overline\tau(l)=\left\{\begin{matrix}
\tau(l)&\,\,\,\,\,\,\,\,\,\,\,{\rm if}\,\,\,\,\,\,\,\,\,\,\,1\leq l\leq k+1\,\,\,\,{\rm and}\,\,\,\,\tau(l)< \alpha,\\\tau(l)+1&\,\,\,\,\,\,\,\,\,\,\,{\rm if}\,\,\,\,\,\,\,\,\,\,\,1\leq l\leq k+1\,\,\,\,{\rm and}\,\,\,\,\tau(l)\geq \alpha;\\
\end{matrix}\right.
\end{equation}
define a bijection 
\begin{equation}
\begin{split}
[\bullet]_{\alpha}\,\,\,\,\,:\,\,\,\,\,\,S^{\alpha}_{k+1}&\longrightarrow \big\{\eta\big|\eta\in S_{k+2},\eta(k+2)=\alpha\big\} \\
\overline\tau&\longmapsto[\overline\tau]_{\alpha}
\end{split}
\end{equation}
by
\begin{equation}\label{lift2}
[\overline\tau]_{\alpha}(l)=\left\{\begin{matrix}
\overline\tau(l)&\,\,\,\,\,\,\,\,\,{\rm if}\,\,\,\,\,\,\,\,\,\,\,\,\,\,\,\,\,\,\,\,\,\,\,\,\,\,1\leq l\leq k+1,\\\alpha&\,\,\,\,\,\,\,\,\,\,\,\,\,\,\,\,\,\,\,\,\,\,\,\,\,\,\,\,\,\,\,\,\,\,\,\,\,\,\,\,\,\,\,\,{\rm if}\,\,\,\,\,\,\,\,\,\,\,\,\,\,\,\,\,\,\,\,\,\,\,\,\,\,\,\,\,\,\,l=k+2\;,\,\,\,\,\,\,\,\,\,\,\,\,\,\,\,\,\,\,\,\,\,\,\,\,\,\,\,\,\,\,\,\,\,\,\,\,\,\,\,\,\\
\end{matrix}\right.
\end{equation}
the following  equality holds,
\begin{equation}
sgn(\tau)=(-1)^{k+2-\alpha}\cdot sgn([\overline\tau]_{\alpha}).
\end{equation} 

Rewriting formula $(\ref{dis})$ under the above notation, we have 
\begin{equation}\label{dis2}
\begin{split}
\Delta^{\sigma}_{u_1\cdots u_j,v_{1}}&(\widetilde t_{i_1\cdots\widehat{i_\beta}\cdots i_{k+2}})=tr\big\{\prod_{l=1}^k\big(g_{i_{\widehat\sigma(l)} i_{\widehat\sigma(k+1)}}^{-1}dg_{i_{\widehat\sigma(l)} i_{\widehat\sigma(k+1)}}\big)^{1-\delta_{u_1\cdots u_j}(l)-\delta_{v_1}(l)}\\
&\cdot\big(g_{i_{\widehat\sigma(l)}i_{\widehat\sigma(k+1)}}^{-1}h_{i_{\widehat\sigma(l)}}dh_{i_{\widehat\sigma(l)}}^{-1}g_{i_{\widehat\sigma(l)}i_{\widehat\sigma(k+1)}}\big)^{\delta_{u_1\cdots u_j}(l)}\cdot\big(-h_{i_{\widehat\sigma(k+1)}}dh_{i_{\widehat\sigma(k+1)}}^{-1}\big)^{\delta_{v_1}(l)}\big\},
\end{split}
\end{equation}
and
\begin{equation}\label{dis3}
\begin{split}
\Delta^{\tau}_{u_1\cdots u_j,v_{1}}&(\widetilde t_{i_1\cdots\widehat{i_\alpha}\cdots i_{k+2}})=tr\big\{\prod_{l=1}^k\big(g_{i_{\overline\tau(l)} i_{\overline\tau(k+1)}}^{-1}dg_{i_{\overline\tau(l)} i_{\overline\tau(k+1)}}\big)^{1-\delta_{u_1\cdots u_j}(l)-\delta_{v_1}(l)}\\
&\cdot\big(g_{i_{\overline\tau(l)}i_{\overline\tau(k+1)}}^{-1}h_{i_{\overline\tau(l)}}dh_{i_{\overline\tau(l)}}^{-1}g_{i_{\overline\tau(l)}i_{\overline\tau(k+1)}}\big)^{\delta_{u_1\cdots u_j}(l)}\cdot\big(-h_{i_{\overline\tau(k+1)}}dh_{i_{\overline\tau(k+1)}}^{-1}\big)^{\delta_{v_1}(l)}\big\}.
\end{split}
\end{equation}

In order to prove equality $(\ref{k1})$, it suffices to construct a bijection
\begin{equation}
\begin{split}
P_{\alpha,\beta,u_1,\cdots,u_j,v_1}:\{\sigma\in S_{k+1}\big|\sigma(v_1)=\alpha\}\longrightarrow\{\tau\in S_{k+1}\big|\tau(v_1)=\beta-1\},
\end{split}
\end{equation}
such that 
\begin{equation}\label{k11}
\begin{split}
&sgn(\sigma)\cdot\Delta^{\sigma}_{u_1\cdots u_j,v_{1}}(\widetilde t_{i_1\cdots\widehat{i_\beta}\cdots i_{k+2}})=\\
&\,\,\,\,\,\,\,\,\,\,\,\,\,\,\,\,\,(-1)^{\alpha-\beta+1}\cdot sgn(P_{\alpha,\beta,u_1,\cdots,u_j,v_1}(\sigma))\cdot\Delta^{P_{\alpha,\beta,u_1,\cdots,u_j,v_1}(\sigma)}_{u_1\cdots u_j,v_{1}}(\widetilde t_{i_1\cdots\widehat{i_\alpha}\cdots i_{k+2}}).
\end{split}
\end{equation}
For each $\sigma\in S_{k+1}$ such that $\sigma(v_1)=\alpha$, we define $P_{\alpha,\beta,u_1,\cdots,u_j,v_1}({\sigma})$ as follows. Take $\widehat\sigma$ as in formula (\ref{hlift1}) and define $\eta\in S_{k+1}^{\alpha}$ by
\begin{equation}\label{hlift4}
\eta(l)=\left\{\begin{matrix}
\widehat\sigma(l)&\,\,\,\,\,\,\,\,\,\,\,\,\,\,\,\,\,\,\,\,\,\,\,\,\,{\rm if}\,\,\,\,\,\,\,\,\,\,\,1\leq l\leq v_1-1,\\\beta&\,\,\,\,\,\,\,\,\,\,\,{\rm if}\,\,\,\,\,\,\,\,\,\,\,\,\,\,\,\,\,\,l=v_1,\\
\widehat\sigma(l)&\,\,\,\,\,\,\,\,\,\,\,\,\,\,\,\,\,\,\,\,\,\,\,\,\,\,\,\,\,\,\,\,\,\,\,{\rm if}\,\,\,\,\,\,\,\,\,\,\,v_1+1\leq l\leq k+1;\\
\end{matrix}\right.
\end{equation}
let $P_{\alpha,\beta,u_1,\cdots,u_j,v_1}({\sigma})$ be  the inverse image of $\eta$ under the bijection $\overline{\bullet}$\,\,, that is, $\overline{P_{\alpha,\beta,u_1,\cdots,u_j,v_1}({\sigma})}=\eta$. It is easy to verify that $P_{\alpha,\beta,u_1,\cdots,u_j,v_1}({\sigma})\in \{\tau\in S_{k+1}\big|\tau(v_1)=\beta-1\}$. 

Then, by formulas $(\ref{dis2})$ and $(\ref{dis3})$, we have
\begin{equation}\label{dis4}
\begin{split}
&\Delta^{P_{\alpha,\beta,u_1,\cdots,u_j,v_1}({\sigma})}_{u_1\cdots u_j,v_{1}}(\widetilde t_{i_1\cdots\widehat{i_\alpha}\cdots i_{k+2}})=tr\big\{\prod_{l=1}^k\big(g_{i_{\eta(l)} i_{\eta(k+1)}}^{-1}dg_{i_{\eta(l)} i_{\eta(k+1)}}\big)^{1-\delta_{u_1\cdots u_j}(l)-\delta_{v_1}(l)}\\
&\,\,\,\,\,\,\,\,\,\,\,\,\,\,\,\,\,\,\,\,\,\,\,\,\,\,\,\,\cdot\big(g_{i_{\eta(l)}i_{\eta(k+1)}}^{-1}h_{i_{\eta(l)}}dh_{i_{\eta(l)}}^{-1}g_{i_{\eta(l)}i_{\eta(k+1)}}\big)^{\delta_{u_1\cdots u_j}(l)}\cdot\big(-h_{i_{\eta(k+1)}}dh_{i_{\eta(k+1)}}^{-1}\big)^{\delta_{v_1}(l)}\big\}\\
&=tr\big\{\prod_{l=1}^k\big(g_{i_{\widehat\sigma(l)} i_{\widehat\sigma(k+1)}}^{-1}dg_{i_{\widehat\sigma(l)} i_{\widehat\sigma(k+1)}}\big)^{1-\delta_{u_1\cdots u_j}(l)-\delta_{v_1}(l)}\cdot\big(g_{i_{\widehat\sigma(l)}i_{\widehat\sigma(k+1)}}^{-1}h_{i_{\widehat\sigma(l)}}dh_{i_{\widehat\sigma(l)}}^{-1}g_{i_{\widehat\sigma(l)}i_{\widehat\sigma(k+1)}}\big)^{\delta_{u_1\cdots u_j}(l)}\\
&\,\,\,\,\,\,\,\,\,\,\,\,\,\,\,\,\,\,\,\,\,\,\,\,\,\,\,\,\,\,\,\,\,\,\,\cdot\big(-h_{i_{\widehat\sigma(k+1)}}dh_{i_{\widehat\sigma(k+1)}}^{-1}\big)^{\delta_{v_1}(l)}\big\}=\Delta^{\sigma}_{u_1\cdots u_j,v_{1}}(\widetilde t_{i_1\cdots\widehat{i_\beta}\cdots i_{k+2}}).
\end{split}
\end{equation}

To complete the proof of Claim IV in this case, it suffices to prove the following equality,
\begin{equation}
sgn(P_{\alpha,\beta,u_1,\cdots,u_j,v_1}({\sigma}))=(-1)^{\alpha-\beta+1}\cdot sgn(\sigma).
\end{equation}
Since 
\begin{equation}
sgn(P_{\alpha,\beta,u_1,\cdots,u_j,v_1}({\sigma}))=(-1)^{k+2-\alpha}\cdot sgn([\overline{P_{\alpha,\beta,u_1,\cdots,u_j,v_1}({\sigma})}]_{\alpha})=(-1)^{k+2-\alpha}\cdot sgn([\eta]_{\alpha})
\end{equation}
and
\begin{equation}
(-1)^{\alpha-\beta+1}\cdot sgn(\sigma)=(-1)^{\alpha-\beta+1}\cdot(-1)^{k+2-\beta}\cdot sgn([\sigma]_{\beta}),
\end{equation}
it suffices to prove that
\begin{equation}\label{sgn}
sgn([\eta]_{\alpha})=-sgn([\sigma]_{\beta}).
\end{equation}
Since $[\eta]_{\alpha}=\iota\circ[\sigma]_{\beta}\in S_{k+2}$ where $\iota\in S_{k+2}$ is the permutation $(\alpha,\beta)$,  equality $(\ref{sgn})$ holds.

Therefore, $\kappa^{\alpha}_1=(-1)^{\alpha-\beta+1}\cdot\kappa^{\beta}_1$.
\medskip

\noindent{\bf Case II ($\kappa^{\alpha}_2=(-1)^{\alpha-\beta+1}\cdot\kappa^{\beta}_2$) :} Similarly, it suffices to prove that for fixed integer $v_1$ such that $1\leq v_1\leq k$ the following equality holds:
\begin{equation}\label{k2}
\sum_{\substack{\sigma\in S_{k+1},\\\sigma(v_1)=\alpha}}{sgn(\sigma)}\cdot\Delta^{\sigma}_{\emptyset,v_{1}}(\widetilde t_{i_1\cdots\widehat{i_\beta}\cdots i_{k+2}})=(-1)^{\alpha-\beta+1}\cdot\sum_{\substack{\tau\in S_{k+1},\\\tau(v_1)=\beta-1}}{sgn(\tau)}\cdot\Delta^{\tau}_{\emptyset,v_{1}}(\widetilde t_{i_1\cdots\widehat{i_\alpha}\cdots i_{k+2}}).
\end{equation}
Define $\widehat\sigma$, $[\widehat\sigma]_{\beta}$, $\overline\tau$ and $[\overline\tau]_{\alpha}$ in the same way as formulas (\ref{hlift1}), (\ref{lift1}), (\ref{hlift2}) and (\ref{lift2}).  Rewriting formula (\ref{dis2}) under the above notation, we have 
\begin{equation}\label{di2}
\begin{split}
\Delta^{\sigma}_{\emptyset,v_{1}}&(\widetilde t_{i_1\cdots\widehat{i_\beta}\cdots i_{k+2}})=tr\big\{\prod_{l=1}^k\big(g_{i_{\widehat\sigma(l)} i_{\widehat\sigma(k+1)}}^{-1}dg_{i_{\widehat\sigma(l)} i_{\widehat\sigma(k+1)}}\big)^{1-\delta_{v_1}(l)}\cdot\big(-h_{i_{\widehat\sigma(k+1)}}dh_{i_{\widehat\sigma(k+1)}}^{-1}\big)^{\delta_{v_1}(l)}\big\},
\end{split}
\end{equation}
and 
\begin{equation}\label{di3}
\begin{split}
\Delta^{\tau}_{\emptyset,v_{1}}&(\widetilde t_{i_1\cdots\widehat{i_\alpha}\cdots i_{k+2}})=tr\big\{\prod_{l=1}^k\big(g_{i_{\overline\tau(l)} i_{\overline\tau(k+1)}}^{-1}dg_{i_{\overline\tau(l)} i_{\overline\tau(k+1)}}\big)^{1-\delta_{v_1}(l)}\cdot\big(-h_{i_{\overline\tau(k+1)}}dh_{i_{\overline\tau(k+1)}}^{-1}\big)^{\delta_{v_1}(l)}\big\}.
\end{split}
\end{equation}

Notice that equality (\ref{k2}) holds if there exists a bijection
\begin{equation}
\begin{split}
P_{\alpha,\beta,v_1}:\{\sigma\in S_{k+1}\big|\sigma(v_1)=\alpha\}\rightarrow\{\tau\in S_{k+1}\big|\tau(v_1)=\beta-1\},
\end{split}
\end{equation}
such that 
\begin{equation}\label{k12}
sgn(\sigma)\cdot\Delta^{\sigma}_{\emptyset,v_{1}}(\widetilde t_{i_1\cdots\widehat{i_\beta}\cdots i_{k+2}})=(-1)^{\alpha-\beta+1}\cdot sgn(P_{\alpha,\beta,v_1}(\sigma))\cdot\Delta^{P_{\alpha,\beta,v_1}(\sigma)}_{\emptyset,v_{1}}(\widetilde t_{i_1\cdots\widehat{i_\alpha}\cdots i_{k+2}}).
\end{equation}
For each $\sigma\in S_{k+1}$ such that $\sigma(v_1)=\alpha$, we define an element $P_{\alpha,\beta,v_1}({\sigma})\in S_{k+1}$ as follows. Take $\widehat\sigma$ as in formula (\ref{hlift1}) and define $\eta\in S_{k+1}^{\alpha}$ by
\begin{equation}\label{hlift42}
\eta(l)=\left\{\begin{matrix}
\widehat\sigma(l)&\,\,\,\,\,\,\,\,\,\,\,{\rm if}\,\,\,\,\,\,\,\,\,\,\,1\leq l\leq v_1-1,\\\beta&\,\,\,\,\,\,\,\,\,\,\,{\rm if}\,\,\,\,\,\,\,\,\,\,\,l=v_1,\\
\widehat\sigma(l)&\,\,\,\,\,\,\,\,\,\,\,{\rm if}\,\,\,\,\,\,\,\,\,\,\,v_1+1\leq l\leq k+1;\\
\end{matrix}\right.
\end{equation}
let $P_{\alpha,\beta,v_1}({\sigma})$ be  the inverse image of $\eta$ under the map $\overline{\bullet}$\,\,, that is, $\overline{P_{\alpha,\beta,v_1}({\sigma})}=\eta$. It is easy to verify that $P_{\alpha,\beta,v_1}({\sigma})\in \{\tau\in S_{k+1}\big|\tau(v_1)=\beta-1\}$.

Then, by formulas $(\ref{dis2})$ and $(\ref{dis3})$, we have that
\begin{equation}\label{dis42}
\begin{split}
\Delta^{P_{\alpha,\beta,v_1}({\sigma})}_{\emptyset,v_{1}}(\widetilde t_{i_1\cdots\widehat{i_\alpha}\cdots i_{k+2}})&=tr\big\{\prod_{l=1}^k\big(g_{i_{\eta(l)} i_{\eta(k+1)}}^{-1}dg_{i_{\eta(l)} i_{\eta(k+1)}}\big)^{1-\delta_{v_1}(l)}\cdot\big((-1)\cdot h_{i_{\eta(k+1)}}dh_{i_{\eta(k+1)}}^{-1}\big)^{\delta_{v_1}(l)}\big\}\\
&=tr\big\{\prod_{l=1}^k\big(g_{i_{\widehat\sigma(l)} i_{\widehat\sigma(k+1)}}^{-1}dg_{i_{\widehat\sigma(l)} i_{\widehat\sigma(k+1)}}\big)^{1-\delta_{v_1}(l)}\cdot\big((-1)\cdot h_{i_{\widehat\sigma(k+1)}}dh_{i_{\widehat\sigma(k+1)}}^{-1}\big)^{\delta_{v_1}(l)}\big\}\\
&=\Delta^{\sigma}_{\emptyset,v_{1}}(\widetilde t_{i_1\cdots\widehat{i_\beta}\cdots i_{k+2}}).
\end{split}
\end{equation}

Notice that 
\begin{equation}
sgn(P_{\alpha,\beta,v_1}({\sigma}))=(-1)^{k+2-\alpha}\cdot sgn([\overline{P_{\alpha,\beta,v_1}({\sigma})}]_{\alpha})=(-1)^{k+2-\alpha}\cdot sgn([\eta]_{\alpha}),
\end{equation}
\begin{equation}
(-1)^{\alpha-\beta+1}\cdot sgn(\sigma)=(-1)^{\alpha-\beta+1}\cdot(-1)^{k+2-\beta}\cdot sgn([\widehat\sigma]_{\beta}).
\end{equation}
Since $[\eta]_{\alpha}=\iota\circ[\widehat\sigma]_{\beta}\in S_{k+2}$ where $\iota\in S_{k+2}$ is the permutation $(\alpha,\beta)$, $sgn([\sigma]_{\beta})=-sgn([\eta]_{\alpha})$.

Then, $\kappa^{\alpha}_2=(-1)^{\alpha-\beta+1}\cdot\kappa^{\beta}_2$.
\medskip

\noindent{\bf Case III ($\kappa^{\alpha}_3=(-1)^{\alpha-\beta+1}\cdot\kappa^{\beta}_3$) :} Similarly we will show that for fixed integers $j,u_1,\cdots,u_j$ such that $1\leq j\leq k$ and $1\leq u_1<\cdots<u_j\leq k$,  the following equality holds:
\begin{equation}\label{k3}
\sum_{\substack{\sigma\in S_{k+1},\\\sigma(k+1)=\alpha}}{sgn(\sigma)}\cdot\Delta^{\sigma}_{(u_1\cdots u_j;\,\emptyset)}(\widetilde t_{i_1\cdots\widehat{i_\beta}\cdots i_{k+2}})=(-1)^{\alpha-\beta+1}\cdot\sum_{\substack{\tau\in S_{k+1},\\\tau(k+1)=\beta-1}}{sgn(\tau)}\cdot\Delta^{\tau}_{(u_1\cdots u_j;\,\emptyset)}(\widetilde t_{i_1\cdots\widehat{i_\alpha}\cdots i_{k+2}}).
\end{equation}

Define $\widehat\sigma$, $[\widehat\sigma]_{\beta}$, $\overline\tau$ and $[\overline\tau]_{\alpha}$ in the same way as formulas (\ref{hlift1}), (\ref{lift1}), (\ref{hlift2}) and (\ref{lift2}). Rewriting formula (\ref{;empty}) under the new notation, we have 
\begin{equation}
\begin{split}
&\Delta^{\sigma}_{(u_1\cdots u_j;\,\emptyset)}(\widetilde t_{i_1\cdots \widehat i_{\beta}\cdots i_{k+1}})=tr\big\{\prod_{l=1}^k\big(g_{i_{\widehat\sigma(l)} i_{\widehat\sigma(u_1)}}^{-1}dg_{i_{\widehat\sigma(l)} i_{\widehat\sigma(u_1)}}\big)^{1-\delta_{u_1\cdots u_j}(l)}\cdot\big(g_{i_{\widehat\sigma(l)}i_{\widehat\sigma(u_1)}}^{-1}h_{i_{\widehat\sigma(l)}}dh_{i_{\widehat\sigma(l)}}^{-1}g_{i_{\widehat\sigma(l)}i_{\widehat\sigma(u_1)}}\big)^{\delta_{u_1\cdots u_j}(l)}\big\},
\end{split}
\end{equation}
and
\begin{equation}
\begin{split}
&\Delta^{\tau}_{(u_1\cdots u_j;\,\emptyset)}(\widetilde t_{i_1\cdots \widehat i_{\alpha}\cdots i_{k+1}})=tr\big\{\prod_{l=1}^k\big(g_{i_{\overline\tau(l)} i_{\overline\tau(u_1)}}^{-1}dg_{i_{\overline\tau(l)} i_{\overline\tau(u_1)}}\big)^{1-\delta_{u_1\cdots u_j}(l)}\cdot\big(g_{i_{\overline\tau(l)}i_{\overline\tau(u_1)}}^{-1}h_{i_{\overline\tau(l)}}dh_{i_{\overline\tau(l)}}^{-1}g_{i_{\overline\tau(l)}i_{\overline\tau(u_1)}}\big)^{\delta_{u_1\cdots u_j}(l)}\big\}.
\end{split}
\end{equation}

Notice that equality (\ref{k3}) holds if there exists a bijection
\begin{equation}
\begin{split}
P_{\alpha,\beta}:\{\sigma\in S_{k+1}\big|\sigma(k+1)=\alpha\}\rightarrow\{\tau\in S_{k+1}\big|\tau(k+1)=\beta-1\},
\end{split}
\end{equation}
such that 
\begin{equation}\label{k12}
{sgn(\sigma)}\cdot\Delta^{\sigma}_{(u_1\cdots u_j;\,\emptyset)}(\widetilde t_{i_1\cdots\widehat{i_\beta}\cdots i_{k+2}})=(-1)^{\alpha-\beta+1}\cdot{sgn(P_{\alpha,\beta}(\sigma))}\cdot \Delta^{P_{\alpha,\beta}(\sigma)}_{(u_1\cdots u_j;\,\emptyset)}(\widetilde t_{i_1\cdots\widehat{i_\alpha}\cdots i_{k+2}}).
\end{equation}
For each $\sigma\in S_{k+1}$ such that $\sigma(k+1)=\alpha$, we define an element $P_{\alpha,\beta}({\sigma})\in S_{k+1}$ as follows. Take $\widehat\sigma$ as in formula (\ref{hlift1}) and define $\eta\in S_{k+1}^{\alpha}$ by
\begin{equation}\label{hlift43}
\eta(l)=\left\{\begin{matrix}
\widehat\sigma(l)&\,\,\,\,\,\,\,\,\,\,\,{\rm if}\,\,\,\,\,\,\,\,\,\,\,1\leq l\leq k,\\\beta&\,\,\,\,\,\,\,\,\,\,\,{\rm if}\,\,\,\,\,\,\,\,\,\,\,l=k+1;
\end{matrix}\right.
\end{equation}
let $P_{\alpha,\beta}({\sigma})$ be  the inverse image of $\eta$ under the map $\overline{\bullet}$\,\,, that is, $\overline{P_{\alpha,\beta}({\sigma})}=\eta$. It is easy to verify that $P_{\alpha,\beta}({\sigma})\in \{\tau\in S_{k+1}\big|\tau(k+1)=\beta-1\}$. 

Then, by formula $(\ref{;empty})$, we have that
\begin{equation}\label{3dis42}
\begin{split}
\Delta^{P_{\alpha,\beta}({\sigma})}_{(u_1\cdots u_j;\,\emptyset)}&(\widetilde t_{i_1\cdots\widehat{i_\alpha}\cdots i_{k+2}})=tr\big\{\prod_{l=1}^k\big(g_{i_{\eta(l)} i_{\eta(u_1)}}^{-1}dg_{i_{\eta(l)} i_{\eta(u_1)}}\big)^{1-\delta_{u_1\cdots u_j}(l)}\cdot\big(g_{i_{\eta(l)}i_{\eta(u_1)}}^{-1}h_{i_{\eta(l)}}dh_{i_{\eta(l)}}^{-1}g_{i_{\eta(l)}i_{\eta(u_1)}}\big)^{\delta_{u_1\cdots u_j}(l)}\big\}\\
&=tr\big\{\prod_{l=1}^k\big(g_{i_{\widehat\sigma(l)} i_{\widehat\sigma(u_1)}}^{-1}dg_{i_{\widehat\sigma(l)} i_{\widehat\sigma(u_1)}}\big)^{1-\delta_{u_1\cdots u_j}(l)}\cdot\big(g_{i_{\widehat\sigma(l)}i_{\widehat\sigma(u_1)}}^{-1}h_{i_{\widehat\sigma(l)}}dh_{i_{\widehat\sigma(l)}}^{-1}g_{i_{\widehat\sigma(l)}i_{\widehat\sigma(u_1)}}\big)^{\delta_{u_1\cdots u_j}(l)}\big\}\\
&=\Delta^{\sigma}_{(u_1\cdots u_j;\,\emptyset)}(\widetilde t_{i_1\cdots\widehat{i_\beta}\cdots i_{k+2}}).
\end{split}
\end{equation}

Since 
\begin{equation}
sgn(P_{\alpha,\beta}({\sigma}))=(-1)^{k+2-\alpha}\cdot sgn([\overline{P_{\alpha,\beta}({\sigma})}]_{\alpha})=(-1)^{k+2-\alpha}\cdot sgn([\eta]_{\alpha})
\end{equation}
and
\begin{equation}
(-1)^{\alpha-\beta+1}\cdot sgn(\sigma)=(-1)^{\alpha-\beta+1}\cdot (-1)^{k+2-\beta}\cdot sgn([\widehat\sigma]_{\beta}),
\end{equation}
it suffices to prove that
\begin{equation}\label{sgn3}
sgn([\eta]_{\alpha})=-sgn([\widehat\sigma]_{\beta}).
\end{equation}
Noticing that $[\eta]_{\alpha}=\iota\circ[\widehat\sigma]_{\beta}\in S_{k+2}$ where $\iota\in S_{k+2}$ is the permutation $(\alpha,\beta)$, formula (\ref{sgn3}) holds.

Therefore, $\kappa^{\alpha}_3=(-1)^{\alpha-\beta+1}\cdot\kappa^{\beta}_3$.
\medskip

\noindent{\bf Case IV ($\kappa^{\alpha}_4=(-1)^{\alpha-\beta+1}\cdot\kappa^{\beta}_4$) :} Finally we will prove that for fixed integer $j$ such that $1\leq j\leq k-1$, and distinct integers $u_1,\cdots,u_j,x_1\in\{1,\cdots,k\}$, the following equality holds:
\begin{equation}
\begin{split}
\sum_{\substack{\sigma\in S_{k+1},\\\sigma(x_1)=\alpha}}{sgn(\sigma)}\cdot\Delta^{\sigma}_{(u_1\cdots u_j;\,x_1)}(\widetilde t_{i_1\cdots\widehat{i_\beta}\cdots i_{k+2}})=(-1)^{\alpha-\beta+1}\cdot\sum_{\substack{\tau\in S_{k+1},\\\tau(x_1)=\beta-1}}{sgn(\tau)}\cdot\Delta^{\tau}_{(u_1\cdots u_j;\,x_1)}(\widetilde t_{i_1\cdots\widehat{i_\alpha}\cdots i_{k+2}}).
\end{split}
\end{equation}

It suffices to prove that there exists a bijection
\begin{equation}
\begin{split}
P_{\alpha,\beta,u_1\cdots u_j,x_1}:\{\sigma\in S_{k+1}\big|\sigma(x_1)=\alpha\}\rightarrow\{\tau\in S_{k+1}\big|\tau(x_1)=\beta-1\},
\end{split}
\end{equation}
such that 
\begin{equation}\label{k42}
{sgn(\sigma)}\cdot\Delta^{\sigma}_{(u_1\cdots u_j;\,x_1)}(\widetilde t_{i_1\cdots\widehat{i_\beta}\cdots i_{k+2}})=(-1)^{\alpha-\beta+1}\cdot{sgn(P_{\alpha,\beta,u_1\cdots u_j,x_1}({\sigma}))}\cdot \Delta^{P_{\alpha,\beta,u_1\cdots u_j,x_1}({\sigma})}_{(u_1\cdots u_j;\,x_1)}(\widetilde t_{i_1\cdots\widehat{i_\alpha}\cdots i_{k+2}}).
\end{equation}

For each $\sigma\in S_{k+1}$ such that $\sigma(x_1)=\alpha$, we define an element $P_{\alpha,\beta,u_1\cdots u_j,x_1}({\sigma})\in S_{k+1}$ as follows. Take $\widehat\sigma$ as in formula (\ref{hlift1}) and define $\eta\in S_{k+1}^{\alpha}$ by
\begin{equation}\label{hlift44}
\eta(l)=\left\{\begin{matrix}
\widehat\sigma(l)&\,\,\,\,\,\,\,\,\,\,\,{\rm if}\,\,\,\,\,\,\,\,\,\,\,1\leq l\leq x_1-1,\\\beta&\,\,\,\,\,\,\,\,\,\,\,{\rm if}\,\,\,\,\,\,\,\,\,\,\,l=x_1,\\
\widehat\sigma(l)&\,\,\,\,\,\,\,\,\,\,\,{\rm if}\,\,\,\,\,\,\,\,\,\,\,x_1+1\leq l\leq k+1;\\
\end{matrix}\right.
\end{equation}
let $P_{\alpha,\beta,u_1\cdots u_j,x_1}$ be  the inverse image of $\eta$ under the map $\overline{\bullet}$\,\,, that is, $\overline{P_{\alpha,\beta,u_1\cdots u_j,x_1}}=\eta$. It is easy to verify that $P_{\alpha,\beta,u_1\cdots u_j,x_1}\in \{\tau\in S_{k+1}\big|\tau(x_1)=\beta-1\}$. 

Then by formula $(\ref{;empty})$, we have
\begin{equation}
\begin{split}
&\Delta^{\sigma}_{(u_1\cdots u_j;\,x_1)}(\widetilde t_{i_1\cdots\widehat{i_\beta}\cdots i_{k+2}})=tr\big\{\prod_{l=1}^k\big(g_{i_{\widehat\sigma(l)} i_{\widehat\sigma(u_1)}}^{-1}dg_{i_{\widehat\sigma(l)} i_{\widehat\sigma(u_1)}}\big)^{1-\delta_{x_1}(l)-\delta_{u_1\cdots u_j}(l)}\\
&\,\,\,\,\,\,\,\cdot\big(- g_{i_{\widehat\sigma(u_1)} i_{\widehat\sigma(k+1)}}dg_{i_{\widehat\sigma(u_1)} i_{\widehat\sigma(k+1)}}^{-1}\big)^{\delta_{x_1}(l)}\cdot\big(g_{i_{\widehat\sigma(l)}i_{\widehat\sigma(u_1)}}^{-1}h_{i_{\widehat\sigma(l)}}dh_{i_{\widehat\sigma(l)}}^{-1}g_{i_{\widehat\sigma(l)}i_{\widehat\sigma(u_1)}}\big)^{-\delta_{u_1\cdots u_j}(l)}\big\},
\end{split}
\end{equation}
and 
\begin{equation}
\begin{split}
&\Delta^{P_{\alpha,\beta,u_1\cdots u_j,x_1}(\sigma)}_{(u_1\cdots u_j;\,x_1)}(\widetilde t_{i_1\cdots\widehat{i_\alpha}\cdots i_{k+2}})=tr\big\{\prod_{l=1}^k\big(g_{i_{\eta(l)} i_{\eta(u_1)}}^{-1}dg_{i_{\eta(l)} i_{\eta(u_1)}}\big)^{1-\delta_{x_1}(l)-\delta_{u_1\cdots u_j}(l)}\\
&\,\,\,\,\,\,\,\cdot\big(-g_{i_{\eta(u_1)} i_{\eta(k+1)}}dg_{i_{\eta(u_1)} i_{\eta(k+1)}}^{-1}\big)^{\delta_{x_1}(l)}\cdot\big(g_{i_{\eta(l)}i_{\eta(u_1)}}^{-1}h_{i_{\eta(l)}}dh_{i_{\eta(l)}}^{-1}g_{i_{\eta(l)}i_{\eta(u_1)}}\big)^{\delta_{u_1\cdots u_j}(l)}\big\}.
\end{split}
\end{equation}
Hence,
\begin{equation}\label{dis44}
\begin{split}
\Delta^{\sigma}_{(u_1\cdots u_j;\,x_1)}(\widetilde t_{i_1\cdots\widehat{i_\beta}\cdots i_{k+2}})=\Delta^{P_{\alpha,\beta,u_1\cdots u_j,x_1}(\sigma)}_{(u_1\cdots u_j;\,x_1)}(\widetilde t_{i_1\cdots\widehat{i_\alpha}\cdots i_{k+2}}).
\end{split}
\end{equation}

Next we are going to prove that 
\begin{equation}
sgn(P_{\alpha,\beta,u_1\cdots u_j,x_1}({\sigma}))=(-1)^{\alpha-\beta+1}\cdot sgn(\sigma).
\end{equation}
Similarly, it suffices to prove that
\begin{equation}\label{sgn4}
sgn([\eta]_{\alpha})=-sgn([\sigma]_{\beta}).
\end{equation}
Since $[\eta]_{\alpha}=\iota\circ[\widehat\sigma]_{\beta}\in S_{k+2}$ where $\iota\in S_{k+2}$ is the permutation $(\alpha,\beta)$, formula (\ref{sgn4}) holds.

Therefore, $\kappa^{\alpha}_4=(-1)^{\alpha-\beta+1}\cdot\kappa^{\beta}_4$.
\medskip

We complete the proof of Claim IV.
\,\,\,\,$\endpf$
\medskip

By Claim III, for integer $k\geq 1$ and elements $j_1,\cdots, j_k\in I$,  we can define 
\begin{equation}\label{s4}
s_{j_1\cdots j_k}:=s_{\widehat{\beta}j_1\cdots j_{k};\,\beta j_1\cdots j_{k}},
\end{equation}
where $\beta$ is any element in $I$.  Noticing that only the transition functions $g_{\gamma\delta}$  where $\gamma,\delta\in\{j_1,\cdots,j_k\}$ appearing in (\ref{s0}),  $s_{j_1\cdots j_{k}}\in\Gamma(U_{j_1\cdots j_{k}},\Omega^k)$. Thus, we can define a \v Cech $(k-1)$-cochain $h_{k-1}(E,g,\widetilde g)$  by 
\begin{equation}
h_{k-1}(E,g,\widetilde g):=\bigoplus\limits_{j_1<\cdots<j_{k}}s_{j_1\cdots j_{k}}\in\bigoplus\limits_{j_1<\cdots<j_{k}}\Gamma(U_{j_1\cdots j_{k}},\Omega^k).
\end{equation} 
\begin{remark}
		Similar to Remark \ref{all}, we can extend the components of $h_{k-1}(E,g,\widetilde g)$  to all $k$-tuples of elements in $I$.  Then in the same way as  Lemma \ref{cochain}, we can show that formulas (\ref{s0}) and (\ref{s4}) is compatible with this extension. 
\end{remark}

Next, we prove that $\widehat f_k(E,\widetilde g)$ is cohomologous to $\widehat f_k(E,g)$
by $h_{k-1}(E,g,\widetilde g)$.
\medskip

\noindent{\bf Claim V : }
$\partial h_{k-1}(E,g,\widetilde g)=\widehat f_k(E,\widetilde g)-\widehat f_k(E,g)$. That is, for any elements $i_1,\cdots,i_{k+1}\in I$,
\begin{equation}\label{cobound1}
\widetilde t_{i_1\cdots i_{k+1}}- t_{i_1\cdots i_{k+1}}=\sum_{j=1}^{k+1}(-1)^{j-1}s_{i_1\cdots\widehat i_j\cdots i_{k+1}}\big|_{U_{i_1\cdots i_{k+1}}}.
\end{equation}

\medskip

\noindent{\bf Proof of Claim V :} 
Recall that 
\begin{equation}\label{t}
\sum_{\sigma\in S_{k+1}}\frac{sgn(\sigma)}{(k+1)!}\cdot \Delta^{\sigma}_{\emptyset,\emptyset}(\widetilde t_{i_1\cdots i_{k+1}})=\sum_{\sigma\in S_{k+1}}\frac{sgn(\sigma)}{(k+1)!}\cdot tr\big\{\prod_{l=1}^k\big(g_{i_{\sigma(l)} i_{\sigma(k+1)}}^{-1}dg_{i_{\sigma(l)} i_{\sigma(k+1)}}\big)\big\}=t_{t_{i_1\cdots i_{k+1}}}.
\end{equation}
Hence by Claims I and II, and formulas (\ref{tildet2}) and (\ref{t}), we have that
\begin{equation}\label{tildet3}
\begin{split}
&\widetilde t_{i_1\cdots i_{k+1}}-t_{i_1\cdots i_{k+1}}=\sum_{\sigma\in S_{k+1}}\frac{sgn(\sigma)}{(k+1)!} \sum_{\substack{1\leq j\leq k-1\,; \\  1\leq u_1<\cdots<u_j\leq k\,;\\ 1\leq v_1\leq k\,;\\u_1,\cdots,u_j,v_1\, {\text{are distinct}}}}\Delta^{\sigma}_{u_1\cdots u_j,v_{1}}(\widetilde t_{i_1\cdots i_{k+1}})\\
&\,\,\,\,+\sum_{\sigma\in S_{k+1}}\frac{sgn(\sigma)}{(k+1)!} \sum_{\substack{ j=0\,; \\ 1\leq v_1\leq k\,;}}\Delta^{\sigma}_{\emptyset,v_{1}}(\widetilde t_{i_1\cdots i_{k+1}})+\sum_{\sigma\in S_{k+1}}\frac{sgn(\sigma)}{(k+1)!} \sum_{\substack{1\leq j\leq k\,; \\  1\leq u_1<u_2<\cdots<u_j\leq k\,;}}\Delta^{\sigma}_{u_1\cdots u_j,\emptyset}(\widetilde t_{i_1\cdots i_{k+1}})\\
&=\sum_{\sigma\in S_{k+1}}\frac{sgn(\sigma)}{(k+1)!} \sum_{\substack{1\leq j\leq k-1\,; \\  1\leq u_1<\cdots<u_j\leq k\,;\\ 1\leq v_1\leq k\,;\\u_1,\cdots,u_j,v_1\, {\text{are distinct}}}}\Delta^{\sigma}_{u_1\cdots u_j,v_{1}}(\widetilde t_{i_1\cdots i_{k+1}})+\sum_{\sigma\in S_{k+1}}\frac{sgn(\sigma)}{(k+1)!} \sum_{\substack{j=0\,; \\ 1\leq v_1\leq k\,;}}\Delta^{\sigma}_{\emptyset,v_{1}}(\widetilde t_{i_1\cdots i_{k+1}})\\
&\,\,\,\,\,\,\,\,\,\,\,\,+\sum_{\sigma\in S_{k+1}}\frac{sgn(\sigma)}{(k+1)!}\sum_{\substack{1\leq j\leq k\,; \\  1\leq u_1<u_2<\cdots<u_j\leq k}}\Delta^{\sigma}_{(u_1\cdots u_j;\,\emptyset)}(\widetilde t_{i_1\cdots i_{k+1}})\\
&\,\,\,\,\,\,\,\,\,\,\,\,+\sum_{\sigma\in S_{k+1}}\frac{sgn(\sigma)}{(k+1)!}\sum_{\substack{1\leq j\leq k-1\,; \\  1\leq u_1<u_2<\cdots<u_j\leq k;\\   1\leq x_1\leq k,\, x_1\notin\{u_1,\cdots,u_j\}}}\Delta^{\sigma}_{(u_1\cdots u_j;\,x_1)}(\widetilde t_{i_1\cdots i_{k+1}}).\\
\end{split}
\end{equation}

Notice that for each fixed $\xi$  such that $1\leq \xi\leq k+1$,  $S_{k+1}$ is the disjoint union of the following sets, 
\begin{equation}
\begin{split}
\big\{\sigma\in S_{k+1}\big|\sigma(\xi)=\alpha\big\}\,\,\,{\rm where}\,\,\alpha=1,\cdots,k+1.
\end{split}
\end{equation}
Then, by interchanging the order of summation in formula (\ref{t}), we have
\begin{equation}
\begin{split}
\widetilde t_{i_1\cdots i_{k+1}}-&t_{i_1\cdots i_{k+1}}=\sum_{\alpha=1}^{k+1}\,\,\,\,\frac{1}{(k+1)!} \sum_{\substack{1\leq j\leq k-1\,; \\  1\leq u_1<\cdots<u_j\leq k\,;\\ 1\leq v_1\leq k\,;\\u_1,\cdots,u_j,v_1\, {\text{are distinct}}}}\sum_{\substack{\sigma\in S_{k+1},\\\sigma(v_1)=\alpha}}sgn(\sigma)\cdot\Delta^{\sigma}_{u_1\cdots u_j,v_{1}}(\widetilde t_{i_1\cdots i_{k+1}})\\
&+\sum_{\alpha=1}^{k+1}\,\,\,\,\frac{1}{(k+1)!}\sum_{\substack{j=0\,; \\ 1\leq v_1\leq k\,;}}
\sum_{\substack{\sigma\in S_{k+1},\\\sigma(v_1)=\alpha}}sgn(\sigma)\cdot\Delta^{\sigma}_{\emptyset,v_{1}}(\widetilde t_{i_1\cdots i_{k+1}})\\
&+\sum_{\alpha=1}^{k+1}\,\,\,\,\frac{1}{(k+1)!}\sum_{\substack{1\leq j\leq k\,; \\  1\leq u_1<u_2<\cdots<u_j\leq k;}}\sum_{\substack{\sigma\in S_{k+1},\\\sigma(k+1)=\alpha}}{sgn(\sigma)}\cdot\Delta^{\sigma}_{(u_1\cdots u_j;\,\emptyset)}(\widetilde t_{i_1\cdots i_{k+1}})\\
&+\sum_{\alpha=1}^{k+1}\,\,\,\,\frac{1}{(k+1)!}\sum_{\substack{1\leq j\leq k-1\,; \\  1\leq u_1<u_2<\cdots<u_j\leq k;\\   1\leq x_1\leq k,\, x_1\notin\{u_1,\cdots,u_j\}}}\sum_{\substack{\sigma\in S_{k+1},\\\sigma(x_1)=\alpha}}{sgn(\sigma)}\cdot\Delta^{\sigma}_{(u_1\cdots u_j;\,x_1)}(\widetilde t_{i_1\cdots i_{k+1}}).\\
\end{split}
\end{equation}
Recalling formula (\ref{s0}), we conclude that
\begin{equation}
\widetilde t_{i_1\cdots i_{k+1}}- t_{i_1\cdots i_{k+1}}=\sum_{\alpha=1}^{k+1}(-1)^{\alpha-1}s_{i_1\cdots\widehat i_\alpha\cdots i_{k+1}}\big|_{U_{i_1\cdots i_{k+1}}}.
\end{equation}

We complete the proof of Claim V.\,\,\,\,$\endpf$
\medskip

Therefore, we finish the proof of Lemma \ref{invariant}.
\,\,\,\,$\endpf$
\medskip

\noindent{\bf Proof of Theorem \ref{tcf} :}  Theorem \ref{tcf} follows from Lemmas \ref{cocycle} and \ref{invariant}.
\,\,\,\,$\endpf$
\begin{remark}
	If instead of $\Omega^k$ we consider the sheaf $\mathcal A^k$ consisting of germs of smooth differential $k$-forms, the above construction still holds. However, since $H^k(X,\mathcal A^k)=0$ for $k\geq 1$, the resulting invariant cocycle is trivial.
	That seems the reason why one has to go to the \v Cech-de Rham complex as \cite{Sh}, in order to define the Chern classes of a smooth vector bundle.
\end{remark}
\begin{question}
	Can we give a combinatorial formula for the  Chern classes of a smooth vector bundle as an elemnt in the \v Cech-de Rham complex explicitly?
\end{question}

\section {Invariants of holomorphic vector bundles I : T invariants in $H^k(X,\Omega^k)$ and $H^1(X,d\Omega^0)$}
Let $X$ be a complex manifold.  In this section, we will define the $T$ invariants and the refined first $T$ invariant  of a holomorphic vector bundle over $X$. 
\subsection{T invariants in the Dolbeault cohomolgy}

\noindent{\bf Proof of Theorem \ref{cfd} :} Let $\mathcal U:=\{U_i\}_{i\in I}$ be an open cover of $X$, where $I$ is an ordered set,  and $g:=\{g_{i_1i_2}\}$ be a system of transition functions associated with a certain trivialization of $E$ with respect to $\mathcal U$. By Theorem \ref{tcf}, $f_{k,\,\mathcal U}(E)$ is an well-defined element in $\check H^k(\mathcal U,\Omega^k)$.

Suppose $\mathcal V:=\{V_j\}_{j\in J}$ is a refinement of $\mathcal U$. Choose a map $c:J\longrightarrow I$ such that $V_j\subset U_{c(j)}$ for all $j\in J$. For each $j_1,j_2\in J$, define a matrix-valued function $\widetilde g_{j_1j_2}$ as follows,
\begin{equation}\label{wg}
\begin{split}
\widetilde g_{j_1j_2}:V_{j_1j_2}&\longrightarrow GL(M,\mathbb C)\\
x&\longmapsto g_{c(j_1)c(j_2)}(x);
\end{split}
\end{equation}
then $\widetilde g:=\{\widetilde g_{j_1j_2}\}$ is a system of transition functions associated with a certain trivialization of $E$ with respect to $\mathcal V$.  On the other hand, recall that the natural restriction $I_{\mathcal U\mathcal V}:\check H^k(\mathcal U,\Omega^k)\rightarrow \check H^k(\mathcal V,\Omega^k)$ is induced by the following map between \v Cech complexes
\begin{equation}\label{ga}
\begin{split}
\Gamma_{\mathcal U\mathcal V}:\widehat{\mathcal C}^{\bullet}(\mathcal U,\Omega^k)&\longrightarrow \widehat{\mathcal C}^{\bullet}(\mathcal V,\Omega^k)\\
(\xi_{i_1\cdots i_p})&\longmapsto(\xi_{c(j_1)\cdots c(j_p)}\big|_{V_{j_1\cdots j_p}}).
\end{split}
\end{equation}

In the following, we will show that $I_{\mathcal U\mathcal V}(f_{k,\,\mathcal U}(E))=f_{k,\,\mathcal V}(E)$. Let $\widehat f_{k,\,\mathcal U}(E)=(t_{i_1\cdots i_{k+1}})$ be the \v Cech cocycle constructed in 
Theorem \ref{tcf} corresponding to $f_{k,\,\mathcal U}(E)$, where
\begin{equation}
\begin{split}
t_{i_1\cdots i_{k+1}}=\sum_{\sigma\in S_{k+1}}\frac{sgn(\sigma)}{(k+1)!}&\cdot tr\big(g^{-1}_{i_{\sigma(1)}i_{\sigma(k+1)}}dg_{i_{\sigma(1)}i_{\sigma(k+1)}}g^{-1}_{i_{\sigma(2)}i_{\sigma(k+1)}}dg_{i_{\sigma(2)}i_{\sigma(k+1)}}\cdot\\
&\cdot g^{-1}_{i_{\sigma(3)}i_{\sigma(k+1)}}dg_{i_{\sigma(3)}i_{\sigma(k+1)}}\cdots g^{-1}_{i_{\sigma(k)}i_{\sigma(k+1)}}dg_{i_{\sigma(k)}i_{\sigma(k+1)}}\big).
\end{split}
\end{equation}
Then by formula (\ref{ga}), $\Gamma_{\mathcal U\mathcal V}(E)$ maps $\widehat f_{k,\,\mathcal U}(E)$ to a \v Cech cocycle $\Gamma_{\mathcal U\mathcal V}(\widehat f_{k,\,\mathcal U}(E)):=(\widetilde t_{j_1\cdots j_{k+1}})$, where
\begin{equation}\label{wt}
\begin{split}
\widetilde t_{j_1\cdots j_{k+1}}=\sum_{\sigma\in S_{k+1}}&\frac{sgn(\sigma)}{(k+1)!}\cdot tr\big(g^{-1}_{c(j_{\sigma(1)})c(j_{\sigma(k+1)})}dg_{c(j_{\sigma(1)})c(j_{\sigma(k+1)})}g^{-1}_{c(j_{\sigma(2)})c(j_{\sigma(k+1)})}dg_{c(j_{\sigma(2)})c(j_{\sigma(k+1)})}\cdot\\
&\cdot g^{-1}_{c(j_{\sigma(3)})c(j_{\sigma(k+1)})}dg_{c(j_{\sigma(3)})c(j_{\sigma(k+1)})}\cdots g^{-1}_{c(j_{\sigma(k)})c(j_{\sigma(k+1)})}dg_{c(j_{\sigma(k)})c(j_{\sigma(k+1)})}\big).
\end{split}
\end{equation} 
Noticing formula (\ref{wg}), we can rewrite formula (\ref{wt}) as follows,
\begin{equation}\label{wwt}
\begin{split}
\widetilde t_{j_1\cdots j_{k+1}}=\sum_{\sigma\in S_{k+1}}\frac{sgn(\sigma)}{(k+1)!}&\cdot tr\big(\widetilde g^{-1}_{j_{\sigma(1)}j_{\sigma(k+1)}}d\widetilde g_{j_{\sigma(1)}j_{\sigma(k+1)}}\widetilde g^{-1}_{j_{\sigma(2)}j_{\sigma(k+1)}}d\widetilde g_{j_{\sigma(2)}j_{\sigma(k+1)}}\cdot\\
&\cdot \widetilde g^{-1}_{j_{\sigma(3)}j_{\sigma(k+1)}}d\widetilde g_{j_{\sigma(3)}j_{\sigma(k+1)}}\cdots \widetilde g^{-1}_{j_{\sigma(k)}j_{\sigma(k+1)}}d\widetilde g_{j_{\sigma(k)}j_{\sigma(k+1)}}\big).
\end{split}
\end{equation}
Then $I_{\mathcal U\mathcal V}(f_{k,\,\mathcal U}(E))= f_{k,\,\mathcal V}(E)$.
 
Since any two open covers $\mathcal U^1$ and $\mathcal U^2$ have a common refinement $\mathcal V$, 
\begin{equation}
[f_{k,\,\mathcal U^1}(E)]=[f_{k,\,\mathcal U^2}(E)]=[f_{k,\,\mathcal V}(E)],
\end{equation}
where $[f_{k,\,\bullet}(E)]$ is the equivalence class of $f_{k,\,\bullet}(E)$ in $\check H^k(X,\Omega^k)$.
Then $f_k(E)$ is independent of the choice of the  cover $\mathcal U$. We complete the proof of Theorem \ref{cfd}. \,\,\,\,\,\,$\endpf$
\medskip

\noindent{\bf Proof of Corollary \ref{comt} :} 
Since $X$ is compact,  each open cover of $X$ has a finite Stein cover as its refinement. In particular, there is a finite Stein cover $\mathcal U$ of $X$.

Recall the following natural homomorpshims,
\begin{equation}
\check H^k({\mathcal U},\Omega^k)\xrightarrow{L_{\mathcal U}}\check H^k(X,\Omega^k)\xrightarrow{L} H^k(X,\Omega^k).
\end{equation}
Notice that $\check H^k({\mathcal U},\Omega^k)\cong H^k(X,\Omega^k)$, for coherent sheaves are acyclic over Stein sets. Then, $L$ is surjective.  

Suppose $\mathcal V$ is an open cover of $X$ and $a\in\check H^k(\mathcal V,\Omega^k)$ such that $L([a])=0$, where $[a]$ is the equivalence class of $a$ in $\check H^k(X,\Omega^k)$. Then, there is a finite Stein open cover $\mathcal W$ of $X$ as a refinement of $\mathcal V$. Since $(L\circ L_{\mathcal W})([I_{\mathcal V\mathcal W}(a)])=L([a])=0$ and $L\circ L_{\mathcal W}$ is an isomorphism, $[a]=[I_{\mathcal V\mathcal W}(a)]=0$. Therefore, $\check H^k(X,\Omega^k)\cong H^k(X,\Omega^k)$.  Then $f_k(E)=\widehat f_{k,\,\mathcal U}(E)$ under the natural isomorphisms $\check H^k({\mathcal U},\Omega^k)\cong\check H^k(X,\Omega^k)\cong H^k(X,\Omega^k)$.
\,\,\,\,\,\,$\endpf$
\medskip

Before proving Theorem \ref{compare}, we recall the notion of the Chern character  of  a vector bundle (see \cite{Hi} for reference). Let $E$ be a vector bundle over $X$ and $\gamma_1,\cdots,\gamma_M$ be the Chern roots of $E$. Then, ch$(E)$ the chern character of $E$ is defined as follows,
\begin{equation}
\begin{split}
{\rm ch}(E):=\sum_{i=1}^M\exp{\gamma_i}=\sum_{k=0}^{\infty}\frac{1}{k!}\cdot(\gamma^k_i+\cdots+\gamma^k_M).
\end{split}
\end{equation}
Notice that ch$(E)$ can be written as a polynomial whose variables are the Chern classes of $E$ by Newton's identities. Moreover, when $E$ is a holomorphic vector bundle, one can view ch$(E)$ as an element in the ring $\oplus_{i=0}^{\infty}H^k(X,\Omega^k)$.

We first prove Theorem \ref{compare} in the case when $E$ is a line bundle.
\begin{lemma}\label{line}
	Let $X$ be a compact complex manifold.  Let $\pi:E\rightarrow X$ be a holomorphic  line bundle of $X$.  $f_k(E)$ is defined as  in Theorem \ref{cfd}.
	Then, 
	\begin{equation}\label{dudu}
	{\rm ch}(E)=\sum_{k=0}^{\infty}\frac{1}{k!\cdot(2\pi\sqrt{-1})^{k}}\cdot f_{k}(E).
	\end{equation}
\end{lemma}
\noindent{\bf Proof of Lemma \ref{line} :}  Let  $\mathcal U:=\{U_i\}_{i=1}^l$ be a Stein  open cover of $X$ and $g:=\{g_{i_1i_2}\}$ be a system of transition functions associated with a certain trivialization of $E$ with respect to $\mathcal U$.

When $k=1$, formula (\ref{t11}) reads
\begin{equation}\label{t1}
t_{i_1i_2}=\frac{1}{2}tr(g^{-1}_{i_1i_2}dg_{i_1i_2})-\frac{1}{2}tr(g^{-1}_{i_2i_1}dg_{i_2i_1}).
\end{equation}
Since $E$ is a line bundle and hence $g_{i_1i_2}$ is a function, we have that
\begin{equation}\label{d-1}
dg_{i_2i_1}=d(g^{-1}_{i_1i_2})=-g^{-1}_{i_1i_2}dg_{i_1i_2}g^{-1}_{i_1i_2}.
\end{equation}
Substituting formula (\ref{d-1}) into (\ref{t1}), we derive that
\begin{equation}
t_{i_1i_2}=\frac{1}{2}g^{-1}_{i_1i_2}dg_{i_1i_2}-\frac{1}{2}g^{-1}_{i_2i_1}(-g^{-1}_{i_1i_2}dg_{i_1i_2}g^{-1}_{i_1i_2})=d\log(g_{i_1i_2}).
\end{equation}
Similar to the proof of Proposition $1$ in \S 1.1 of \cite{GH}, we conclude that $f_1(E)=(2\pi\sqrt{-1})\cdot c_1(E)$ in $H^1(X,\Omega^1)$.

When $k\geq 2$, we can simplify formula (\ref{t11}) as follows,
\begin{equation}\label{k>1}
\begin{split}
t_{i_1\cdots i_{k+1}}&=\frac{1}{k+1}\sum_{\alpha=1}^{k+1}\sum_{\substack{1\leq u^{\alpha}_1<\cdots<u^{\alpha}_k\leq k+1\,; \\ u^{\alpha}_1,\cdots,u^{\alpha}_k,\alpha\, {\text{are distinct}}}} sgn{\Large(}\Big(\begin{matrix}1&2&\cdots&k&k+1\\u^{\alpha}_1&u^{\alpha}_2&\cdots&u^{\alpha}_k&\alpha\end{matrix}\Big){\Large)}\\
&\cdot\big(d\log g_{i_{u^{\alpha}_1}i_\alpha}\wedge d\log g_{i_{u^{\alpha}_2}i_\alpha}\wedge d\log g_{i_{u^{\alpha}_3}i_\alpha}\wedge\cdots \wedge d\log g_{i_{u^{\alpha}_k}i_\alpha}\big).
\end{split}
\end{equation}

Since $\log g_{ik}=\log g_{ij}+\log g_{jk}$, we have that
\begin{equation}\label{k>11}
\begin{split}
t_{i_1\cdots i_{k+1}}&=\frac{1}{k+1}\sum_{\alpha=1}^{k+1}\sum_{\substack{1\leq u^{\alpha}_1<\cdots<u^{\alpha}_k\leq k+1\,; \\ u^{\alpha}_1,\cdots,u^{\alpha}_k,\alpha\, {\text{are distinct}}}} sgn{\Large(}\Big(\begin{matrix}1&2&\cdots&k&k+1\\u^{\alpha}_1&u^{\alpha}_2&\cdots&u^{\alpha}_k&\alpha\end{matrix}\Big){\Large)}\\
&\cdot\big((d\log g_{i_{u^{\alpha}_1}i_{(k+1)}}+d\log g_{i_{(k+1)}i_\alpha})\wedge\cdots \wedge (d\log g_{i_{u^{\alpha}_k}i_{(k+1)}}+d\log g_{i_{(k+1)}i_\alpha})\big)\\
&=d\log g_{i_{1}i_{(k+1)}}\wedge d\log g_{i_{2}i_{(k+1)}}\wedge\cdots\wedge d\log g_{i_{k}i_{(k+1)}}\\
&=d\log g_{i_{1}i_{2}}\wedge d\log g_{i_{2}i_{3}}\wedge\cdots\wedge d\log g_{i_{k}i_{(k+1)}}.
\end{split}
\end{equation}
By formula ($14.24$) on Page 174 of \cite{BT},  $f_k(E)$ is the $k$-fold cup product of $f_1(E)$ with iteself, that is, 
\begin{equation}
f_k(E)=f_1(E)\cup f_1(E)\cup\cdots\cup f_1(E)\in \check H^k(\mathcal U,\Omega^k).
\end{equation}

Then,
\begin{equation}\label{cup}
\begin{split}
{\rm ch}(E)&=\sum_{k=0}^{\infty}\frac{1}{k!}\cdot c_{1}^k(E)=\sum_{k=0}^{\infty}\frac{1}{k!}\cdot\Big(\frac{ f_{1}(E)}{2\pi\sqrt{-1}}\Big)^k\\
&=\sum_{k=0}^{\infty}\frac{1}{k!\cdot(2\pi\sqrt{-1})^{k}}\cdot f_{k}(E).
\end{split}
\end{equation}
We complete the proof of Lemma \ref{line}. \,\,\,\,\,$\endpf$
\medskip

Next we prove Theorem \ref{compare} for the vector bundle with a full flag structure.

\begin{lemma}\label{split}
	Let $X$ be a compact complex manifold.  Let $\pi:E\rightarrow X$ be a holomorphic   vector bundle of $X$ with a filtration by holomoprhic subbundles
	\begin{equation}
	E=E_M\supset E_{M-1}\supset\cdots\supset E_2\supset E_1\supset E_0=0
	\end{equation}
	with line bundle quotients $L_i=E_i/E_{i-1}$. $f_k(E)$ is defined as  in Theorem \ref{cfd}.
	Then, 
	\begin{equation}\label{cf1}
	{\rm ch}(E)=\sum_{k=0}^{\infty}\frac{1}{k!\cdot(2\pi\sqrt{-1})^{k}}\cdot f_{k}(E).
	\end{equation}
\end{lemma}
\noindent{\bf Proof of Lemma \ref{split} :} Since $E$ has a full flag structure, there is a Stein open cover   $\mathcal U:=\{U_i\}_{i=1}^l$ of $X$ and a system of transition functions $g:=\{g_{ij}\}$  associated with a certain trivialization of $E$ with respect to $\mathcal U$  such that $g_{ij}$ is an upper triangular matrix as follows,
\begin{equation}
g_{ij}=\left(\begin{matrix}
g^1_{ij}&*&\cdots&*&*\\
0&g^2_{ij}&\cdots&*&*\\
\cdots&\cdots&\cdots&\cdots&\cdots\\
0&0&\cdots&g^{M-1}_{ij}&*\\
0&0&\cdots&0&g_{ij}^M
\end{matrix}\right)\,\,\,{\rm for}\,\,\,1\leq i,j\leq l.
\end{equation}

Similar to formulas (\ref{k>1}) and (\ref{k>11}), computation yields that
\begin{equation}
\begin{split}
t_{i_1\cdots i_{k+1}}&=\sum_{\sigma\in S_{k+1}}\frac{sgn(\sigma)}{(k+1)!}\cdot tr\big(g^{-1}_{i_{\sigma(1)\sigma(k+1)}}dg_{i_{\sigma(1)\sigma(k+1)}}g^{-1}_{i_{\sigma(2)\sigma(k+1)}}dg_{i_{\sigma(2)\sigma(k+1)}}\cdot\\
&\,\,\,\,\,\,\,\,\,\,\,\,\,\,\,\,\,\,\,\,\,\,\,\,\,\,\cdot g^{-1}_{i_{\sigma(3)}i_{\sigma(k+1)}}dg_{i_{\sigma(3)\sigma(k+1)}}\cdots g^{-1}_{i_{\sigma(k)}i_{\sigma(k+1)}}dg_{i_{\sigma(k)}i_{\sigma(k+1)}}\big)\\
&=\sum_{j=1}^M\,\,\,\frac{1}{k+1}\,\,\,\sum_{\alpha=1}^{k+1}\sum_{\substack{1\leq u^{\alpha}_1<\cdots<u^{\alpha}_k\leq k+1\,; \\ u^{\alpha}_1,\cdots,u^{\alpha}_k,\alpha\, {\text{are distinct}}}} sgn{\Large(}\Big(\begin{matrix}1&2&\cdots&k&k+1\\u^{\alpha}_1&u^{\alpha}_2&\cdots&u^{\alpha}_k&\alpha\end{matrix}\Big){\Large)}\\
&\,\,\,\,\,\,\,\,\,\,\,\,\,\,\,\,\,\,\,\,\,\,\,\,\cdot\big(d\log g^j_{i_{u^{\alpha}_1}i_\alpha}\wedge d\log g^j_{i_{u^{\alpha}_2}i_\alpha}\wedge d\log g^j_{i_{u^{\alpha}_3}i_\alpha}\wedge\cdots \wedge d\log g^j_{i_{u^{\alpha}_k}i_\alpha}\big)\\
&=\sum_{j=1}^M\,\,\,\,d\log g^j_{i_{1}i_{(k+1)}}\wedge d\log g^j_{i_{2}i_{(k+1)}}\wedge\cdots\wedge d\log g^j_{i_{k}i_{(k+1)}}\\
&=\sum_{j=1}^M\,\,\,\,d\log g^j_{i_{1}i_{2}}\wedge d\log g^j_{i_{2}i_{3}}\wedge\cdots\wedge d\log g^j_{i_{k}i_{(k+1)}}\\
&=:t^1_{i_1\cdots i_{k+1}}+t^2_{i_1\cdots i_{k+1}}+\cdots+t^M_{i_1\cdots i_{k+1}}.
\end{split}
\end{equation}
Notice that by Lemma \ref{line}, $f_k(L_j)$, $j=1,\cdots,M$, is given by
\begin{equation}
f_k(L_j)=\bigoplus\limits_{1\leq i_1<\cdots<i_{k+1}\leq l}t^j_{i_1\cdots i_{k+1}}\in\bigoplus\limits_{1\leq i_1<\cdots<i_{k+1}\leq l}\Gamma(U_{i_1\cdots i_{k+1}},\Omega^k);
\end{equation}
hence,
\begin{equation}
f_k(E)=\sum_{j=1}^Mf_k(L_j)=\sum_{j=1}^M\big(f_1(L_j)\big)^k.
\end{equation}

Therefore,
\begin{equation}
\begin{split}
{\rm ch}(E)&=\sum_{k=0}^{\infty}\frac{1}{k!}\sum_{j=1}^Mc_{1}^k(L_j)=\sum_{k=0}^{\infty}\frac{1}{k!}\sum_{j=1}^M\Big(\frac{ f_{1}(L_j)}{2\pi\sqrt{-1}}\Big)^k\\
&=\sum_{k=0}^{\infty}\frac{1}{k!\cdot(2\pi\sqrt{-1})^{k}}\sum_{j=1}^M f_{1}^k(L_j)=\sum_{k=0}^{\infty}\frac{1}{k!\cdot(2\pi\sqrt{-1})^{k}}\cdot f_{k}(E).
\end{split}
\end{equation}
We complete the proof of Lemma \ref{split}.
\,\,\,\,$\endpf$
\medskip

\noindent{\bf Proof of Theorem \ref{compare} :} The idea is to use the splitting principle.

As in \cite{Fu}, we construct a split manifold for $(X,E)$ as follows. Let $P(E)$ be the projective bundle of $E$; denote the projection map by $\sigma:P(E)\rightarrow X$ and denote the pullback of $E$ to $P(E)$ by $\sigma^{-1}(E)$. Repeating this procedure, we can find a compact complex manifold $X^{\prime}$ and a holomorphic map $\pi:X^{\prime}\rightarrow X$ such that
the pullback of $E$ to $X^{\prime}$ has a filtration of holomorphic subbundles:
\begin{equation}
\pi^{-1}(E)=E_M\supset E_{M-1}\supset\cdots\supset E_2\supset E_1\supset E_0=0
\end{equation}
with line bundle quotients $L_i=E_i/E_{i-1}$.

Let $\mathcal U:=\{U_i\}_{i=1}^l$ be a Stein cover of $X$; let $\mathcal V$ be a Stein cover of $X^{\prime}$ such that $\mathcal V$ is a refinement of the cover $\mathcal U^*:=\pi^*\mathcal U=\{\pi^{-1}(U_i)\}_{i=1}^l$ of $X^{\prime}$. Then, we have the following commutative diagram for each positive integer $k$:
\begin{equation}
\xymatrix{
	&\check H^k({\mathcal U},\Omega_X^k)\ar[d]^{ L_{\mathcal U}} \ar[r]^{T_k \,\,\,\,\,\,}&\check H^k({\mathcal V},\Omega_{X^{\prime}}^k) \ar[d]^{L_{\mathcal V}}\\
	&\check H^k(X,\Omega_X^k)\ar[d]^{L_X} \ar[r]^{\Pi^*_k \,\,\,\,\,\,}&\check H^k(X^{\prime},\Omega_{X^{\prime}}^k) \ar[d]^{L_{X^{\prime}}}\\
	& H^k(X,\Omega_X^k)\ar[r]^{ \pi^*_{k}\,\,\,\,\,\,}&H^k(X^{\prime},\Omega_{X^{\prime}}^k).}
\end{equation}
Here $T_{k}=I_{\mathcal U^*V}\circ\check\pi^*_{k}$, where $\check\pi^*_{k}$ is the natural pullback homorphism from $\check H^k({\mathcal U},\Omega_{X}^k)$ to $\check H^k({\mathcal U^*},\Omega_{X^{\prime}}^k)$; $\Pi^*_k$ is the induced homomorphism between the direct limits; $\pi^*_{k}$ is the natural pullback map from $H^k(X,\Omega_X^k)$ to $H^k(X^{\prime},\Omega_{X^{\prime}}^k)$. Moreover, $L_{\mathcal U},L_{\mathcal V},L_{X}$ and $L_{X^{\prime}}$ are all isomorphisms.

Noticing Corollary \ref{comt} and Remark \ref{rcomt}, let $\widehat f_{k,\,\mathcal U}(E)$ be the \v Cech cocycle representing  $f_k(E)\in H^k(X,\Omega^k)$ and let $\widehat f_{k,\,\mathcal V}(\pi^{-1}(E))$ be the \v Cech cocycle representing  $f_k(\pi^{-1}(E))\in H^k(X^{\prime},\Omega_{X^{\prime}}^k)$.  Recall that the pullback of a system of transition functions of $E$ with respect to $\mathcal U$ gives a system of transition functions of $\pi^{-1}(E)$ with respect to $\mathcal U^*$. In the same way as the proof of Theorem \ref{cfd} (see formulas (\ref{wt}) and (\ref{wwt})), we can prove the following identity,
\begin{equation}
(\Pi^*_{k}\circ L_{\mathcal U})\big(\widehat f_{k,\,\mathcal U}(E)\big)=(L_{\mathcal V})\big(\widehat f_{k,\,\mathcal V}(\pi^{-1}(E))\big),
\end{equation}
and hence
\begin{equation}
(\pi^*_{k})(f_k(E))=f_k(\pi^{-1}(E)).
\end{equation}

On the other hand, by representing the Chern classes by the integrations of invariant polynomials of the curvature forms of a Hermitian metric, it is clear that 
\begin{equation}
(\pi^*_{k})({\rm ch}(E))={\rm ch}(\pi^{-1}(E))\in \oplus_{i=0}^{\infty}H^k(X^{\prime},\Omega_{X^{\prime}}^k).
\end{equation}
By the Hirsch lemma (see Lemma 18 of \cite{CFGU} or Lemma 3.3 of \cite{RYY}),
the following natural homomorphism is an embedding
\begin{equation}
\pi^*:\oplus_{i=0}^{\infty}H^k(X,\Omega_X^k)\hookrightarrow\oplus_{i=0}^{\infty}H^k(X^{\prime},\Omega_{X^{\prime}}^k);
\end{equation}
in particular,  $\pi^*_{k}$ is injective. Therefore, in order to prove Theorem \ref{compare}, it suffices to prove that
\begin{equation}
{\rm ch}(\pi^{-1}(E))=\sum_{k=0}^{\infty}\frac{1}{k!\cdot(2\pi\sqrt{-1})^{k}}\cdot f_{k}(\pi^{-1}(E)).
\end{equation}
By Lemma \ref{split}, we draw the conclusion. \,\,\,\,\,\,\,$\endpf$
\medskip

\subsection{Refinement of the first Chern class for a holomorphic vector bundle}

In the following, we will show that the first $T$ invariant of a holomorphic vector bundle can be lift to the group $\check H^1(X,d\Omega^0)$.
Recall that by Theorem \ref{compare}, the first $T$ invariant coincides with the first Chern class, up to a certain normalized factor, in the Dolbeault cohomology group $H^1(X,\Omega^1)$. Hence, the lifted first $T$ invariant can be viewed as a refinement of the first Chern class in  $H^1(X,\Omega^1)$. 
\medskip

\noindent{\bf Proof of Theorem \ref{fla} :} The proof is similar to the proof of Theorem \ref{cfd}. We only need to show that the coefficients of the \v Cech cochains take values in the sheaf $d\Omega^0$. 

Firstly, we have that
\begin{equation}\label{dg}
d\Big(tr\big(g^{-1}_{i_{1}i_{2}}dg_{i_1i_2}\big)\Big)=tr\big(-g^{-1}_{i_{1}i_{2}}dg_{i_1i_2}g^{-1}_{i_{1}i_{2}}dg_{i_1i_2}\big)=0,
\end{equation}
for $tr(AA)=0$ for any differential $1$-form-valued $M\times M$ matrix $A$. Therefore, $t_{i_1i_{2}}\in\Gamma(U_{i_1i_{2}},d\Omega^0)$.

Next suppose that $\{\widetilde g_{ij}\}$ is another system of transition functions associated with a certain trivialization of $E$ with respect to $\mathcal U$. Then, there exists a \v Cech $0$-cochain
\begin{equation}
h:=\bigoplus\limits_{i_1}h_{i_1}\in\bigoplus\limits_{i_1}\Gamma_{hol} (U_{i_1},GL(M,\mathbb C)),
\end{equation}
such that $\widetilde g_{ij}=h_i^{-1}g_{ij}h_j$ for $i,j\in I.$  As in the proof of Lemma \ref{invariant}, we can define a \v Cech $0$-cochain
\begin{equation}\label{h0}
h_{0}(E,g,\widetilde g):=\bigoplus\limits_{i_1}s_{i_1}\in\bigoplus\limits_{i_1}\Gamma(U_{i_1},\Omega^1),
\end{equation}
such that, for any elements $i_1,i_2\in I$, 
\begin{equation}\label{hcobound}
\widetilde t_{i_1i_2}- t_{i_1i_{2}}=\big(s_{i_1}-s_{i_2}\big)\big|_{U_{i_1i_2}}.
\end{equation}
Recalling formula (\ref{s0}), computation yields that
\begin{equation}\label{si1}
\begin{split}
s_{i_1}=-tr\big(h^{-1}_{i_{1}}dh_{i_1}\big).
\end{split}
\end{equation}

In order to show that $\widehat f_1(E)$ is a well-defined element of $\check H^1({\mathcal U},d\Omega^0)$, which is  independent of the choice of the system of the transition functions, it suffices to to prove that, in addition to formulas (\ref{h0}) and (\ref{hcobound}),
\begin{equation}\label{sd}
s_{i_1}\in\Gamma(U_{i_1},d\Omega^0)\,\, {\rm  for }\,\,i_1\in I.
\end{equation}
By formula (\ref{si1}), formula (\ref{sd}) holds in the same way as (\ref{dg}).

The remaining of the proof is similar to the proof of Theorem \ref{cfd} which we omit here to avoid repetition.   We complete the proof of Theorem \ref{fla}. \,\,\,\,\,\,$\endpf$
\medskip

\noindent{\bf Proof of Corollary \ref{comqf} :} Since $X$ is a compact complex manifold,  each open cover of $X$ has a finite good cover as its refinement (see Appendix I in \cite{F} for instance). In particular, there is a finite good cover $\mathcal U$ of $X$.

Recall the following natural homomorpshims,
\begin{equation}
\check H^1({\mathcal U},d\Omega^0)\xrightarrow{L_{\mathcal U}}\check H^1(X,d\Omega^0)\xrightarrow{L} H^1(X,d\Omega^0),
\end{equation}
Then $\check H^1({\mathcal U},d\Omega^0)\cong H^1(X,d\Omega^0)$ (see \S 3.3 of \cite{F}). The remaining of the proof is similar to that of Corollary \ref{comt}, which we omit it to avoid repetition. \,\,\,\,\,\,$\endpf$
\medskip                                                                                                                                                                          

Next, we will show that a holomorphic line bundle is $\mathcal Q$-flat if and only if its first $T$ invariant is trivial in $H^1(X,d\Omega^0)$. 
\medskip

\noindent{\bf Proof of Theorem \ref{flat} :} 
If $E$ is $\mathcal Q$-flat, then for a certain positive integer $m$ the line bundle $mE$ has a system of transition functions which are all constant. Formula (\ref{t12}) yields that $f_1^r(mE)=0$. Since $H^1(X,d\Omega^0)$ is a vector space over $\mathbb C$, $f_1^r(E)=\frac{1}{m}\cdot f_1^r(mE)=0.$

Next assume $f_1^r(E)=0$ in $H^1(X,d\Omega^0)$.  Recall the following short exact sequence of $\mathbb C$-sheaves (see \cite{HA} or \cite{F}),
\begin{equation}\label{COME}
0\rightarrow\mathbb C\rightarrow\Omega^0\rightarrow d\Omega^0\rightarrow 0\,.
\end{equation}
Combined with the exponential sheaf sequence $0\rightarrow\mathbb Z\xrightarrow{2\pi\sqrt{-1}}\mathcal O_X\xrightarrow{\exp} \mathcal O_X^*\rightarrow 0$ (see \cite{GH}), we have the following commutative diagram of sheaves
\begin{equation}\label{cd}
\begin{array}{ccccccccc}
0&\rightarrow&\mathbb Z &\xrightarrow{2\pi\sqrt{-1}} &\mathcal O_X&\xrightarrow{\exp} 
&\mathcal O_X^*&\xrightarrow{}&0\\
&&\Big\downarrow\rlap{$\scriptstyle j$}&&\Big\downarrow\rlap{$\scriptstyle Id$}&&\Big\downarrow\rlap{$\scriptstyle G$}&\\
0&\rightarrow&\mathbb C&\rightarrow & \Omega^0&\rightarrow& d\Omega^0&\xrightarrow{} &0\\
\end{array}.
\end{equation}
Here $j(n)=2\pi\sqrt{-1}\cdot n$, $n\in\mathbb Z$; by convention, $\mathcal O_X=\Omega^0$ and $Id$ is the identity map; $G(g)=d\log(g)$ for each germ $g$ of sheaf $\mathcal O_X^*$.

Taking the corresponding long exact sequences of cohomology groups, we derive that
\begin{equation}
\rightarrow H^0(X,\Omega^0)\rightarrow H^0(X,d\Omega^0)\rightarrow H^1(X,\mathbb C)\rightarrow H^1(X,\Omega^0)\rightarrow H^1(X,d\Omega^0)\rightarrow H^2(X,\mathbb C),
\end{equation}
\begin{equation}\label{iH1d}
\rightarrow H^0(X,\mathcal O_X )\rightarrow H^0(X,\mathcal O_X^*)\rightarrow H^1(X,\mathbb Z )\rightarrow H^1(X,\mathcal O_X)\rightarrow H^1(X,\mathcal O_X^*)\rightarrow H^2(X,\mathbb Z ).
\end{equation}
Notice that the group structre of the sheaf $\mathcal O_X^*$ is given by multiplication instead of addition, and  $H^*(X,\mathcal O_X^*)$ should be understood as \v Cech cohomology groups (see \cite{GH}).  On the other hand, since the other groups in (\ref{iH1d}) are isomorphic to their corresponding \v Cech cohomology groups, we will view all groups in (\ref{iH1d}) as \v Cech cohomology groups in the following.

Let   $\mathcal U:=\{U_i\}_{i=1}^l$ be a good cover of $X$. We have that (see \S 3.3 of \cite{F} for instance),
\begin{equation}
\begin{split}
H^p(X,&\Omega^0)\cong\check H^p({\mathcal U},\Omega^0),\,\,\,\,\,H^p(X,d\Omega^0)\cong\check H^p({\mathcal U},d\Omega^0),\,\,\,\,\,H^p(X,\mathbb C)\cong\check H^p({\mathcal U},\mathbb C),\\
H^p(X,&\mathcal O_X^*)\cong\check H^p({\mathcal U},\mathcal O_X^*),\,\,\,\,\,H^p(X,\mathbb Z)\cong\check H^p({\mathcal U},\mathbb Z)\,\,{\rm for}\,\,p\geq 0.
\end{split}
\end{equation}
Moreover,  we can derive the following commutative diagram from (\ref{cd}),
\begin{equation}\label{rmapstodo}
\begin{array}{ccccccc}
\check H^1({\mathcal U},\mathbb Z)&\rightarrow &\check H^1({\mathcal U},\mathcal O_X)&\xrightarrow{h}
&\check H^1({\mathcal U},\mathcal O_X^*)&\xrightarrow{c_1}&\check H^2({\mathcal U},\mathbb Z)\\
\Big\downarrow\rlap{$\scriptstyle j_{*} $ }&&\Big\downarrow\rlap{$\scriptstyle Id$}&&\Big\downarrow\rlap{$\scriptstyle \widehat f_{1,\mathcal U}^r$}&&\Big\downarrow\rlap{$\scriptstyle j_{*}$}\\
\check H^1({\mathcal U},\mathbb C)&\xrightarrow{B} & \check H^1({\mathcal U},\Omega^0)&\xrightarrow{D}&\check H^1({\mathcal U},d\Omega^0)&\xrightarrow{\delta^1} &\check H^2({\mathcal U},\mathbb C)\\
\end{array},
\end{equation}
Here  $c_1$ is the first Chern class map with the image in  $\check H^2({\mathcal U},\mathbb Z)$ and $Id$ is the identity map; since each element of  $\check H^1({\mathcal U},\mathcal O_X^*)$ can be one to one identified with a holomoprhic line bundle over $X$ up to biholomorphisms, the map $\widehat f_{1,\mathcal U}^r$ is induced by the refined first $T$ invariant map.

By a slight abuse of notation, we still denote by $E$ the element in $\check H^1({\mathcal U},\mathcal O_X^*)$ corresponding to the holomoprhic line bundle $E$.  Since  $\check H^1({\mathcal U},d\Omega^0)\cong\check H^1(X,d\Omega^0)\cong H^1(X,d\Omega^0)$, $\widehat f_{1,\mathcal U}^r(E)=0$ in $\check H^1({\mathcal U},d\Omega^0)$; hence, $\delta^1(f_{1,\mathcal U}^r(E))=0$ in $\check H^2({\mathcal U},\mathbb C)$. Then, $c_1(E)$ is a torsion element in $\check H^2({\mathcal U},\mathbb Z)\cong H^2(X,\mathbb Z)$.  Take a positive integer $m$ such that $c_1(mE)=0$. Since the horizontal lines of (\ref{rmapstodo}) are exact, there is an element  $e\in\check H^1({\mathcal U},\mathcal O_X)$ such that $h(e)=mE$. Since $D(Id(e))=\widehat f_{1,\mathcal U}^r(mE)=0$, there is an element $\xi\in\check H^1(X,\mathbb C)$ such that $B(\xi)=Id(e)$. Note that $\xi$ is a \v Cech $1$-cocycle with coefficients in $\mathbb C$, and hence it can be represented as follows
\begin{equation}
\xi=\bigoplus\limits_{1\leq j_1<j_2\leq l}\xi_{j_1j_2}\in\bigoplus\limits_{1\leq j_1<j_2\leq l}\Gamma(U_{j_1j_2},\mathbb C).
\end{equation}
By chasing diagram (\ref{rmapstodo}) in the reverse order, we have a  system of transition functions associated with a certain trivialization of $mE$ with respect to $\mathcal U$ as follows:
\begin{equation}
g_{j_1j_2}=\exp(\xi_{j_1j_2})\in\Gamma(U_{j_1j_2},\mathbb C)\,\,{\rm for}\,\,1\leq j_1<j_2\leq l.
\end{equation}

We complete the proof of Lemma \ref{flat}.
\,\,\,\,\,$\endpf$
\medskip

To indicate the geometric meaning of the refined first $T$ invariants, we introduce the following notation whose consistency is ensured by Theorem \ref{flat}.
\begin{definition}\label{dfla}
	Let $X$ be a compact complex manifold and $E$ be a holomorphic line bundle over $X$. Let $f_1^r(E)\in H^1(X,d\Omega^0)$ be the refined first $T$ invariant  defined as  Theorem \ref{fla}.  Then we call $\frac{1}{2\pi\sqrt{-1}}f_1^r(E)$ the $\mathcal Q$-flat class of $E$.	Moreover, the map associating each line bundle with its $\mathcal Q$-flat class is called the $\mathcal Q$-flat class map and is denoted by $F$. 
\end{definition}

\begin{remark}
	It is clear that $F$ is a homomorphism from $H^1(X,\mathcal O_X^*)$ to $H^1(X,d\Omega^0)$  such that $F(\mathcal O_X)=0$ and $F(E_1\otimes E_2)=F(E_1)+F(E_2)$ for line bundles $E_1$ and $E_2$.
\end{remark}

Denote by $i_2$ the natural homormohpism from $H^1(X,\Omega^1)$ to $H^2(X,\mathbb C)$. We then have the following proposition. 
\begin{proposition}\label{chern} Let $X$ be a compact complex manifold. The first Chern class maps factor through the $\mathcal Q$-flat class map $F$ as follows:
	\begin{equation}\label{fact}
		H^1(X,\mathcal O^*)\xrightarrow{F} H^1(X,d\Omega^0)\xrightarrow {j_1} H^1(X,\Omega^1)\xrightarrow{i_2} H^2(X,\mathbb C)\,\,.
	\end{equation}
	That is, for each holomorphic line bundle $E$ of $X$,  $(j_1\circ F)(E)$ is the first Chern class of $E$ in the Dolbeault cohomology group $H^1(X,\Omega^1)$; $(i_2\circ j_1\circ F)(E)$ is the first Chern class of $E$ in the De Rham cohomology group $H^2(X,\mathbb C)$.
\end{proposition}

\noindent{\bf Proof of Proposition \ref{chern} :} 
Let   $\mathcal U:=\{U_i\}_{i=1}^l$ be a good cover of $X$ and view all the above groups as \v Cech cohomology groups $H^*({\mathcal U},\cdot)$. Then, Proposition \ref{chern} is an easy consequence of Theorem \ref{compare} and Proposition $1$ in \S 1.1 of \cite{GH}. \,\,\,\,\,$\endpf$
\medskip

\subsection{Strict refinement criterion}
It is well known that when $X$ is a non-K\"ahler manifold the homomorphism $i_2$ in (\ref{fact}) may not be injective, so that the first Chern class in the Dolbeault cohomology may be finer than the first Chern class in the De Rham cohomology (see Example \ref{calabi}).

In the following, we will give an criterion for when the refined $T$ invariant (or equivalently, the $\mathcal Q$-flat class) is strictly finer than the first $T$ invariant (or equivalently, the first Chern class in the Dolbeault cohomology).  

First, let us recall some notation of Fr\"olicher spectral sequence following \cite{Fr}. For $p,q\geq 0$, denote by ${}^{p}T^{n}$ the vector space consisting of differential forms locally given by 
\begin{equation}
\sum_{r\geq p;\,r+s=n} \phi_{\alpha_1\cdots\alpha_r\beta_1\cdots\beta_s} dz_{\alpha_1}\wedge dz_{\alpha_2}\wedge\cdots dz_{\alpha_r}\wedge d\bar z_{\beta_1}\wedge d\bar z_{\beta_2}\wedge\cdots d\bar z_{\beta_s},
\end{equation}
where $z_{v}^{\prime}\,s$ are local holomorphic coordinates. For  $r\geq 0$, define  ${}^{p}T^{n}_r$, ${}^{p}B^{n}_r$ and $E^{p,q}_r$ by
\begin{equation}
{}^{p}T^{n}_r:=\{t\big|t\in {}^{p}T^{n};\,dt\in {}^{p+r}T^{n+1}\},
\end{equation}
\begin{equation}
{}^{p}B^{n}_r:=\{t\big|t\in {}^{p}T^{n};\,t=ds;\,s\in {}^{p-r}T^{n-1}\},
\end{equation}
\begin{equation}
E^{p,q}_r:=\frac{{}^{p}T^{p+q}_r}{{}^{p+1}T^{p+q}_{r-1}+{}^{p}B^{p+q}_{r-1}}\,\,\,\,\,.
\end{equation}
Notice that the exterior differentiation induces the following cochain complex
\begin{equation}
\cdots\rightarrow E^{p-r,q+r-1}_r\xrightarrow{d_r} E^{p,q}_r\xrightarrow{d_r} E^{p+r,q-r+1}_r\rightarrow\cdots.
\end{equation}
Moreover, if we denote by $ K^{p,q}_r\subset  E^{p,q}_r$ the $d_r$-kernel, then
\begin{equation}\label{r=+1}
E^{p,q}_{r+1}=\frac{K^{p,q}_r}{d_r(E^{p+r,q-r-1}_r)}\,\,.
\end{equation}
We also collect some complexes for small $p,q,r$ as follows:
\begin{equation}
\begin{split}
&E_1^{0,0}\xrightarrow{d_1} E_{1}^{1,0}\xrightarrow{d_1} E_{1}^{2,0}\xrightarrow{d_1} E_{1}^{3,0};\\
&0\xrightarrow{d_2} E_2^{1,0}\xrightarrow{d_2} 0,\,\,\,\,0\xrightarrow{d_3} E_3^{1,0}\xrightarrow{d_3} 0,\,\,\,0\xrightarrow{d_4}E_4^{1,0}\xrightarrow{d_4} 0,\cdots;\\
&0\xrightarrow{d_1} E_1^{0,1}\xrightarrow{d_1} E_{1}^{1,1}\xrightarrow{d_1} E_{1}^{2,1}\xrightarrow{d_1} E_{1}^{3,1};\\
&0\xrightarrow{d_2} E_2^{0,1}\xrightarrow{d_2}E_{2}^{2,0}\xrightarrow{d_2} 0,\,\,\,
0\xrightarrow{d_2}E_3^{0,1}\xrightarrow{d_2}  0,\,\,\,\, 0\xrightarrow{d_2} E_4^{0,1}\xrightarrow{d_2} 0,\,\,\,0\xrightarrow{d_2} E_5^{0,1}\xrightarrow{d_2} 0,\cdots.
\end{split}
\end{equation}

Under the above notation, we state the criterion as follows.

\begin{theorem}\label{coincides} Let $X$ be a compact complex manifold. The first Chern classes factor through the $\mathcal Q$-flat class map $F$ as follows:
	\begin{equation}
	H^1(X,\mathcal O^*)\xrightarrow{F} H^1(X,d\Omega^0)\xrightarrow {j_1} H^1(X,\Omega^1)\xrightarrow{i_2} H^2(X,\mathbb C)\,\,.
	\end{equation}
	Denote by  $j_1\big|_{{\rm Im}(F)}$  the restriction of $j_1$ to the image of $F$. Then $j_1\big|_{{\rm Im}(F)}$ is injective if and only if
	\begin{equation}\label{zero}
	E_2^{0,1}=E_3^{0,1}\,\,.
	\end{equation}
\end{theorem}

\noindent{\bf {Proof of Theorem \ref{coincides} :}}  We first prove that $j_1\big|_{{\rm Im}(F)}$ is injective if and only if $j_1$ is. If $j_1$ is injective, then it is cleat that $j_1\big|_{{\rm Im}(F)}$ is injective. Now suppose $j_1\big|_{{\rm Im}(F)}$ is injective.  Recalling diagram (\ref{rmapstodo}), we have the following commutative diagram:
\begin{equation}\label{rmap}
\begin{tikzcd}
	&H^1(X,\mathcal O_X)\ar[d,"Id"]\ar[r,"h"]&H^1(X,\mathcal O_X^*) \ar[d,"F"] &\\
   &H^1(X,\Omega^0)\ar[r,"D"]&H^1(X,d\Omega^0)\ar[d,"j_1"]\ar[r,"i_2\circ j_1"]&H^2(X,\mathbb C)\\
	&&H^1(X,\Omega^1)\ar[ur,"i_2"]&.
\end{tikzcd}
\end{equation}
Assume that there is an element $\xi\in H^1(X,d\Omega^0)$ such that $j_1(\xi)=0\in H^1(X,\Omega^1)$. Then, $(i_2\circ j_1)(\xi)=0\in H^2(X,\mathbb C)$ and hence there is an element $\tau\in H^1(X,\Omega^0)$ such that $D(\tau)=\xi$. Thus, $(h\circ(Id)^{-1})(\tau)$ gives a line bundle $E\in H^1(X,\mathcal O_X^*)$. Since the diagram is commutative, we have that $F(E)=\xi$ and $(j_1\circ F)(E)=j_1(\xi)=0$. Because $j_1\big|_{{\rm Im}(F)}$ is injective, we conclude that $\xi=F(E)=0$. Therefore, $j_1$ is injective.

Next we will apply the Fr\"olicher spectral sequence to prove that $j_1$ is injective if and only if (\ref{zero}) holds.  
Combining diagrams $(\ref{rmapstodo})$ and (\ref{rmap}), we have
\begin{equation}\label{ddd}
\begin{tikzcd}
	0\ar[r] &H^0(X,d\Omega^0)\ar[r]&H^1(X,\mathbb C)\ar[r]& H^1(X,\Omega^0)\ar[r,"D"]\ar[dr,"j_1\circ D"]&H^1(X,d\Omega^0)\ar[d,"j_1"]\ar[r,"i_2\circ j_1"]&H^2(X,\mathbb C)\\
	&&&&H^1(X,\Omega^1)\ar[ur,"i_2"]&.
\end{tikzcd}
\end{equation} 
It is easy to see that
\begin{equation}
\ker j_1\subset \Ima D,
\end{equation}
and hence $\ker j_1=0$ if and only if
\begin{equation}
\dim \Ima (j_1\circ D)=\dim \Ima D.
\end{equation}

By formula  $(\ref{r=+1})$, we have
\begin{equation}
E_{2}^{0,1}=\ker(E_1^{0,1}\xrightarrow{d_1} E_{1}^{1,1}),\,\,\,\,\,\,\,
E_{3}^{0,1}=E_{4}^{0,1}=\cdots=\ker (E_2^{0,1}\xrightarrow{d_2} E_{2}^{2,0}),
\end{equation}
\begin{equation}\label{10}
E_{2}^{1,0}=E_{3}^{1,0}=\cdots=\frac{\ker (E_1^{1,0}\xrightarrow{d_1} E_{1}^{2,0})}{{\rm Im} (E_1^{0,0}\xrightarrow{d_1} E_{1}^{1,0})},
\end{equation}
and hence (see Theorem $3$ in \cite{Fr})
\begin{equation}\label{dd}
\dim H^1(X,\mathbb C)=\dim E_{3}^{0,1}+\dim E_{3}^{1,0}.
\end{equation}
Since (see Lemma $1$ in \cite{Fr}) 
\begin{equation}
E_{1}^{p,q}\cong H^q(X,\Omega^p)\cong H^{p,q}(X),
\end{equation}
it is easy to see that $H^0(X,d\Omega^0)=E_{2}^{1,0}=E_{3}^{1,0}$.

We will introduce the double complex interpretation of $E^{p,q}_r$. Recall the following double complex  in \S 14 of \cite{BT}:
\begin{equation}
\begin{tikzcd}
0&&&\\
Z^{0,1}\ar[u,"\bar\partial"]\ar[r,"\partial"] &Z^{1,1}&\\
&Z^{1,0}\ar[u,"\bar\partial"]\ar[r,"\partial"] &Z^{2,0}&\\
&&0\ar[u,"\bar\partial"]
\end{tikzcd},
\end{equation}
where $Z^{p,q}$ are differential forms on $X$ of bidegree $(p,q)$ ($p$-th wedge of $dz^{\nu}$ and $q$-th wedge of $d\bar z^{\mu}$).  Notice that an element in $E_1^{0,1}$ can be represented by an element $b\in Z^{0,1}$ such that $\bar\partial b=0$; an element of $E_2^{0,1}$ can be represented by an element $b\in Z^{0,1}$ such that $\bar\partial b=0$ and there is an element $c\in  Z^{1,0}$ with $\partial b=\bar\partial c$; an element of $E_3^{0,1}$ can be represented by an element $b\in Z^{0,1}$ such that $\bar\partial b=0$, and there is an element $c\in  Z^{1,0}$ with $\partial b=\bar\partial c$ and $\partial c=0$.  Therefore, $H^1(X,\Omega^0)$ has a filtration
\begin{equation}
H^1(X,\Omega^0)=E_{1}^{0,1}\supset E_{2}^{0,1} \supset E_{3}^{0,1};
\end{equation}
by formula (\ref{dd}), $\Ima (D)\cong E_{1}^{0,1}/E_{3}^{0,1}$.

Now in order to Theorem \ref{coincides}, it remains to show that $\Ima (j_1\circ D)\cong E_{1}^{0,1}/E_{2}^{0,1}$. Notice that the map $j_1\circ D$ is induced by the following differential map between sheaves
\begin{equation}
d:\Omega^0\rightarrow \Omega^1,
\end{equation} 
and $H^q(X,\Omega^p)\cong H^{p,q}(X)$.
We derive that
\begin{equation}
(j_1\circ D)(b)=\partial(b),
\end{equation} 
for each $b\in Z^{0,1}$ such that $\bar\partial b=0$. Noticing that $\partial b=0\in H^1(X,\Omega^1)$ if and only if there is an element $c\in Z^{1,0}$ such that $\bar\partial c=\partial b=0$, we conclude that $\Ima (j_1\circ D)\cong E_{1}^{0,1}/E_{2}^{0,1}$.

Therefore, we complete the proof of Theorem \ref{coincides}.
\,\,\,\,\,\,\,\,$\endpf$
\medskip

\noindent{\bf Proof of Theorem \ref{criterion} :} Theorem \ref{criterion} follows from Theorem \ref{coincides}.
\,\,\,\,\,\,\,\,$\endpf$
\medskip

\begin{remark}
	If $X$ has Property $(H)$ as defined in \cite{F}, then the restriction of $j_1$ to the image of $F$ and the restriction of $i_2\circ j_1$  to the image of $F$ are both injective. In particular, this is always true for K\"ahler manifolds and manifolds of Fujiki class $\mathcal C$.
\end{remark}

\begin{example}[See \cite{Bor}]\label{calabi}
	Let $X$ be a Calabi-Eckmann manifold which is diffeomorphic to the product manifold $S^{2u+1}\times S^{2v+1}$ ($u,v\geq 1$). Then, $\dim H^1(X,d\Omega^0)=\dim H^1(X,\Omega^1)=1$ and $\dim H^2(X,\mathbb C)=0$ . Therefore, the first Chern class in the Dolbeault cohomology is strictly finer than the first Chern class in the De Rham cohomology for a certain line bundle. However, the refined first $T$ invariant in $H^1(X,d\Omega^0)$  is isomorphic to its image in $H^1(X,\Omega^1)$ for every holomoprhic vector bundle of $X$.
\end{example}
\begin{example}[See \cite{RB}]\label{exam}
Consider the real, nilpotent subgroup of $GL(6,\mathbb C)$
\begin{equation}\label{g2}
G_2:=\left\{\left(
\begin{matrix}
1&0&0&0&-\bar y_1&w_1\\
&1&\bar z_1&- x_1&0&w_2\\
&&1&0&0&y_1\\
&&&1&0&y_2\\
\vdots&&&&1&z_1\\
0&\cdots&&&&1\\
\end{matrix}\right):x_1,y_1,y_2,z_1,w_1,w_2\in\mathbb C\right\}.
\end{equation}
Regarding $x_1,y_1,y_2,z_1,w_1,w_2$ as complex coordinates we can identify $G_2$ with $\mathbb C^6$ and then multiplication on the left with a fixed element is holomorphic. Taking the
quotient with respect to the discrete subgroup $\Gamma:=G_2\cap GL(6,\mathbb Z[i])$ acting on the left yields a (compact) nilmanifold with left-invariant complex structure
$X_2:=\Gamma\backslash G_2$.

If we call the matrix in (\ref{g2}) $A$, then the space of left-invariant differential
$(1,0)$-forms is spanned by the components of $A^{-1}dA$ which yields the following:
\begin{equation}
U:= \{dx_1, dy_1,dy_2,dz_1, \omega_1,\omega_2\}
\end{equation}
where 
\begin{equation}
\begin{split}
&\omega_1=dw_1-\bar y_1dz_1,\\
&\omega_2=dw_2-\bar z_{1}dy_1+x_1dy_2\,\,.
\end{split}
\end{equation}
The differential form $\bar\omega_1$ defines a class $[\bar\omega_1]_2$ in $E^{0,1}_2$ and
\begin{equation}
d_2([\bar\omega_1]_2)=[dx_1\wedge dy_2]_2\neq 0\,\,{\rm in}\,\,E_2^{2,0}.
\end{equation}
Therefore, $E_2^{0,1}\neq E_3^{0,1}$. By Theorem \ref{coincides}, we have a line bundle whose $\mathcal Q$-flat class is not zero but its first Chern class in the Dolbeault cohomology is zero.

\end{example}
\begin{remark}[\cite{A}] Let $X$ be a compact complex manifold and $E$ be a vector  bundle over $X$. Then, $E$ has a holomorphic connection if and only if the Atiyah class $a(E)$ of $E$ is trivial, where
\begin{equation}
a(E)\in H^1(X,\Omega^1\otimes End (E)).
\end{equation}
Notice that when $E$ is a line bundle $a(E)$ is the first Chern class in the Dolbeault cohomology. Therefore, when $E_2^{0,1}=E_3^{0,1}$ in the Fr\"olicher spectral sequence of $X$ a line bundle has a holomorphic connection if and only if it is flat up to a  positive  multiple.
\end{remark}

\section{Invariants of holomorphic vector bundles II : refinement of the higher $T$ invariants in $H^k(X,d\Omega^{k-1})$, $k\geq 2$.}
In this section, we will generalize the discussion in \S 3.3 to the higher $T$ invariants for holomorphic vector bundles. Notice that when $k=1$ the refinement is already given in \S 3. However,  in general cases when $k\geq 2$, the $k$-th $T$ invariant can not be refined to be an element of $H^k(X,d\Omega^{k-1})$, for the coefficients of the cocycle defined by formula $(\ref{t11})$ are not $d$-closed.  We will show that when a vector bundle has a full flag structure there is a refinement of the $T$ invariants when $k\geq 2$.
\subsection{The ring structure of $\Phi_X:=\oplus_{i=0}^{\infty}H^i(X,d\Omega^{i-1})$}
In this subsection, we will introduce the cohomological ring $\Phi_X$ where the refined $T$ invariants live.

Let $X$ be a compact complex manifold and let $\Phi_X$ be the vector space defined by
	\begin{equation}
	\Phi_X:=\oplus_{i=0}^{\infty}H^i(X,d\Omega^{i-1}),
	\end{equation}  
where $H^0(X,d\Omega^{-1})=\mathbb C$ by convention.
Recall the following short exact sequences
\begin{equation}\label{short}
\begin{split}
&0\rightarrow \mathbb C\rightarrow \Omega^{0}\rightarrow d\Omega^{0}\rightarrow 0\\
&0\rightarrow d\Omega^{0}\rightarrow \Omega^{1}\rightarrow d\Omega^{1}\rightarrow 0\\
&\,\,\,\,\,\,\,\,\,\,\,\,\,\,\,\,\,\,\,\,\,\,\,\,\,\,\,\,\,\,\,\,\,\cdots\\
&0\rightarrow d\Omega^{q}\rightarrow \Omega^{q+1}\rightarrow d\Omega^{q+1}\rightarrow 0\\
&\,\,\,\,\,\,\,\,\,\,\,\,\,\,\,\,\,\,\,\,\,\,\,\,\,\,\,\,\,\,\,\,\cdots;
\end{split}
\end{equation}
and the associated long exact sequence
\begin{equation}\label{long}
\begin{split}
&\rightarrow H^k(X,\mathbb C)\rightarrow H^k(X,\Omega^0)\rightarrow H^k(X,d\Omega^0)\rightarrow H^{k+1}(X,\mathbb C)\rightarrow,\\
&\rightarrow H^k(X,\Omega^1)\rightarrow H^k(X,d\Omega^1)\rightarrow H^{k+1}(X,d\Omega^0)\rightarrow  H^{k+1}(X,\Omega^1)\rightarrow,\\
&\,\,\,\,\,\,\,\,\,\,\,\,\,\,\,\,\,\,\,\,\,\,\,\,\,\,\,\,\,\,\,\,\cdots\cdots\cdots\\
&\rightarrow H^k(X,\Omega^{q+1})\rightarrow H^k(X,d\Omega^{q+1})\rightarrow H^{k+1}(X,d\Omega^q)\rightarrow  H^{k+1}(X,\Omega^{q+1})\rightarrow,\\
&\,\,\,\,\,\,\,\,\,\,\,\,\,\,\,\,\,\,\,\,\,\,\,\,\,\,\,\,\,\,\,\,\cdots\cdots\cdots.\\
\end{split}
\end{equation}

Then, we can prove that
\begin{lemma}\label{van}
	Denote the complex dimension of $X$ by $N$. Then,
	\begin{equation}
	H^{p}(X,d\Omega^q)=0\,\, {\rm for}\,\,p\geq N+1\,\, {\rm and}\,\,q\geq 0.
	\end{equation}
\end{lemma}
\noindent{\bf Proof of Lemma \ref{van} :} Notice that $H^{p}(X,\mathbb C)=H^{p}(X,\Omega^q)=0$ for $p\geq N+1$ and $q\geq 0$. Then, the lemma follows from the long exact sequences (\ref{long}).\,\,\,\,$\endpf$
\medskip

Hence, we have 
\begin{equation}
\Phi_X=\oplus_{i=0}^{N}H^i(X,d\Omega^{i-1}).
\end{equation}

Next we will give a ring structure for the vector space $\Phi_X$ as follows. 
Let $\mathcal U:=\{U_i\}_{i=1}^l$ of $X$ be  a finite, good open cover in the sense that each $U_i$ is Stein and each  nonempty intersection $U_{i_1\cdots i_p}$ is contractible.  Then, 
\begin{equation}
H^{q+1}(X,d\Omega^{q})\cong\check H^{q+1}({\mathcal U},d\Omega^{q})\,\,{\rm for}\,\,q\geq 0,
\end{equation}
where $\check H^{q+1}({\mathcal U},d\Omega^{q})$ is the $(q+1)$-th \v Cech cohomology group of the sheaf $d\Omega^q$ respect to $\mathcal U$. In the following we will identify the sheaf cohomology groups and their corresponding \v Cech cohomology groups.

Let $\xi\in \check H^{r}(\mathcal U,d\Omega^{r-1})$ and $\eta\in \check H^{s}(\mathcal U,d\Omega^{s-1})$ be \v Cech cocycles as follows:
\begin{equation}
\begin{split}
\xi=\bigoplus\limits_{1\leq i_1<\cdots<i_{r+1}\leq l}u_{i_1\cdots i_{r+1}}\in\bigoplus\limits_{1\leq i_1<\cdots<i_{r+1}\leq l}\Gamma(U_{i_1\cdots i_{r+1}},d\Omega^{r-1}),\\
\eta=\bigoplus\limits_{1\leq i_1<\cdots<i_{s+1}\leq l}v_{i_1\cdots i_{s+1}}\in\bigoplus\limits_{1\leq i_1<\cdots<i_{s+1}\leq l}\Gamma(U_{i_1\cdots i_{s+1}},d\Omega^{s-1}),
\end{split}
\end{equation}
 where $0\leq r,s\leq N$.
 We can define the cup product $\xi\cup\eta\in \check H^{r+s}({\mathcal U},d\Omega^{r+s-1})$ by (see \cite{BT})
\begin{equation}
\xi\cup\eta:=\bigoplus\limits_{1\leq i_1<\cdots<i_{r+s+1}\leq l}w_{i_1\cdots i_{r+s+1}}\in\bigoplus\limits_{1\leq i_1<\cdots<i_{r+s+1}\leq l}\Gamma(U_{i_1\cdots i_{r+s+1}},d\Omega^{r+s-1}),
\end{equation}
where 
\begin{equation}
w_{i_1\cdots i_{r+s+1}}=u_{i_1\cdots i_{r+1}}\wedge v_{i_{r+1}\cdots i_{r+s+1}}.
\end{equation}
Since $dw_{i_1\cdots i_{r+s+1}}=0$, 
\begin{equation}
w_{i_1\cdots i_{r+s+1}}\in\Gamma(U_{i_1\cdots i_{r+s+1}},d\Omega^{r+s-1})\,\,\,{\rm for}\,\,1\leq i_1<\cdots<i_{r+s+1}\leq l.
\end{equation}
Moreover, it is easy to verify that $\xi\cup\eta=\eta\cup\xi$ and hence $\Phi_X$ is a commutative ring.
 
By the natural inclusion we have the following proposition.
\begin{proposition}
	Denote by $\Psi_X$ the cohomology ring
	\begin{equation}
	\Psi_X:=\oplus_{i=1}^{N}H^i(X,\Omega^{i}).
	\end{equation}
	Then, there is a ring homomorphism
	\begin{equation}\label{phipsi}
	j:\Phi_X\rightarrow\Psi_X
	\end{equation}
	induced by the natural group homomorphisms
	\begin{equation}
	 j_p:H^p(X,d\Omega^{p-1})\rightarrow H^p(X,\Omega^p)\,\,{\rm for}\,\,p\geq 0.
	\end{equation}
\end{proposition}

Although the homomorhpism (\ref{phipsi}) is neither injective nor surjective in general even if $X$ is a projective manifold, it is true that $j$ is an isomorphism for $X=\mathbb {CP}^n$.
\begin{proposition}\label{cpn}  The ring homomorphism (\ref{phipsi}) is an isomorphism for $X=\mathbb {CP}^n$.
\end{proposition}
\noindent{\bf Proof of Proposition \ref{cpn} :} We will prove by induction on $k=1,2,\cdots$ that the following properties hold:
\begin{equation}\label{induction}
\begin{split}
& (1)\,\, H^{p}(\mathbb{CP}^n,d\Omega^{k-1})=0\,\,\,\,\,\,\,\,\,\,\,\,\,\,\,\,\,\,\,\,\,\,\,\,\,\,\,\,\,\,\,\,\,\,\,\,{\rm for}\,\,\,\,\,\,\,\,0\leq p\leq k-1\,;\\
& (2)\,\, H^{p}(\mathbb{CP}^n,d\Omega^{k-1})\cong H^{p+1}(\mathbb{CP}^n,d\Omega^{k-2})\cong \cdots\cong H^{p+k}(\mathbb{CP}^n,\mathbb C)\,\,\,\,\,\,\,\,\, {\rm for}\,\,\,\,\,\,\,\, p\geq k\,;\\
& (3)\,\,{\rm the\,\, natural\,\, homomorphism } \,\,H^{k}(\mathbb{CP}^n,d\Omega^{k-1})\rightarrow H^{k}(\mathbb{CP}^n,\Omega^k)\,\,{\rm is\,\, an\,\, isomorphism}.
\end{split}
\end{equation}

For $k=1$, we first consider the following short exact seuqence of sheaves (see (\ref{short})):
\begin{equation}
0\rightarrow \mathbb C\rightarrow \Omega^{0}\rightarrow d\Omega^{0}\rightarrow 0,
\end{equation}
and its associated long exact sequence of cohomology groups (see (\ref{long}))
\begin{equation}\label{l0}
\begin{split}
0&\rightarrow H^0(\mathbb{CP}^n,\mathbb C)\rightarrow H^0(\mathbb{CP}^n,\Omega^0)\rightarrow H^0(\mathbb{CP}^n,d\Omega^0)\rightarrow H^{1}(\mathbb{CP}^n,\mathbb C)\rightarrow\\
\rightarrow\cdots&\rightarrow H^p(\mathbb{CP}^n,\mathbb C)\rightarrow H^p(\mathbb{CP}^n,\Omega^0)\rightarrow H^p(\mathbb{CP}^n,d\Omega^0)\rightarrow H^{p+1}(\mathbb{CP}^n,\mathbb C)\rightarrow\,\cdots.\\
\end{split}
\end{equation}
Noticing that $H^0(\mathbb{CP}^n,\mathbb C)=H^0(\mathbb{CP}^n,\Omega^0)$ and $H^p(\mathbb{CP}^n,\Omega^0)=0$ for $p\geq 1$, Properties ($1$) and ($2$) hold.  

Moreover, by diagram chasing in the \v Cech-algebraic De Rham double complex, we have the following commutative diagram
\begin{equation}\label{ft}
\begin{tikzcd}
H^{1}(\mathbb{CP}^n,d\Omega^{0})\ar[r,"\cong"]\ar[dr]&H^{2}(\mathbb{CP}^n,\mathbb C)\\
&H^{1}(\mathbb{CP}^n,\Omega^{1})\ar[u,"\cong"].
\end{tikzcd}
\end{equation}
Therefore, Property ($3$) holds for $k=1$.

Now suppose that Properties ($1$), ($2$) and ($3$) hold for all $k$ such that $1\leq k\leq l$. We are going to prove that Properties ($1$), ($2$) and ($3$) hold for $k=l+1$. Similarly, considering the following long exact sequence
\begin{equation}\label{ll}
\begin{split}
0&\rightarrow H^0(\mathbb{CP}^n,d\Omega^{l-1})\rightarrow H^0(\mathbb{CP}^n,\Omega^l)\rightarrow H^0(\mathbb{CP}^n,d\Omega^l)\rightarrow H^{1}(\mathbb{CP}^n,d\Omega^{l-1})\rightarrow\\
\cdots&\rightarrow H^p(\mathbb{CP}^n,d\Omega^{l-1})\rightarrow H^p(\mathbb{CP}^n,\Omega^l)\rightarrow H^p(\mathbb{CP}^n,d\Omega^l)\rightarrow H^{p+1}(\mathbb{CP}^n,d\Omega^{l-1})\rightarrow\,\cdots.\\
\end{split}
\end{equation}

For an integer $p$ such that $0\leq p\leq l-2$, noticing that $H^{p}(X,\Omega^l)=H^{p+1}(X,\Omega^l)=0$, we have $H^p(\mathbb{CP}^n,d\Omega^l)\cong H^{p+1}(\mathbb{CP}^n,d\Omega^{l-1})=0$. For an integer $p$ such that $p\geq l+1$, noticing that $H^{p}(X,\Omega^l)=H^{p+1}(X,\Omega^l)=0$, we have $H^p(\mathbb{CP}^n,d\Omega^l)\cong H^{p+1}(\mathbb{CP}^n,d\Omega^{l-1})\cong \cdots\cong H^{p+l+1}(\mathbb{CP}^n,\mathbb C)$.  Moreover, we have the following long exact sequences,
\begin{equation}\label{l-1}
\begin{split}
H^{l-1}(\mathbb{CP}^n,\Omega^{l})&\rightarrow H^{l-1}(\mathbb{CP}^n,d\Omega^l)\rightarrow H^{l}(\mathbb{CP}^n,d\Omega^{l-1})\rightarrow H^l(\mathbb{CP}^n,\Omega^{l})\rightarrow\\
&\rightarrow H^l(\mathbb{CP}^n,d\Omega^l)\rightarrow H^{l+1}(\mathbb{CP}^n,d\Omega^{l-1})\,,\\
\end{split}
\end{equation}
and
\begin{equation}\label{l-2}
H^{l-2}(\mathbb{CP}^n,\Omega^{l+1})\rightarrow H^{l-2}(\mathbb{CP}^n,d\Omega^{l+1})\rightarrow H^{l-1}(\mathbb{CP}^n,d\Omega^{l})\rightarrow H^{l-1}(\mathbb{CP}^n,\Omega^{l+1}),
\end{equation}
\begin{equation}\label{l-3}
H^{l-3}(\mathbb{CP}^n,\Omega^{l+2})\rightarrow H^{l-3}(\mathbb{CP}^n,d\Omega^{l+2})\rightarrow H^{l-2}(\mathbb{CP}^n,d\Omega^{l+1})\rightarrow H^{l-2}(\mathbb{CP}^n,\Omega^{l+2}),
\end{equation}
$$\cdots$$
\begin{equation}\label{l-4}
H^{0}(\mathbb{CP}^n,\Omega^{2l-1})\rightarrow H^{0}(\mathbb{CP}^n,d\Omega^{2l-1})\rightarrow H^{1}(\mathbb{CP}^n,d\Omega^{2l-2})\rightarrow H^{1}(\mathbb{CP}^n,\Omega^{2l-1}),
\end{equation}
\begin{equation}\label{l-5}
0\rightarrow H^{0}(\mathbb{CP}^n,d\Omega^{2l-1})\rightarrow H^{0}(\mathbb{CP}^n,\Omega^{2l})\,.
\end{equation}
By a backward induction based on $(\ref{l-2})$,  $(\ref{l-3})$,  $(\ref{l-4})$ and $(\ref{l-5})$, we conclude that
\begin{equation}\label{va}
H^{l-1}(\mathbb{CP}^n,d\Omega^{l})=H^{l-2}(\mathbb{CP}^n,d\Omega^{l+1})=\cdots=H^{0}(\mathbb{CP}^n,d\Omega^{2l-1})=0.
\end{equation} 
Substituting $(\ref{va})$ into $(\ref{l-1})$, we have that
\begin{equation}\label{l-1}
0\rightarrow H^{l}(\mathbb{CP}^n,d\Omega^{l-1})\rightarrow H^l(\mathbb{CP}^n,\Omega^{l})\rightarrow H^l(\mathbb{CP}^n,d\Omega^l)\rightarrow 0\,,
\end{equation}
Since $\dim H^{l}(\mathbb{CP}^n,d\Omega^{l-1})=\dim H^{2l}(\mathbb{CP}^n,\mathbb C)=H^{l}(\mathbb{CP}^n,\Omega^{l})$, we have that $H^l(\mathbb{CP}^n,d\Omega^l)=0$. 

Notice that Properties ($1$) and ($2$) in (\ref{induction}) hold for $k=l+1$ as above.
To prove Property ($3$) $k=l+1$, we recall the following natural homomorphisms 
\begin{equation}
H^{l+1}(\mathbb{CP}^n,d\Omega^{l})\rightarrow H^{l+2}(\mathbb{CP}^n,d\Omega^{l-1})\rightarrow \cdots\rightarrow H^{2l+2}(\mathbb{CP}^n,\mathbb C)
\end{equation}
Moreover, by diagram chasing in the \v Cech-algebraic De Rham double complex, we have the following commutative diagram
\begin{equation}\label{ft}
\begin{tikzcd}
H^{l+1}(\mathbb{CP}^n,d\Omega^{l})\ar[r,"\cong"]\ar[drr] &H^{l+2}(\mathbb{CP}^n,d\Omega^{l-1})\ar[r,"\cong"]&\cdots\ar[r,"\cong"]&H^{2l+2}(\mathbb{CP}^n,\mathbb C)\\
&&H^{l+1}(\mathbb{CP}^n,\Omega^{l+1})\ar[ur,"\cong"]&.
\end{tikzcd}
\end{equation}
Therefore, Property ($3$) holds for $k=l+1$.

We conlcude Proposition \ref{cpn} by Property $(3)$.
\,\,\,\,\,\,$\endpf$
\medskip

\subsection{Vector bundles with a full flag structure}
In the following, we will define a refined $T$ invariant $f_k(E,\mathcal F)$ for each pair $(E,\mathcal F)$ so that $f_k(E,\mathcal F)\in H^k(X,d\Omega^{k-1})$ for $k\geq 1$.
\medskip

\noindent{\bf Proof of Theorem \ref{f} :} The proof is similar to the proof of Theorem \ref{fla}. We first prove that the coefficients of the \v Cech cochains take values in the sheaf $d\Omega^{k-1}$. 

Notice that by assumption $g_{ij}$ is an upper triangular matrix as follows,
\begin{equation}
g_{ij}=\left(\begin{matrix}
g^1_{ij}&*&\cdots&*&*\\
0&g^2_{ij}&\cdots&*&*\\
\cdots&\cdots&\cdots&\cdots&\cdots\\
0&0&\cdots&g^{M-1}_{ij}&*\\
0&0&\cdots&0&g_{ij}^M
\end{matrix}\right)\,\,\,{\rm for}\,\,\, i,j\in I.
\end{equation}
Then,
\begin{equation}
\begin{split}
t_{i_1\cdots i_{k+1}}&=\sum_{j=1}^M\,\,\,\,d\log g^j_{i_{1}i_{2}}\wedge d\log g^j_{i_{2}i_{3}}\wedge\cdots\wedge d\log g^j_{i_{k}i_{(k+1)}}\in\Gamma(U_{i_1\cdots i_{k+1}},d\Omega^{k-1}).
\end{split}
\end{equation}

Suppose that $\{\widetilde g_{ij}\}$ is another system of transition functions of $E$ with respect to $\mathcal F$ and $\mathcal U$. Then, there exists a \v Cech $0$-cochain
\begin{equation}
h:=\bigoplus\limits_{1\leq i_1\leq l}h_{i_1}\in\bigoplus\limits_{1\leq i_1\leq l}\Gamma_{hol} (U_{i_1},GL(M,\mathbb C)),
\end{equation}
such that $\widetilde g_{ij}=h_i^{-1}g_{ij}h_j$ for $i,j\in I$, where $h_i$ is an upper triangular matrix for $i\in I$.  As in the proof of Lemma \ref{invariant}, we can define a \v Cech $(k-1)$-cochain by formula (\ref{s0})
\begin{equation}\label{cob}
h_{k-1}(E, \mathcal F, g,\widetilde g):=\bigoplus\limits_{j_1<\cdots<j_{k}}s_{j_1\cdots j_{k}}\in\bigoplus\limits_{ j_1<\cdots<j_{k}}\Gamma(U_{j_1\cdots j_{k}},\Omega^k),
\end{equation} 
such that $\partial h_{k-1}(E,g,\widetilde g)=f_k(E,\mathcal F, \widetilde g)-f_k(E,\mathcal F, g)$; that is, for any elements $i_1,\cdots,i_{k+1}\in I$, we have
\begin{equation}\label{rcobound1}
\widetilde t_{i_1\cdots i_{k+1}}- t_{i_1\cdots i_{k+1}}=\sum_{j=1}^{k+1}(-1)^{j-1}s_{i_1\cdots\widehat i_j\cdots i_{k+1}}\big|_{U_{i_1\cdots i_{k+1}}}.
\end{equation}

To show that $f_k(E,\mathcal F)$ is a well-defined element of $\check H^k({\mathcal U},d\Omega^{k-1})$  (independent of the choice of the system of the transition functions which preserves filtration $\mathcal F$), it suffices to to prove that, in the formula (\ref{cob}),
\begin{equation}\label{dcl}
s_{j_1\cdots j_{k}}\in\Gamma(U_{j_1\cdots j_{k}},d\Omega^{k-1})\,\,\, {\rm  for }\,\,j_1<\cdots<j_{k}\in I.
\end{equation}
Recalling formula (\ref{s0}) and the fact that all the matrices therein are upper triangular, we can reduce to the case when $E$ is a line bundle which is proved in Theorem \ref{fla}.

The remaining of the proof is similar to the proof of Theorem \ref{cfd} which we omit here to avoid repetition.

We complete the proof of Theorem \ref{f}.\,\,\,\,\,\,$\endpf$
\medskip

\noindent{\bf Proof of Corollary \ref{comf} :}
The proof is similar to the proofs of Corollary \ref{comqf} which we omit here.
\,\,\,\,\,\,$\endpf$
\medskip

Then, we can define the refined Chern character as follows.
\begin{definition}\label{rcc}
		Let $X$ be a compact complex manifold.  Let $\pi:E\rightarrow X$ be a holomorphic vector bundle of rank $M$ over $X$ with a full flag structure $\mathcal F$.   Define the refined Chern character ch$(E,\mathcal F)$ in the ring $\Phi_X=\oplus_{i=0}^{\infty}H^k(X,d\Omega^{k-1})$ as follows:
	\begin{equation}
	{\rm ch}(E,\mathcal F)=\sum_{k=0}^{\infty}\frac{1}{k!\cdot(-2\pi\sqrt{-1})^{2k}}\cdot f_{k}(E,\mathcal F).
	\end{equation}
\end{definition}

We raise the following questions.

\begin{question}
	If the refined $T$ invariants depend on the flag structure? 
\end{question}
\begin{question}
	Can we refine the $T$ invariants for all holomorphic vector bundles?
\end{question}

\section {$T$ invariants for locally free sheaves of schemes over general fields}
In this section, we will discuss the analogues of the $T$ invariants for the locally free sheaves of schemes over general fields. 
\subsection{Locally free sheaves of schemes}
Let $X$ be a scheme over a field $K$. Denote by $\Omega_{X/{\rm Spec\,}K}$ the sheaf of K\"ahler differential of $X$ over Spec $K$. Denote by $\Omega^k$ the $k$-th exterior power sheaf of $\Omega_{X/{\rm Spec\,}K}$ for $k\geq 1$. Let $E$ be a locally free sheaf over $X$ of rank $M$. In the following, we will prove Theorem \ref{sch}.
\medskip

\noindent{\bf Proof of Theorem \ref{sch} :} We denote the characteristic of $K$ by char$(K)$.

The proof is the same as the proof of Theorem \ref{cfd}. The only thing we need to verify is that when we go through the proof of Theorem \ref{cfd} the equalities therein are non-trivial even if char$(K)$ is not zero. 

We can verify that based on the assumption that  $(k+2)!$ is not divisible by char$(K)$. For instance, formula (\ref{scheme}) is well-defined since $(k+1)!$ is not divisible by char$(K)$; formulas (\ref{tij}) (\ref{t1j}) (\ref{t2j}) (\ref{t3j}) hold nontrivially since char$(K)\neq 2$; formula (\ref{tkj}) holds nontrivially since $(k+2)$ is not divisible by char$(K)$.

Therefore, we complete the proof of Theorem \ref{sch}
\,\,\,\,\,$\endpf$
\medskip

Recalling Definition \ref{cdo} of the cohomological Chern character, we would like to ask the following question.
\begin{question} Do we have the following Whitney sum formula for holomorphic vector bundles $E$ and $F$:
\begin{equation}
{\rm ch}_{coh}(E\oplus F)={\rm ch}_{coh}(E)\cdot{\rm ch}_{coh}(F)?
\end{equation}
\end{question}
\begin{remark}
	The above defined cohomological Chern character is different from  the usual Chern character  in the Chow ring; it is in the cohomology ring $\sum_{i=0}^{\infty}H^i(X,\Omega^i)$. 
\end{remark}

Similarly, we can define the refined $T$ invariants in the same way as \S 3.  However, note that in general there is no good cover of $X$ such that each open set is affine and each intersection is contractible. Hence, the \v Cech cohomology of the sheaf $d\Omega^p$ ($p\geq 0$) may be different from its sheaf cohomology. We state the corresponding definition and theorems as follows.

\begin{dt}\label{adholt}
	Let $X$ be a scheme over a field $K$ char$(K)\neq 2,3$.  Let $E$ be a locally free sheaf over $X$. Let $\mathcal U:=\{U_i\}_{i\in I}$ be an open cover of $X$ and $g:=\{g_{ij}\}$ be a system of transition functions of $E$ with respect to $\mathcal U$. 
	For any indices $i_1,i_2\in I$, define $t_{i_1i_{2}}\in\Gamma(U_{i_1i_{2}},\Omega^k)$ by
	\begin{equation}
	\begin{split}
	t_{i_1 i_{2}}:=\sum_{\sigma\in S_{2}}&\frac{sgn(\sigma)}{2!}\cdot tr\big(g^{-1}_{i_{\sigma(1)}i_{\sigma(2)}}dg_{i_{\sigma(1)}i_{\sigma(2)}}\big)=tr\big(g^{-1}_{i_1i_2}dg_{i_1i_2}\big),
	\end{split}
	\end{equation}
	define a \v Cech $1$-cochain $\widehat f^r_{1,\,\mathcal U}(E)$ by
	\begin{equation}
	\widehat f^r_{1,\,\mathcal U}(E):=\bigoplus\limits_{ i_1<i_{2}}t_{i_1i_{2}}\in\bigoplus\limits_{i_1<i_{2}}\Gamma(U_{i_1i_{2}},d\Omega^0),
	\end{equation}
	Then $	\widehat f^r_{1,\,\mathcal U}(E)$ is a \v Cech $1$-cocycle and hence defines an element $f^r_{1,\,\mathcal U}(E)\in\check H^k(\mathcal U,\Omega^k)$; moreover, $f^r_{1,\,\mathcal U}(E)$ is independent of the choice of $g$.  
	Denote by  $f^r_1(E)$ the image of $f^r_{1,\,\mathcal U}(E)$ in the \v Cech cohomology group $\check H^1(X,d\Omega^0)$ under the canonical homomorphism.
	Then $f_1^r(E)$ is independent of the choice of $\mathcal U$. We call  $f_1^r(E)$ the refined first $T$ invariant of $E$.
\end{dt}

Similarly, we can make the definition of full flag structures for locally free sheaves over schemes, and hence have the following theorem.
\begin{dt}\label{drholt}
	Let $X$ be a scheme over a field $K$ char$(K)\neq 2,3$.  Let $E$ be a locally free sheaf over $X$ with a full flag structure $\mathcal F$. Let $\mathcal U:=\{U_i\}_{i\in I}$ be an open cover of $X$ and $g:=\{g_{ij}\}$ be a system of transition functions of $E$ with respect to $\mathcal U$. 
	Then, there is a well-defined \v Cech $k$-cocycle $\widehat f^r_{k,\,\mathcal U}(E,\mathcal F)$,
	\begin{equation}
	\widehat f^r_{k,\,\mathcal U}(E,\mathcal F):=\bigoplus\limits_{i_1<\cdots<i_{k+1}}t_{i_1\cdots i_{k+1}}\in\bigoplus\limits_{i_1<\cdots<i_{k+1}}\Gamma(U_{i_1\cdots i_{k+1}},d\Omega^{k-1}),
	\end{equation}
	\begin{equation}
	\begin{split}
	t_{i_1\cdots i_{k+1}}:=\sum_{\sigma\in S_{k+1}}&\frac{sgn(\sigma)}{(k+1)!}\cdot tr\big(g^{-1}_{i_{\sigma(1)}i_{\sigma(k+1)}}dg_{i_{\sigma(1)}i_{\sigma(k+1)}}g^{-1}_{i_{\sigma(2)}i_{\sigma(k+1)}}dg_{i_{\sigma(2)}i_{\sigma(k+1)}}\cdot\\
	&\cdot g^{-1}_{i_{\sigma(3)}i_{\sigma(k+1)}}dg_{i_{\sigma(3)}i_{\sigma(k+1)}}\cdots g^{-1}_{i_{\sigma(k)}i_{\sigma(k+1)}}dg_{i_{\sigma(k)}i_{\sigma(k+1)}}\big).
	\end{split}
	\end{equation}
	Denote  by $f^r_k(E,\mathcal F)$ the image of $\widehat f^r_{k,\,\mathcal U}(E,\mathcal F)$ in $\check H^k(X,d\Omega^{k-1})$ under the canonical homomorphism.  Then $f_k^r(E,\mathcal F)$ is independent of the choice of $\mathcal U$. We call  $f_k^r(E,\mathcal F)$ the refined $k$-th $T$ invariant of $E$ with respect to $\mathcal F$.
\end{dt}

Since the proofs of the above theorems are similar to the complex manifold case, we omit them here.  We end up with this section by the following questions.
\begin{question} What is the relation between the cohomological Chern character and the Chern character in the Chow ring when $X$ is a projective variety?
\end{question}
\begin{question} Let $X$ be an arithmetic variety (see \cite{GS1}). Can we inteprete the arithmetic Chern character defined by Gillet and Soul\'e in \v Cech cohomology?   
\end{question}

\subsection{Real topological vector bundles}

Since the characteristic classes of real vector bundles over smooth manifold live only in $\mathbb Z/2\mathbb Z$ coefficients, the cocycles constructed as in (\ref{t11}) are trivial in $H^*(X,\mathbb C)$. 

Notice that one can represent  the first Stiefel-Whitney class in \v Cech cohomology by the following well known fact.

\begin{fact}
	Let $X$ be a compact real manifold.  Let $\pi:E\rightarrow X$ be a real vector bundle over $X$. Let $\mathcal U:=\{U_i\}_{i=1}^l$ be a finite good open cover of $X$ (in the sense that each nonempty intersection is contractible) and $g:=\{g_{ij}\}$ be a system of transition functions of $E$ with respect to $\mathcal U$. Then the  first Stiefel-Whitney class  $w_1(E)$ of $E$ is represented by the following \v Cech $1$-cocycle 
	\begin{equation}
	w_1(E):=\bigoplus\limits_{1\leq i_1<i_2\leq l}t_{i_1i_{2}}\in\bigoplus\limits_{1\leq i_1<i_{2}\leq l}\Gamma(U_{i_1 i_{2}},\mathbb Z/2\mathbb Z),
	\end{equation}
	\begin{equation}
	\begin{split}
    t_{i_1i_{2}}:=\tau(\det g_{i_1i_2}),
	\end{split}
	\end{equation}
	where function  $\tau:\mathbb R\backslash \{0\}\longmapsto \{0,1\}$ takes value $0$ for positive numbers and takes value $1$ for negative numbers. 
\end{fact}

A natural question then is as follows. 
\begin{question} Can we represent the Stiefel-Whitney class of a real vector bundle in  \v Cech cohomology? 
\end{question}

\bigskip
\noindent H. Fang, Department of Mathematics, University of Wisconsin-Madison,
Madison, WI 53706, USA.  (hfang35$ @$wisc.edu)

\end{document}